\documentclass[11pt, oneside]{amsart}   	
\usepackage{geometry}                		
\usepackage{graphicx}				
\usepackage{amssymb}
\usepackage{epstopdf}
\usepackage{amsmath}
\usepackage{amsfonts}
\usepackage{fancyhdr} 
\usepackage{mathrsfs}
\usepackage{wrapfig}
\usepackage{bbold}
\usepackage{amsthm} 
\usepackage{enumitem}
\usepackage[toc]{appendix}

\usepackage[colorlinks=true, pdfstartview=FitV, linkcolor=blue, 
        citecolor=blue, urlcolor=blue]{hyperref}

\newtheorem{thm}{Theorem}[]
\newtheorem{lem}{Lemma}[]
\newtheorem{prop}{Proposition}[section]

\newtheorem{cor}[prop]{Corollary}

\newtheorem{alg}{Algorithm}[]

\newenvironment{pf}{\noindent {\bf\em Proof}.\ \ }{\hspace*{\fill}\rule{1.4ex}{1.4ex}\,}

\newcommand{\ddd}{\mathbb{D}}
\newcommand{\real}{\mathbb{R}}
\newcommand{\complex}{\mathbb{C}}

\newcommand{\pro}{\mathfrak{p}\!}

\newcommand{\integer}{\mathbb{Z}}
\newcommand{\one}{\mathbb{1}}

\newcommand{\simp}{\mathscr{S}}

\newcommand{\hilbert}{W^{1,2}}

\newcommand{\reg}{\mathbf{Reg}}
\newcommand{\step}{\mathbf{Step}}

\newcommand{\pwplus}{C^1_{{\rm pw},+}}
\newcommand{\Lop}{L}
\newcommand{\calE}{\mathcal{E}}
\newcommand{\bv}{{\rm BV}}

\newcommand{\yop}{\boldsymbol{\mathcal{D}}_{\! y\,}}
\newcommand{\omop}{\boldsymbol{\mathcal{A}}_{\omega}}
\newcommand{\omopone}{\boldsymbol{\mathcal{A}}_{\omega,1}}
\newcommand{\omoptwo}{\boldsymbol{\mathcal{A}}_{\omega,2}}
\newcommand{\somop}{\boldsymbol{\mathcal{Z}}_{\omega}}
\newcommand{\somopn}{\boldsymbol{\mathcal{Z}}_{n,\omega}}

\newcommand{\pwac}{C_{\rm pw}}
\newcommand{\pwacplus}{C_{{\rm pw},+}}

\DeclareMathOperator{\supp}{supp}

\DeclareMathOperator{\sing}{sing}
\DeclareMathOperator{\regular}{reg}
\DeclareMathOperator{\vol}{Vol}

\DeclareMathOperator{\esssup}{ess\,sup}
\DeclareMathOperator{\spec}{spec}
\DeclareMathOperator{\aut}{{\rm Aut}}
\DeclareMathOperator{\ran}{Ran}
\DeclareMathOperator{\Th}{Th}

\usepackage{scalerel,stackengine}
\stackMath
\newcommand\reallywidecheck[1]{%
\savestack{\tmpbox}{\stretchto{%
  \scaleto{%
    \scalerel*[\widthof{\ensuremath{#1}}]{\negmedspace\kern0pt\bigwedge\kern0pt}%
    {\rule[-\textheight/2]{1ex}{\textheight}}
  }{.5\textheight}%
}{0.5ex}}%
\stackon[1pt]{#1}{\scalebox{-1}{\tmpbox}}%
}
\newcommand{\rwc}{\reallywidecheck}

\title{Scattering on the line via singular approximation}
\author{Peter C.~Gibson}
\address{Dept.~of Mathematics \& Statistics, York University, 4700 Keele St., Toronto, Ontario, Canada, M3J~1P3}
\email{pcgibson@yorku.ca}

\date{February 22, 2022}							

\begin{document}

\begin{abstract}
Motivated by applications to acoustic imaging, the present work establishes a framework to analyze scattering for the one-dimensional wave, Helmholtz, Schr\"odinger and Riccati equations that allows for coefficients which are more singular than can be accommodated by previous theory. 
In place of the standard scattering matrix or the Weyl-Titchmarsh $m$-function, 
the analysis centres on a new object, the generalized reflection coefficient, which maps frequency (or the spectral parameter) to automorphisms of the Poincar\'e disk.  Purely singular versions of the generalized reflection coefficient, which are amenable to direct analysis, serve to approximate the general case. 
Orthogonal polynomials on both the unit circle and unit disk play a key technical role, as does an exotic Riemannian structure on PSL$(2,\mathbb{R})$. 
A central role is also played by the newly-defined harmonic exponential operator, introduced to mediate between impedance (or index of refraction) and the reflection coefficient.  The approach leads to new, explicit formulas and effective algorithms for both forward and inverse scattering. The algorithms may be viewed as nonlinear analogues of the FFT.  
In addition, the scattering relation is shown to be elementary in a precise sense at or below the critical threshold of continuous impedance.  For discontinuous impedance, however, the reflection coefficient ceases to decay at infinity, the classical trace formula breaks down, and the scattering relation is complicated by the emergence of almost periodic structure. 
\end{abstract}

\maketitle

\begin{center}
MSC 34L25, 35P25, 42C05, 65L09, 34B24, 34L40;
\end{center}

\newpage

\tableofcontents

\newpage

\section{Introduction\label{sec-introduction}}

The present paper develops a new approach to forward and inverse scattering on the line based on purely singular approximation. Motivated by applications relating to the wave and Helmholtz equations, the methods developed are applicable also to Schr\"odinger and Riccati equations having singular potentials and coefficients, respectively, beyond the scope of previously established theory.  An extensive literature on one-dimensional scattering, especially concerning the linear Schr\"odinger equation, dates back more than sixty years.  Nevertheless, several new results of a fundamental nature appear here for the first time, including:
\begin{itemize}
\item the relation of scattering to parameterized automorphisms of the Poincar\'e disk;
\item the scattering law for concatenation of media;
\item explicit formulas for both forward and inverse scattering;
\item almost periodic structure encoded by orthogonal polynomials on the unit disk;
\item a singular trace formula valid where the classical trace formula breaks down;
\item forward and inverse scattering algorithms that are non-linear analogues of the FFT.
\end{itemize}
In addition, it is shown that the measurement operator, which in the context of the Helmholtz equation maps acoustic impedance to the reflection coefficient, is elementary in a precise sense if impedance is continuous.  Beyond this critical threshold of continuity, discontinuities in the impedance function engender almost periodic structure in the reflection coefficient, resulting in blowup of the classical trace formula. 
The results lay the foundation for a quantitative analysis of inverse stability, to be taken up in separate work.

The approach to scattering taken in the present paper centres on a new object called the generalized reflection coefficient. It serves as an alternative to the Weyl-Titchmarsh $m$-function which has historically been
used to analyze scattering and spectral theory for the one-dimensional Schr\"odinger equation \cite{DeTr:1979,Si:1999,GeSi:2000}. 
The generalized reflection coefficient is naturally defined in terms of the Helmholtz equation, but its essential distinction from the $m$-function may be sketched in terms of the Schr\"odinger equation 
$$
-y^{\prime\prime}+qy=\lambda y
$$
on an interval $(a,b)$, as follows. The $m$-function arises from the idea, due to Weyl \cite{We:1910}, of varying the spectral parameter $\lambda$ in the upper half plane $\Im\lambda>0$, while fixing a right-hand boundary condition of the form 
$$y^\prime(b)=h y(b),$$ 
where $h\in\mathbb{R}$. (In \cite{We:1910} Weyl was concerned, not with scattering, but with Sturm-Liouville equations on an interval $(a,b)$ and the dichotomy of limit point versus limit circle behaviour as $b\rightarrow\infty$.) 
The corresponding ratio $y^\prime(a;\lambda)/y(a;\lambda)$, which is the $m$-function at $a$, is analytic in the upper half plane and determines the potential $q$, provided $q$ is sufficiently regular (at least locally integrable).  

By contrast, the generalized reflection coefficient is constructed by restricting the spectral parameter $\lambda=\omega^2$ to real values, and letting the above boundary condition at $b$ depend on both $\omega$ and a complex parameter, setting
$$
h=i\omega\frac{1-\xi}{1+\xi}\qquad\left(\xi\in\mathbb{C}, |\xi|<1\right).
$$
For each $\omega\in\mathbb{R}$, the generalized reflection coefficient at $a$, defined as
$$
g^\omega_q(\xi)=\frac{i\omega y(a;\omega,\xi)-y^\prime(a;\omega,\xi)}{i\omega y(a;\omega,\xi)+y^\prime(a;\omega,\xi)},
$$
is thus a function of a complex parameter $\xi$; it is formulated such that setting $\xi=0$ recovers the standard (left-hand) reflection coefficient (i.e., the upper-right entry of the $S$-matrix). 

The generalized reflection coefficient brings additional structure to bear on the scattering problem, since as a function of $\xi$ it turns out to be an automorphism of the Poincar\'e disk model of the hyperbolic plane.  The automorphism group of the hyperbolic plane, PSL$(2,\mathbb{R})$, carries an exotic Riemannian structure (first described by Emmanuele and Salvai \cite{EmSa:2012}) relevant to scattering in that eigenfunctions of its Laplace-Beltrami operator are natural basis functions for the reflection coefficient in the purely singular case \cite{Gi:JFAA2017}.  Working in the context of the Helmholtz equation (related to the Schr\"odinger equation by a change of variables), scattering for the Helmholtz equation with piecewise continuous impedance is analyzed here via approximation by the purely singular case of piecewise constant impedance. This yields forward and inverse scattering results for the Schr\"odinger equation whose potential $q$ is two derivatives less regular than the corresponding impedance.  

The remainder of this introduction is organized as follows.  Section~\ref{sec-background} formulates the motivating acoustic imaging problem, states the main objectives, and compiles some preliminary material---including the basic relations among equations, as well as the passage from initial to boundary conditions. The latter material, although elementary, 
helps to make sense of definitions and technical developments that come later. 
The technical framework for the present paper is set in \S\ref{sec-technical-framework}, where the key functions spaces, definitions and notation are laid out. Lastly, \S\ref{sec-strategy} gives an overview of the main results and explains the basic technical strategy. 

\emph{Acknowledgements.} Thanks to the Mathematics department at PUC-Rio for its generous hospitality in December 2019, when part of this paper was written, and especially to Carlos Tomei for ongoing conversations and helpful comments. Thanks also to Yue Zhao for his enthusiastic discussions about Lax-Milgram methods, and to Fritz Gesztesy for kindly pointing out several  references.

\subsection{Background\label{sec-background}}

\subsubsection{Basic problem\label{sec-basic-problem}}
The basic problem motivated by acoustic imaging of layered media is as follows. 
Consider 
the standard one-dimensional wave equation 
\begin{equation}\label{wave-0}
U_{tt}-c^2U_{yy}=0
\end{equation}
where the strictly positive wave speed $c(y)$ is constant on the left half-line $y\leq y_0$ and variable elsewhere.
A Dirac delta initially travelling from the left half-line toward the right half-line, 
\begin{equation}\label{right-moving-delta}
U(y,t)=\delta\left(\frac{y-y_0}{c(y_0)}-t\right)\qquad(t<0),
\end{equation}
which corresponds to initial conditions
\begin{equation}\label{initial-conditions-0}
U(y,0)=\delta\bigl((y-y_0)/c(y_0)\bigr),\qquad U_t(y,0)=-\delta^\prime\bigl((y-y_0)/c(y_0)\bigr),
\end{equation}
scatters according to the variable structure of $c$ to produce scattering data
\begin{equation}\label{data-0}
d(t)=\left\{\begin{array}{cc}
U(y_0,t)&\mbox{ if }t>0\\
0&\mbox{ if }t\leq 0
\end{array}\right. 
\end{equation}
measured at $y=y_0$.  The problem is to reconstruct the unknown function $c(y)$ for $y>y_0$, given the measured data $d$. Thus one seeks to analyze the map $c\mapsto d$ and its inverse for the most general possible class of functions $c$. 

In practical terms, the wave speed $c$ serves as a proxy for the unknown physical structure of a layered medium, for example biological tissue such as skin or the retina, or sedimentary geological strata, or laminated structures in the built environment.  The measured data $d$ corresponds to acoustic plane wave echoes measured by an ultrasonic transducer or geophone, which of course are recorded for a limited time.  The inverse problem inherent in acoustic imaging is to reconstruct physical structure on as large a spatial interval as possible given a time-limited version of  $d$.  From the geophysical perspective it is natural to consider piecewise continuous wave speeds $c$, since these model a typical stratigraphic structure in which several gradually varying sedimentary formations are fused together at discontinuities---corresponding to geologic non-conformities.   Piecewise continuous wave speeds also provide a reasonably general model for ultrasonic imaging of layered biological tissue or non-destructive testing of layered structures in the built environment.  

\subsubsection{Principal objectives\label{sec-principal-objectives}}

The inverse scattering map $d\mapsto c$ for the wave equation relates closely to inverse scattering for several other equations, most notably the one-dimensional Schr\"odinger equation.  But piecewise continuous wave speeds correspond to Schr\"odinger potentials $q$ having singularities on the order of the derivative of a Dirac pulse, and such potentials do not fit most established theory (see \cite{AlGeHo:1992}). Early work of Faddeev \cite{Fa:1964}, the monumental work of Deift and Trubowitz \cite{DeTr:1979}, and also the inverse spectral theory of Gesztesy and Simon \cite{Si:1999,GeSi:2000} all require that $q$ be locally $L^1$---corresponding to $C^1$ wave speed. Related work by Sylvester and Winebrenner \cite{SyWi:1998}  on a version of the Helmholtz equation requires the index of refraction to be $L^2$, which corresponds to absolutely continuous wave speed. Thus a principal goal of the present work is to develop a scattering theory on the line that accommodates waves speeds of lower regularity than can be analyzed via the aforementioned work.  

A line of investigation concerning piecewise constant, or step function, wave speeds has emerged in the last decade, including work by Albeverio et al.~
 \cite{AlHrMyPe:2013,AlHrMy:2018} as well as the present author \cite{Gi:SIAP2014,Gi:JFAA2017,Gi:JCP2018,Gi:JAT2019}.   But this latter work does not allow any continuous variation in wave speed, and is thus in a sense disjoint from results requiring continuous wave speed or local integrability of $q$. Indeed
 the two lines of research into inverse scattering on the line, associated to continuous wave speeds or locally integrable potentials on one hand, versus piecewise constant wave speeds on the other hand, have never been satisfactorily reconciled.  As suggested by Albeverio et al.~
 \cite[p.4]{AlHrMyPe:2013}, it is of interest to establish a unified theory that simultaneously accommodates both classes, to say nothing of the vast middle ground in between.  The present paper establishes such a unified theory.

Another goal of the present work is to link theory to practical inversion methods applicable to digital acoustic reflection data measured over a finite time interval. 
The approach taken in the present paper is well-adapted to computation, leading directly to two algorithms, corresponding to forward and inverse scattering, that encode the nonlinear scattering relation $d\mapsto c$ in a digital context. 

\subsubsection{Additional related literature\label{sec-related-literature}}

Additional work centred on the $m$-function was brought to the author's attention subsequent to the first draft of the present article. 
A paper of Eckhardt et.~al \cite[\S6]{EcGeNiTe:2013b} extends Simon's local inverse spectral theory for the Schr\"odinger operator \cite{Si:1999} to distributional potentials, drawing on the extensive analysis in \cite{EcGeNiTe:2013a} of minimally regular Sturm-Liouville operators, with boundary conditions appropriate to the $m$-function.  The coefficients allowed in the latter are just as singular as those in the present work. However, the boundary conditions in \cite{EcGeNiTe:2013a} are not those of the generalized reflection coefficient. Importantly, the latter work does not address continuous dependence of the solution to the Sturm-Liouville problem on its $L^\infty$ coefficients---which is needed in the present paper, and established below in \S\ref{sec-mapping}.  

Gesztesy and Sakhnovich have recently applied the $A$-function concept, introduced in \cite{Si:1999}, to Dirac-type systems \cite{GeSa:2020}, thereby generalizing considerably the earlier theory for the Schr\"odinger operator.  Inverse spectral theory and inverse scattering for the Schr\"odinger operator are related (see, for example, \cite{GeSi:1997}), but they have different goals.  Moreover, the technical approach to scattering in the present work, which centres on the generalized reflection coefficient, is essentially different from the aforementioned spectral analysis involving the $m$-function. The generalized reflection coefficient naturally lends itself to singular approximation, leading to fast and accurate computational methods not directly accessible via the $m$-function.

\subsubsection{Helmholtz, Schr\"odinger and Riccati equations\label{sec-schroedinger}}

To obtain a more convenient formulation of the scattering problem for (\ref{wave-0}) one may: transform the spatial variable to travel time distance, $x(y)=\int_0^y1/c(s)\,ds$; replace wave speed with impedance $\zeta=1/c$; and take the Fourier transform with respect to time in the sense of tempered distributions, consistent with the formula
\begin{equation}\label{fourier-transform}
\hat{f}(\omega)=\int_{-\infty}^\infty f(t)e^{i\omega t}\,dt.
\end{equation}
Setting $u(x,\omega)=\widehat{U}(x,\omega)$, and writing $u^\prime$ in place of $u_x$, equation (\ref{wave-0}) then becomes 
\begin{equation}\label{equation-0}
(\zeta u^\prime)^\prime+\omega^2\zeta u=0,
\end{equation}
which is known variously as the Helmholtz equation or the impedance form of the Schr\"odinger equation. The reflection coefficient is defined to be the Fourier transform of the scattering data as defined above in (\ref{data-0}), 
\begin{equation}\label{reflection-coefficient}
R(\omega)=\hat{d}(\omega).
\end{equation}
The present paper treats the scattering map $c\mapsto d$ mainly in terms of its equivalent formulation $\zeta\mapsto R$ in the frequency domain, with equation (\ref{equation-0}) as the principal focus.  
The same function $R$ defined by (\ref{reflection-coefficient}) arises in the context of the Schr\"odinger and Riccati equations related to (\ref{equation-0}) by changes of variables, as follows. 

To relate the wave and Schr\"odinger equations, set 
\begin{equation}\label{alpha}
\alpha=-\zeta^\prime/(2\zeta)=-\textstyle\frac{1}{2}\bigl(\log\zeta\bigr)^\prime,
\end{equation}
introduce a new dependent variable and coefficient,
\begin{equation}\label{y-q}
y=\zeta^{1/2}u\quad\mbox{ and }\quad q=\alpha^2-\alpha^\prime=\bigl(\zeta^{1/2}\bigr)^{\prime\prime}/\zeta^{1/2},
\end{equation}
and express (\ref{equation-0})
in terms of $y$. This yields the standard Schr\"odinger equation
\begin{equation}\label{schroedinger}
-y^{\prime\prime}+qy=\omega^2y. 
\end{equation}
Up to rescaling, the mapping $\alpha\mapsto q=\alpha^2-\alpha^\prime$ in (\ref{y-q}) is the Miura transform, originally important in relating the modified and unmodified KdV equations \cite{Mi:1968}.  Note that if $\zeta$ is absolutely continuous, (\ref{equation-0}) may be re-written in terms of $\alpha$ as
\begin{equation}\label{equation-alpha}
-u^{\prime\prime}+2\alpha u^\prime=\omega^2u.
\end{equation}
The negative sign in (\ref{alpha}) is natural from the point of view of singular approximation, but is not standard. Other authors, including Sylvester and Winebrenner  \cite{SyWi:1998}, who refer to $\alpha$ as the index of refraction, use the opposite sign convention. The present work frequently invokes the association (\ref{alpha}) between $\zeta$ and $\alpha$.

Another well-known transformation of (\ref{equation-0}) results by considering  
\begin{equation}\label{riccati-r}
\mathbf{r}=\frac{u-\frac{1}{i\omega}u^\prime}{u+\frac{1}{i\omega}u^\prime}.
\end{equation}
Straightforward manipulations verify that $\mathbf{r}$ satisfies the Riccati equation
\begin{equation}\label{Riccati}
\mathbf{r}^\prime+2i\omega \mathbf{r}+\alpha(1-\mathbf{r}^2)=0 
\end{equation}
when $u$ satisfies (\ref{equation-0}) or (\ref{equation-alpha}), provided $u^\prime$ and $\alpha$ are regular functions.  

The above transformations ultimately relate one-dimensional scattering problems for the wave, Schr\"odinger and Riccati equations to one another, as follows. The inverse scattering problem for the Schr\"odinger equation as formulated in \cite{Fa:1964} and \cite{DeTr:1979}, namely to recover $q$ given the scattering matrix, or $S$-matrix, 
\begin{equation}\label{scattering-matrix}
S=\begin{pmatrix}T_1&R_2\\ R_1& T_2\end{pmatrix},
\end{equation}
is equivalent to the problem of recovering $q$ given just $R_2$ (or just $R_1$), since each of $R_1$ and $R_2$ alone determines $S$ \cite[\S2, p.~150]{Fa:1964}. In the case where $q$ has the form (\ref{y-q}) and is compactly supported, the reflection coefficient $R_2$ is a multiple of $R=\widehat{d}$ as defined in (\ref{reflection-coefficient}), the precise relationship being
\begin{equation}\label{R2-R}
R_2=e^{2ix_0\omega}R\quad\mbox{ where }\quad x_0=\int_0^{y_0}\frac{1}{c(s)}\,ds
\end{equation}
(see \S\ref{sec-transforming-initial} below). The travel-time coordinate $x_0$ of the data measurement location therefore determines $R$ given $R_2$ (and vice versa). 
If one solves the inverse scattering problem for the wave equation, determining $\zeta$ from $R$, then this also determines $q$ by (\ref{y-q}), thereby solving the inverse scattering problem for the Schr\"odinger equation.  For the Riccati equation, if $\alpha$ is supported on the closure of a bounded interval $X=(x_0,x_1)$, the boundary condition $\mathbf{r}(x_1,\omega)=0$ determines a unique solution on $X$, for which it turns out that $\mathbf{r}(x_0,\omega)=R(\omega)$, the same reflection coefficient defined in (\ref{reflection-coefficient}).  Thus the scattering problem for the Riccati equation of reconstructing $\alpha$ given $\mathbf{r}(x_0,\omega)$ reduces to determination of $\zeta$ given $R$, or in other words, inverse scattering for the wave equation. 

Suppose more generally that $R_2$ is the reflection coefficient for a compactly supported integrable potential $q$, which is not of the form (\ref{y-q}), but for which the bound states $\beta_1,\ldots,\beta_n$ and corresponding norming constants are known.  Then, by the work of Deift and Trubowitz \cite[Corollary, p.~188]{DeTr:1979}, multiplication of $R_2$ by an appropriate Blaschke factor transforms it into the reflection coefficient 
\[
\widetilde{R}(\omega)=(-1)^n\prod_{j=1}^n\frac{\omega-i\beta_j}{\omega+i\beta_j}R_2(\omega),
\]
corresponding to a positive definite potential $\widetilde{q}$.  Kappeler et al.~
 \cite[Thm.1.1, p.3092]{KaPeShTo:2005} have shown that such a $\widetilde{q}$ is necessarily of the form (\ref{y-q}), so that solving the associated scattering problem for the wave equation recovers $\widetilde{q}$ from $\widetilde{R}$.  And  \cite[Corollary, p.~188]{DeTr:1979} then prescribes how to reconstruct $q$ from $\widetilde{q}$.   Thus the requirement that $q$ be of the form (\ref{y-q}) does not limit the applicability of inverse scattering for the wave equation to that of the Schr\"odinger equation.  Indeed the aforementioned result of Kappeler et al.~
 applies to more singular potentials, those which are locally $H^{-1}$, which is a class of strictly lower regularity than required by \cite{DeTr:1979}. 

\subsubsection{An elementary operator\label{sec-simple-example}}
The following toy example illustrates another goal of the present work, namely, to characterize the nature of the scattering map beyond the mere observation that it is nonlinear.  Consider the Riccati equation (\ref{Riccati}) on $X=(0,1)$ with $\omega=0$ fixed and $\alpha=\alpha_0$ constant (but unknown), subject as before to the boundary condition $\mathbf{r}(1)=0$.  In this toy example the scattering map $\alpha_0\mapsto\mathbf{r}(0)$ is of course a nonlinear scalar function.  But it is elementary in the sense that it is a rational function of the exponential; more precisely,
\begin{equation}\label{toy}
\mathbf{r}(0)=\tanh\alpha_0,
\end{equation}
and so inverse scattering is simply evaluation of the inverse hyperbolic tangent function.  Stability of the inverse map 
\[
\mathbf{r}(0)\mapsto \alpha_0=\tanh^{-1}\bigl(\mathbf{r}(0)\bigr)
\]
decreases as the magnitude of $\bigl(\tanh^{-1}\bigr)^\prime\bigl(\mathbf{r}(0)\bigr)$ increases, with inversion becoming arbitrarily unstable as $|\mathbf{r}(0)|\rightarrow 1$. 

Returning to the full Riccati equation (\ref{Riccati}) with $\omega$-dependence restored and $\alpha\in L^1(X)$, is there a sense in which the full scattering map $\alpha\mapsto R$ is elementary?  Here the scattering map is a nonlinear operator, rather than a function of a scalar variable.  To say what it means for an operator to be elementary requires a reference operator analogous to the exponential function to serve as archetype. 
A formula for the harmonic exponential operator $E^X:L_\real^1(X)\rightarrow C^\infty(\real)$ is given below in \S\ref{sec-harmonic-exponential-0}.   
Taking the harmonic exponential as the archetypical elementary operator, the restriction of the scattering map to $\alpha\in L_\real^1(X)$ is shown in \S\ref{sec-forward-scattering} to be elementary.  In fact it is the operator theoretic analogue of hyperbolic tangent: 
\[
R=\frac{E^X_\alpha-E^X_{-\alpha}}{E^X_\alpha+E^X_{-\alpha}}.  
\]
The designation ``elementary" for the harmonic exponential, while justified to a certain extent by its basic properties as described in \S\ref{sec-regular-harmonic}, is perhaps open to debate---but such a label is ultimately a matter of convention rather than necessity. What is more important is the central role the harmonic exponential operator plays in scattering.

\subsubsection{Transforming initial conditions to boundary conditions\label{sec-transforming-initial}}
The hypothesis in \S\ref{sec-basic-problem} of constant wave speed on the interval $y\leq y_0$ corresponds in the context of equation (\ref{equation-0}) to constant impedance $\zeta(x)$ on the negative half-line $x\leq x_0$, where $u$ consequently has the form 
\[
u(x,\omega)=A(\omega)e^{-i(x-x_0)\omega}+B(\omega)e^{i(x-x_0)\omega}\qquad(x\leq x_0).
\]
The respective components
\begin{equation}\label{respective-components}
A(\omega)e^{-i(x-x_0)\omega}=\frac{1}{2}\left(u-\frac{1}{i\omega}u^\prime\right),\qquad B(\omega)e^{i(x-x_0)\omega}=\frac{1}{2}\left(u+\frac{1}{i\omega}u^\prime\right)
\end{equation}
are the Fourier transforms with respect to time of the left and right-moving wave forms
\begin{equation}\label{left-and-right-moving}
\rwc{A}(t+x-x_0)=d(t+x-x_0),\qquad\rwc{B}\bigl(t-(x-x_0)\bigr)=\delta(x-x_0-t)\qquad(x\leq x_0),
\end{equation}
since the measured data (\ref{data-0}) is left moving, and the initially right-moving delta (\ref{right-moving-delta}) transforms in travel time coordinates to $\delta(x-x_0-t)$.  (The formula (\ref{right-moving-delta}) is valid on the whole line at negative times, and at all times represents the restriction of the right-moving component of the wave field to the left half line $x\leq x_0$.) Thus $A=R$, $B=1$ and 
\begin{equation}\label{delayed-u-expression}
e^{ix_0\omega}u(x,\omega)=e^{2ix_0\omega}R(\omega)e^{-ix\omega}+e^{ix\omega}\qquad(x\leq x_0).
\end{equation}
From the perspective of acoustic imaging, where $d$ is measured only for a finite time, there is no loss of generality in assuming $\zeta^\prime$ has compact support, since in finite time $0<t<t_{\max}$ the data $d$ only manifests interaction of the initially right-moving impulse with the restriction of $\zeta$ to the spatial interval $x_0<x<x_0+t_{\max}/2$.  Thus suppose $\zeta^\prime$ is compactly supported and $q=\sqrt{\zeta}^{\,\prime\prime}/\sqrt{\zeta}$ is sufficiently smooth for standard existence and uniqueness results for (\ref{schroedinger}) to hold. 
Letting $f_1(x,\omega)$ denote the unique solution to (\ref{schroedinger}) such that $f_1(x,\omega)\sim e^{i\omega x}$ as $x\rightarrow\infty$, the entries $R_2$ and $T_2$ in the scattering matrix (\ref{scattering-matrix}) are defined by the asymptotic formula
\begin{equation}\label{R-T-definition}
T_2(\omega)f_1(x,\omega)\sim R_2(\omega)e^{-i\omega x}+e^{i\omega x}\quad\mbox{ as }\quad x\rightarrow -\infty.
\end{equation}
On the negative half-line $x\leq x_0$ the expression (\ref{delayed-u-expression}) satisfies (\ref{schroedinger}); comparison with (\ref{R-T-definition}) therefore yields (\ref{R2-R}).   Moreover, since $A=R$ and $B=1$, it follows from (\ref{respective-components}) that 
\begin{align}
R(\omega)&=\frac{1}{2}\left(u(x_0+,\omega)-\frac{1}{i\omega}u^\prime(x_0+,\omega)\right)\qquad and \label{R-formulation-0}\\
1&=\frac{1}{2}\left(u(x_0+,\omega)+\frac{1}{i\omega}u^\prime(x_0+,\omega)\right),\label{left-hand-boundary-0}
\end{align}
provided $u$ and $u^\prime$ are continuous at $x_0$.
(The notation $f(x+)$ and $f(x-)$ indicates one-sided limits.)  If $\zeta^\prime$ is supported on the closure of the interval $(x_0,x_1)$, then the initial conditions imply there is no left-moving wave on the right half-line $x\geq x_1$, and therefore 
\begin{equation}\label{right-hand-boundary-0}
\frac{1}{2}\left(u(x_1-,\omega)-\frac{1}{i\omega}u^\prime(x_1-,\omega)\right)=0,
\end{equation}
provided $u$ and $u^\prime$ are continuous at $x_1$.  The upshot is that the initial value problem (\ref{wave-0},\ref{initial-conditions-0}) transforms into an equivalent Sturm-Liouville problem with boundary conditions (\ref{left-hand-boundary-0},\ref{right-hand-boundary-0}), whose unique solution $u$ determines the reflection coefficient via (\ref{R-formulation-0}). The problem is not regular however, since the left-hand boundary condition (\ref{left-hand-boundary-0}) is affine and not linear.  

\subsubsection{The classical trace formula\label{sec-classical-trace}}
Hryniv and Mykytyuk \cite{HrMy:2021} have recently proved that for any real-valued potential $q\in L^1(\real)$, the reflection coefficient $R_2$ for the Schr\"odinger equation satisfies the trace formula
\[
-4\sum_{j=1}^n\kappa_j-\frac{1}{\pi}\int_\real\log\bigl(1-\bigl|R_2(\omega)\bigr|^2\bigr)\,d\omega=\int_\real q,
\]
where $\kappa_j$ are norming constants for bound states, thereby generalizing to $L^1(\real)$ a result long known to hold for a more restricted class of potentials \cite{FaZa:1971}.  
Suppose $\alpha\in L^2(\real)$ is real valued, and $q=\alpha^2-\alpha^\prime\in L^1(\real)$.  Then necessarily $\alpha^\prime\in L^1(\real)$, and so $\alpha$ is absolutely continuous; moreover, as mentioned earlier, $q$ has no bound states in this case. The above formula therefore reduces to 
\begin{equation}\label{classical-trace}
-\int_\real\log\bigl(1-\bigl|R_2(\omega)\bigr|^2\bigr)\,d\omega=\pi\int_\real |\alpha|^2.
\end{equation}
The latter is the nonlinear Plancherel theorem asserted in \cite{SyWi:1998} to hold for all $\alpha\in L^2(\real)$.  But there is an essential error in the geometric arguments put forward in \cite{SyWi:1998}, as follows. A bivariate function on the unit bidisk is defined there as 
\[
d_e:\ddd^2\rightarrow\real_+,\quad d_e(z,w)=-\log\left(1-\left|\frac{w-z}{1-\bar{z}w}\right|^2\right),
\]
so that for $0<t<1$,  
\[
d_e(-t,0)=d_e(0,t)=-\log(1-t^2)\quad\mbox{ and }\quad d_e(-t,t)=2\log\bigl((1+t^2)/(1-t^2)\bigr).
\]
See \cite[pp.~673,676]{SyWi:1998}. Hence 
\[
-2\log(1-t^2)+2\log(1+t^2)=d_e(-t,t)>d_e(-t,0)+d_e(0,t)=-2\log(1-t^2),
\]
whereby $d_e$ fails to satisfy the triangle inequality and is thus not a metric.  However, in \cite[Lem.~1.5, Lem.~1.8, Cor.~1.9]{SyWi:1998} the function $d_e$ is explicitly assumed to be a metric, an error which undermines subsequent geometric arguments related to the trace formula.

The present work therefore avoids reliance on \cite{SyWi:1998}.  All that is needed for the inverse scattering formula is the nonlinear Bessel's inequality for (\ref{classical-trace}), 
\[
-\int_\real\log\bigl(1-\bigl|R_2(\omega)\bigr|^2\bigr)\,d\omega\leq\pi\int_\real |\alpha|^2,
\]
and this is proved directly by singular approximation.  The full Plancherel theorem for real-valued $\alpha\in L^2(\real)$ apparently  remains open (the absolutely continuous case being covered by \cite{HrMy:2021}). 

\subsection{Technical framework\label{sec-technical-framework}}
This section describes the key ingredients of regulated functions and the generalized reflection coefficient, and establishes terminology and notation to be used throughout the paper. 
\subsubsection{Regulated functions\label{sec-regulated}}
Regulated functions, introduced by Bourbaki for integration of Banach space-valued functions of a real variable (see \cite[Ch.VII,\S6]{Di:1960} for background) provide just the right level of generality for the coefficient in (\ref{equation-0}).
Given a bounded interval $X=(x_0,x_1)$, let $\step(X)$ denote the set of all step functions on $X$, i.e., piecewise-constant complex-valued functions having finitely many points of discontinuity, and let $\step_+(X)$ denote the strictly positive step functions.  Regulated functions on $X$ are defined to be the uniform closure of step functions, 
\begin{equation}\label{regulated-definition}
\reg(X)=\overline{\step(X)}\subset L^\infty(X). 
\end{equation}
Thus $\reg(X)$ is a closed linear subspace of $L^\infty(X)$.
If $f\in\reg(X)$, then the left and right-hand limits $f(x-)$ and $f(x+)$ exist for every $x\in X$, as do the limits $f(x_0+)$ and $f(x_1-)$.  (Existence of one-sided limits is equivalent to being a uniform limit of step functions.)  
Every function of bounded variation is regulated, but not conversely.
Set 
\begin{equation}\label{positive-regulated}
\reg_+(X)=\left\{f\in\reg(X)\,\left|\,\exists\varepsilon>0\mbox{ s.t. }f(x)>\varepsilon \forall x\in X\right.\right\}.
\end{equation}

Of particular interest is the subclass of $\reg_+(X)$ consisting of piecewise absolutely continuous functions.  (Absolute continuity, as opposed to just continuity, reduces complications from the technical mathematical standpoint; and the difference does not seem to be significant from the physical perspective.)
Fix notation as follows.  Let $C_+(X)\subset\reg_+(X)$ denote absolutely continuous functions which are positive, regulated and bounded away from 0. 
Define $\pwac(X)$ to be the set of all functions $f:X\rightarrow\complex$ for which there exists a finite partition of $X$, 
\[
x_0=y_0<y_1<\cdots<y_n=x_1
\]
such that $f|_{(y_{j-1},y_j)}$ is absolutely continuous for every $1\leq j\leq n$. 
Set 
\begin{equation}\label{positive-piecewise}
\pwacplus(X)=\pwac(X)\cap\reg_+(X). 
\end{equation}
Each set of functions in the increasing sequence
\[
\step_+(X)\subset\pwacplus(X)\subset\reg_+(X)
\]
is easily seen to be a multiplicative group, and every $\zeta\in\pwacplus(X)$ has a unique factorization $\zeta=\zeta_1\zeta_2$ where $\zeta_1$ is absolutely continuous and $\zeta_2\in\step_+(X)$ with $\zeta_2(x_0+)=1$.  (Technical note: uniqueness is in the context of $L^\infty(X)$, where two functions that differ on a set of measure zero are identified.  One is free to work with equivalence class representatives having no removable discontinuities.)

Denote by $C^1_+(X)\subset\reg_+(X)$ the set of positive continuously differentiable regulated functions which are bounded away from 0, and whose derivative is also regulated. The significance of the latter condition is to guarantee that $\alpha=-\zeta^\prime/(2\zeta)$ is bounded if $\zeta\in C^1_+(X)$.

\subsubsection{Regularity of solutions\label{sec-regularity-of-solutions}}

The following proposition gives several equivalent formulations of the appropriate notion of solution for (\ref{equation-0}).
\begin{prop}\label{prop-solution-interpretation}
Given $\zeta\in\reg_+(X)$ and $\omega\in\real$, let $\Lop$ denote the distribution-valued operator on $W^{1,2}(X)$ 
\[
\Lop u=(\zeta u^\prime)^\prime+\omega^2\zeta u.
\]
The following conditions are equivalent.
\begin{enumerate}
\item $\Lop u=0$ almost everywhere on $X$, and $u$ and $\zeta u^\prime$ are absolutely continuous.
\item $u\in W^{1,2}(X)$, $\Lop u\in L^2(X)$, and $\langle\Lop u,v\rangle_{L^2}=0$ for every $v\in W^{1,2}(X)$. 
\item $u\in W^{1,2}(X)$ and $\Lop u=0$ in the sense of distributions. 
\end{enumerate}
\end{prop}
\begin{pf}
(2) $\Rightarrow$ (3) trivially since the set $C_c^\infty(X)$ of test functions is a subset of $W^{1,2}(X)$.\\
(3) $\Rightarrow$ (1).  $\Lop u$ being the zero distribution implies that $\Lop u$ corresponds to a bona fide function that is zero almost everywhere. Since $u\in W^{1,2}(X)$, $u$ is absolutely continuous, $\omega^2\zeta u\in L^2(X)$, and therefore $(\zeta u^\prime)^\prime\in L^2(X)$, since $\Lop u=0$ almost everywhere. It follows that $\zeta u^\prime$ is absolutely continuous. \\
(1) $\Rightarrow$ (2).  $u\in L^2(X)$ if $u$ is absolutely continuous on $X$, since the latter condition implies boundedness.  And if $\zeta u^\prime$ is absolutely continuous, then $u^\prime=(1/\zeta)\zeta u^\prime\in L^2(X)$ since $\zeta$ is bounded away from 0.  Thus $u,u^\prime\in L^2(X)$, meaning that $u\in W^{1,2}(X)$.  Given that $\Lop u=0$ almost everywhere it follows trivially that $\Lop u\in L^2(X)$, and $\langle\Lop u,v\rangle_{L^2}=0$ for every $v\in W^{1,2}(X)$. 
\end{pf}

\subsubsection{Automorphisms of the Poincar\'e disk\label{sec-automorphisms}}

A richly structured object central to this paper is the closure $K$ of the automorphism group of the Poincar\'e disk, endowed with a certain Riemannian structure.  More explicitly, 
\[
K=\overline{\aut\mathbb{P}}
\]
consists of the set of all holomorphic bijections of the unit disk, functions of the form
\begin{equation}\label{disk-automorphisms}
\varphi_{\mu,\rho}:\ddd\rightarrow\ddd\qquad \varphi_{\mu,\rho}(\xi)=\mu\frac{\xi+\rho}{1+\bar\rho\xi}\qquad(\mu\in S^1,\rho\in\ddd),
\end{equation}
together with all constant maps having values in the unit circle
\begin{equation}\label{constant-maps}
\mathbf{c}_\sigma:\ddd\rightarrow S^1\qquad \mathbf{c}_\sigma(\xi)=\sigma\qquad(\sigma\in S^1). 
\end{equation}  
Here closure is with respect to the topology of uniform convergence on compact sets.  Of course, $\aut\mathbb{P}\cong {\rm PSL}(2,\real)\cong S^1\times\ddd$.  
Adjoining the circle $S^1$ of constant maps yields a compact topological space homeomorphic to the 3-sphere, 
\begin{equation}\label{K-three-sphere}
K=\aut\mathbb{P}\cup S^1\cong S^3. 
\end{equation}
Non-invertibility of the constant maps mean that $K$ is a Lie semigroup, not a group---unlike the Lie group spin$(3)$, for example, which is also homeomorphic to $S^3$.  Importantly, $K$ carries a Riemannian structure (degenerate at the circle of constant maps) directly relevant to scattering; see the discussion preceding Theorem~\ref{thm-fourier-series} in \S\ref{sec-OPUD}. 

Let $K^\real$ denote the set of functions $f:\real\rightarrow K$ defined almost everywhere, endowed with the convergence structure of almost-everywhere pointwise convergence. Strictly speaking, this convergence structure does not come from a topology; nevertheless, a mapping from a metric space into $K^\real$ will henceforth be described as \emph{continuous} when the image of every convergent sequence in the domain converges pointwise almost everywhere. 

\subsubsection{The generalized reflection coefficient\label{sec-generalized-reflection-coefficient-2}}

For $\omega\in\real\setminus\{0\}$, denote left and right-moving components of a solution to equation (\ref{equation-0})
by the formulas
\begin{equation}\label{left-right-components}
u^\ell=\frac{1}{2}\left(u-\frac{1}{i\omega}u^\prime\right)\qquad
u^{\rm r}=\frac{1}{2}\left(u+\frac{1}{i\omega}u^\prime\right).
\end{equation}
The terminology of left and right-moving stems from the wave equation in the time domain, where, if $\zeta$ is constant near some point, the inverse Fourier transforms of $u^{\ell}$ and $u^{\rm r}$ locally represent the respective left and right-moving components of the wave field, as in (\ref{left-and-right-moving}).  For brevity $u(x,\omega)$ will henceforth be denoted just $u(x)$, dependency on $\omega$ being understood; the same convention will apply also to $u^\ell$ and $u^{\rm r}$. 

In terms of the above notation, if $\zeta^\prime$ is supported on the closure of $X=(x_0,x_1)$, the initial value problem (\ref{wave-0},\ref{initial-conditions-0}) transforms according to \S\ref{sec-transforming-initial} into an equivalent Sturm-Liouville problem on $X$ of the form
\begin{equation}\label{sturm-liouville-0}
\begin{split}
(\zeta u^\prime)^\prime+\omega^2\zeta u&=0\\
u^{\rm r}(x_0+)=1,\qquad &u^\ell(x_1-)=0,
\end{split}
\end{equation}
with corresponding reflection coefficient $R=u^\ell(x_0+)$.
The generalized reflection coefficient is obtained by introducing a complex variable into the right-hand boundary condition as follows. 
For each $\omega\in\real\setminus\{0\}$, consider the parameterized family of boundary value problems on $X$,
\begin{gather}
(\zeta u^\prime)^\prime+\omega^2\zeta u=0\label{equation}\\[5pt]
u^{\rm r}(x_0+)=1\quad\mbox{ and }\quad u^\ell(x_1-)=\xi u^{\rm r}(x_1-)\qquad(\xi\in\ddd).\label{bc}
\end{gather}
The system (\ref{equation},\ref{bc}) determines what turns out to be a continuous map
\begin{equation}\label{continuous-map}
g:\reg_+(X)\rightarrow K^\real,
\end{equation}
the proof of which is a major objective of \S\ref{sec-mapping}.  The mapping (\ref{continuous-map}) is to be interpreted as follows (the various assertions are proved later, in \S\ref{sec-mapping} and \S\ref{sec-scattering-for-piecewise-continuous}). For any fixed $\zeta\in\reg_+(X)$, for almost every $\omega\in\real$, the system (\ref{equation},\ref{bc}) has a unique solution as stipulated by Proposition~\ref{prop-solution-interpretation}, i.e., in the sense of distributions, with $u\in W^{1,2}(X)$,  for every $\xi\in\ddd$.  For each such solution the value $u^{\ell}(x_0+)$ is well defined.  Set
\begin{equation}\label{g-definition}
g_\zeta^\omega(\xi)=u^{\ell}(x_0+).  
\end{equation}
Then $g_\zeta\in K^\real$ in the sense that $g_\zeta^\omega\in K$ for almost every $\omega\in\real$.  In other words, for fixed $\zeta$ and $\omega$, $u^{\ell}(x_0+)$, viewed as a function of the boundary parameter $\xi$, is a disk automorphism (or, in the most pathological case, constant).  
This $\omega$-parameterized family of disk automorphisms $g_\zeta$ is called the generalized reflection coefficient.
Since for $\xi=0$ the system (\ref{equation},\ref{bc}) reduces to (\ref{sturm-liouville-0}), the usual reflection coefficient is recovered by evaluation at $\xi=0$,
\begin{equation}\label{R-from-g}
R(\omega)=g_\zeta^\omega(0).
\end{equation}

For any $\zeta\in\reg_+(X)$, the values $g_\zeta$ and $g_{1/\zeta}$ are related very simply by the formula
\begin{equation}\label{zeta-reciprocal}
g^\omega_{1/\zeta}(\xi)=-g^\omega_\zeta(-\xi),
\end{equation}
and for any $s>0$, $g_{s\zeta}=g_\zeta$.  Thus $g_\zeta$ is homogeneous of degree zero with respect to $\zeta$. 
The notation just described will be used throughout the paper. To reiterate: $g_\zeta$ denotes a map from $\real$ to $K$; its value at $\omega\in\real$ is denoted $g_\zeta^\omega$; the latter is itself a map on the unit disk having values $g_\zeta^\omega(\xi)$.

If $\zeta(x_0+)$ is fixed, the mapping (\ref{continuous-map}) is injective on $\pwacplus(X)$, and the corresponding values are smooth functions $g_\zeta:\real\rightarrow K$, which may be thought of as curves in $S^3$. The main focus of the present paper is this restricted mapping $g|_{\pwacplus(X)}$.

\subsubsection{Concatenation of media\label{sec-concatenation}}

Consider adjacent intervals $X_1=(x_0,x_1)$ and $X_2=(x_1,x_2)$ for some $x_0<x_1<x_2$, and let $X=(x_0,x_2)$ denote their connected union.  An impedance function $\zeta\in\reg_+(X)$ (representing a layered medium, say) is said to be the concatenation of $\zeta_1\in\reg_+(X_1)$ and $\zeta_2\in\reg_+(X_2)$ if $\zeta_1=\zeta|_{X_1}$ and $\zeta_2=\zeta|_{X_2}$. The notation for concatenation is
\begin{equation}\label{concatenation}
\zeta=\zeta_1^{\;\frown}\zeta_2.
\end{equation}
An advantage of analyzing scattering on a bounded interval as opposed to the half line or the line is to facilitate a description of the relationship between concatenation and scattering:
if $\zeta=\zeta_1^{\;\frown}\zeta_2$ is continuous at $x_1$, then concatenation corresponds to composition in $K$,
\begin{equation}\label{concatenation-scattering}
g_{\zeta_1^{\;\frown}\zeta_2}^\omega=g_{\zeta_1}^\omega\circ g_{\zeta_2}^\omega. 
\end{equation}
Letting $R, R_1, R_2$ denote the reflection coefficients corresponding to $\zeta_1^{\;\frown}\zeta_2, \zeta_1, \zeta_2$, respectively, the compositional relation (\ref{concatenation-scattering}) together with (\ref{R-from-g}) implies
\begin{equation}\label{concatenation-reflection}
R(\omega)=g_{\zeta_1}^\omega\bigl(R_2(\omega)\bigr),
\end{equation}
which cannot be expressed in a simple way in terms of just $R,R_1$ and $R_2$.
If $\zeta$ fails to be continuous at $x_1$ the formula has to be modified accordingly; see Theorem~\ref{thm-g-composition}, \S\ref{sec-general-existence}.  

\subsubsection{OPUC and OPUD\label{sec-OPUC-OPUD}}

Two flavours of orthogonal polynomials play a role in the present work.  Orthogonal polynomials on the unit circle (OPUC) arise in the representation of the reflection coefficient corresponding to evenly-spaced step functions.  The notation used here follows that of \cite{SiOPUC1:2005}: $\Phi_n(z)$ denotes the degree $n$ monic orthogonal polynomial determined by a given sequence of Verblunsky parameters $\rho_j\in\ddd$ $(1\leq j\leq n)$; $\Psi_n(x)$ denotes the associated degree $n$ monic polynomial determined by $-\rho_j$ $(1\leq j\leq n)$; their dual polynomials are 
\[
\Phi_n^\ast(z)=z^j\overline{\Phi_n(1/\bar{z})}\quad\mbox{ and }\quad\Psi_n^\ast(z)=z^j\overline{\Psi_n(1/\bar{z})}. 
\]
See \S\ref{sec-OPUC} for a brief summary of the relevant theory. 

Orthogonal polynomials on the unit disk (OPUD) are bivariate polynomials that also arise in the representation of the reflection coefficient, but in a very different way from OPUC.  A single family of OPUD, called scattering polynomials, comprises the natural basis functions in terms of which to expand the reflection coefficient corresponding to any $\zeta\in\step_+(X)$, irrespective of how the jump points are spaced.  Scattering polynomials are eigenfunctions of the Laplace-Beltrami operator for the scattering metric associated to PSL$(2,\real)$. A complete description of scattering polynomials is given in \S\ref{sec-OPUD}. 

\subsubsection{The harmonic exponential operator\label{sec-harmonic-exponential-0}}
As mentioned already in \S\ref{sec-simple-example}, the regular harmonic exponential, which is treated in detail in \S\ref{sec-regular-harmonic}, emerges as a useful technical intermediary in the context of the scattering map.  Given $X=(x_0,x_1)$ and $y\in(x_0,x_1]$, it is an operator
\[
E^{(x_0,y)}:L^1_\real(X)
\rightarrow C^\infty(\real),\qquad \alpha\mapsto E^{(x_0,y)}_\alpha
\]
defined by the formula
\[
E^{(x_0,y)}_\alpha(\omega)=1+\displaystyle\sum_{j=1}^\infty\thickspace\int\limits_{x_0<s_1<\cdots<s_j<y}\negthickspace\negthickspace\exp\bigl(2i\omega\textstyle\sum_{\nu=1}^j(-1)^{j-\nu}(s_\nu-x_0)\bigr)\displaystyle\prod_{\nu=1}^j\alpha(s_\nu)\,ds_1\cdots ds_j.
\]
Here the subscript $\mathbb{R}$ in $L^1_\mathbb{R}(X)$ indicates real-valued functions.
The term ``exponential'' is justified in part by the facts that $E^{(x_0,y)}_\alpha(0)=\exp\int_{x_0}^y\alpha$ and $E^{(x_0,y)}_\alpha(\omega)$ is the restriction to $\omega\in\real$ of an entire function; ``harmonic" refers to the exponential factor in the defining formula, which is loosely analogous to the exponential factor in the Fourier transform.  In the particular case where $\alpha(x)=\alpha_0$ is constant, the harmonic exponential operator at  $\alpha$, evaluated at $\omega=0$, is just the scalar exponential
$e^{\alpha_0(y-x_0)}$
(cf.~\S\ref{sec-simple-example}). 

A corresponding singular harmonic exponential operator is defined in \S\ref{sec-singular-harmonic}. The latter definition involves the reflectivity function $r$ associated to any $\zeta\in\reg_+(X)$, which is given by the formula
\begin{equation}\label{reflectivity-0} 
r(y)=\frac{\zeta(y-)-\zeta(y+)}{\zeta(y-)+\zeta(y+)}\qquad(y\in X).
\end{equation}
A special case of the singular harmonic exponential may be described in terms of OPUC as follows.  Let $\zeta_n\in\step_+(X)$ have $n$ equally-spaced jump points in the interval $(x_0,y)$,
\[
y_j=x_0+j\Delta_n\qquad(1\leq j\leq n),
\]
where $\Delta_n=(y-x_0)/(n+1)$, and set 
\[
\rho_j=r(y_j)\qquad(1\leq j\leq n).
\]
Let $\Psi^\ast_n(z)$ denote the degree $n$ dual OPUC determined by Verblunsky parameters $-\rho_1,\ldots,-\rho_n$.  The singular harmonic exponential on $(x_0,y)$ associated to $\zeta_n$ is then
\[
\calE_{\zeta_n}^{(x_0,y)}(\omega)=\Psi^\ast_n\bigl(e^{2i\Delta_n\omega}\bigr).
\]

If $\zeta\in C_+(X)$, there is a sequence of equally-spaced step functions $\zeta_n$ $(n\geq 1)$ converging uniformly to $\zeta$ such that 
\[
\calE_{\zeta_n}^{(x_0,y)}(\omega)\rightarrow E_\alpha^{(x_0,y)}(\omega)
\]
uniformly on compact sets, where $\alpha=-\zeta^\prime/(2\zeta)$.  
Thus $E_\alpha^{(x_0,y)}$ is a continuum limit of orthogonal polynomials.  This realization of the harmonic exponential as a continuum limit of OPUC is an essential feature of singular approximation. 

\subsubsection{Almost periodic functions\label{sec-almost-periodic}}  If $\zeta\in\step_+(X)$, then the reflection coefficient $R(\omega)=g_\zeta^\omega(0)$ is almost periodic in the sense of Besicovitch \cite{Be:1926, Be:1955}.  Summation over the continuum provides a useful notation in this context: an expression of the form
\[
\sum_{\lambda\in\real}f(\lambda)
\]
shall be understood to imply that $f(\lambda)=0$ except for at most countably many $\lambda\in\real$,  with only non-zero values of $f(\lambda)$ contributing to the sum.  Almost periodic functions in the sense of Besicovitch comprise a non-separable Hilbert space consisting of equivalence classes of functions $f:\real\rightarrow\complex$ of the form
\begin{equation}\label{almost-periodic-representative}
f(\omega)=\sum_{\lambda\in\real}c_\lambda e^{i\lambda\omega}\quad\mbox{ such that }\quad \sum_{\lambda\in\real}|c_\lambda|^2<\infty.
\end{equation}
The scalar product is
\[
\langle f,h\rangle_{\rm ap}=\lim_{L\rightarrow\infty}\frac{1}{2L}\int_{-L}^Lf(\omega)\overline{h(\omega)}\,d\omega.
\]
If $\|h\|_{\rm ap}<\infty$ for a measurable $h:\real\rightarrow\complex$, then, as an almost periodic function, $h$ is equivalent to a representative $f$ of the form (\ref{almost-periodic-representative}) whose coefficients are determined by 
\begin{equation}\label{almost-periodic-coefficient}
c_\lambda=\langle h,e^{i\lambda\omega}\rangle_{\rm ap}.
\end{equation}
In particular, if $h\in L^2(\real)$, then $\|h\|_{\rm ap}=0$, i.e., as an almost periodic function $h=0$.  In general, the representative $f$ of the form (\ref{almost-periodic-representative}) determined by (\ref{almost-periodic-coefficient}) will be called the almost periodic part of $h$.  

Writing  $\zeta\in\pwacplus(X)$ in terms of its factorization $\zeta=\zeta_1\zeta_2$, where $\zeta_1$ is continuous and $\zeta_2\in\step_+(X)$, it will be proven that 
\[
\|g_\zeta(0)-g_{\zeta_2}(0)\|_{\rm ap}=0.
\]
Thus the almost periodic part of the reflection coefficient $g_\zeta(0)$ is determined by $\zeta_2$; in particular, the almost periodic part of $g_\zeta(0)$ is zero if and only if $\zeta$ is continuous.  The classical trace formula (\ref{classical-trace}) breaks down if $R_2$ has a non-zero almost periodic part, since in this case the left-hand integral diverges and $\alpha$ acquires a purely singular distributional part.

\subsection{Overview of the paper\label{sec-strategy}}

The overall organization of the paper is dictated by its basic strategy of singular approximation, the idea of which is to approximate impedance functions $\zeta\in \reg_+(X)$ by step functions, and then to infer properties of the scattering map $\zeta\mapsto R$ from explicit formulas available for step functions.  The approximation is purely singular in the sense that $\alpha=-\frac{1}{2}\bigl(\log\zeta\bigr)^\prime$ is a purely singular distribution (a finite combination of Dirac pulses) if $\zeta\in\step_+(X)$.  The corresponding potential $q=\alpha^2-\alpha^\prime$ is not even a well-defined distribution, which means the strategy of singular approximation cannot be implemented directly in terms of the standard Schr\"odinger equation---either the Riccati or Helmholtz equation is needed.  Implementation involves two main steps, corresponding to two key approximation lemmas.  
The first step is to prove that the generalized reflection coefficient $g:\reg_+(X)\rightarrow K^\mathbb{R}$ is well defined and continuous. 
The second step is to tease out the structure of $g_\zeta^\omega(\xi)$ as a function of $\omega$ in the case where $\zeta\in \pwacplus(X)$. 

Section~\ref{sec-mapping} is devoted to the first step. This entails, among other things, proving the basic existence and uniqueness results for the given Sturm-Liouville problem. The proof, which follows a standard arc of formulating a notion of weak solution and then proving that weak solutions are strong solutions, leans heavily on the Lax-Milgram lemma. The key approximation result is Lemma~\ref{lem-forward-stability}, asserting continuity of the solution map $\zeta\mapsto u$ relative to respective spaces $\reg_+(X)$ and $W^{1,2}$. Theorems~\ref{thm-general-existence}, \ref{thm-g-in-K} and \ref{thm-continuous}, stated in \S\ref{sec-general-existence}, provide, respectively, existence of solutions, that $g_\zeta^\omega$ is a disk automorphism, and that the map $g:\reg_+(X)\rightarrow K^\mathbb{R}$ is continuous. Theorem~\ref{thm-g-composition} in the same section states the concatenation law, crucial for passing from $C_+(X)$ to $\pwacplus(X)$.  

Section~\ref{sec-harmonic-exponential} concerns the key technical intermediary, the harmonic exponential operator. 
The topology assigned to the codomain of the map $g:\reg_+(X)\rightarrow K^\mathbb{R}$, that of almost-everywhere pointwise convergence, is too weak for the continuity result of \S\ref{sec-mapping} to yield information concerning $\omega$-dependency of the generalized reflection coefficient.  So a second step is required. This involves the harmonic exponential operator, and, more precisely, approximation of the regular harmonic operator by the singular.  The second key approximation result is Lemma~\ref{lem-singular-approximation}. Stated in \S\ref{sec-approximation}, it gives the stronger convergence needed to infer $\omega$-dependence, roughly as follows.  For a special approximating sequence $\zeta_n$ converging uniformly to $\zeta$ (called the standard approximation), the corresponding values of the singular harmonic exponential operator converge uniformly on compact sets to the regular harmonic exponential at $\zeta$, provided that $\zeta$ itself is $C^1$. Another approximation result pertaining to the regular harmonic exponential is needed to pass from $C^1$ impedance to absolutely continuous impedance.  This is Proposition~\ref{prop-harmonic-2}\ref{E-continuous}, which says that the regular harmonic exponential operator $E^X:L_\mathbb{R}^1(X)\rightarrow L^\infty(\mathbb{R})$ is continuous. 
Section~\ref{sec-harmonic-exponential} also compiles results concerning the harmonic exponential operator that are needed for the main scattering results. Subsection \ref{sec-OPUD}, which deals with OPUD, has a different flavour from the other subsections. It analyzes the almost periodic structure that emerges when a given $\zeta\in \pwacplus(X)$ has discontinuities, and facilitates explicit representation for the almost periodic part of $g_\zeta$, for arbitrary $\zeta\in \pwacplus(X)$. 

Roughly speaking, the details of singular approximation are worked out in \S\ref{sec-mapping} and \S\ref{sec-harmonic-exponential}. These develop the technical machinery needed to analyze forward and inverse scattering for piecewise continuous $\zeta$. The development culminates in \S\ref{sec-scattering-for-piecewise-continuous}, which is the heart of the paper. Four theorems describe forward and inverse scattering in terms of explicit formulas, and two additional theorems 
convey
a singular trace  formula. In \S\ref{sec-forward-scattering}, which deals with forward scattering, Theorem~\ref{thm-hyperbolic} asserts that if $\zeta$ is absolutely continuous, then forward scattering is effected by the hyperbolic tangent operator associated to the harmonic exponential. Theorem~\ref{thm-forward-scattering} incorporates the concatenation law to give an explicit formula for the reflection coefficient in the more general case of $\zeta\in \pwacplus(X)$. 
The main results concerning inverse scattering are presented in \S\ref{sec-inverse-scattering}. The short-range inversion formula of Theorem~\ref{thm-short-range-inversion} expresses $\zeta(x)$ directly in terms of the reflection coefficient $R$, for all $x$ sufficiently near the measurement location $x_0$. 
This appears to be the first example of such a formula. 
If $\zeta$ is absolutely continuous, then the short-range inversion formula determines the restriction $\zeta_1$ of $\zeta$ to an initial interval $(x_0,b)$ of the full interval $(x_0,x_1)$, which in turn determines the inverse $g_{\zeta_1}^{-1}$ of the disk automorphism $g_{\zeta_1}$. The reflection coefficient corresponding to the remainder of the initial interval $(b,x_1)$ is then given by $g_{\zeta_1}^{-1}(R)$, to which the short-range inversion formula may be applied once again. The crucial point here is that the length of the interval of validity of the short-range inversion formula in this second application is at least as great as in the first. Thus a finite number of iterations suffices to recover all of $\zeta$. In the case where $\zeta$ has discontinuities the above procedure recovers $\zeta$ up to the first jump point. The location of the jump and its height, however, are encoded in the almost periodic part of the reflection coefficient, allowing iterative application of the short-range inversion formula to be continued past the jump point, to ultimately recover $\zeta$ in full. This is the content of Theorem~\ref{thm-piecewise-injectivity}.  Section~\ref{sec-renormalized-trace} addresses the trace formula discussed in above in \S\ref{sec-classical-trace}, which breaks down if $\zeta$ is discontinuous.  Theorem~\ref{thm-Szego} quotes a Szeg\"o-type result valid in the almost periodic context, which is proved in \cite{Gi:Szego2022}.  Theorem~\ref{thm-renormalized-trace} incorporates this into a singular trace formula valid for discontinuous $\zeta\in\pwacplus(X)$.

Section~\ref{sec-fast-algorithms} links the theoretical results established in the earlier parts of the paper to the practical problem of inverting digitally recorded acoustic reflection data, which motivated the theoretical work in the first place. The main goal of the section is proof of concept, and focus is deliberately limited. Two algorithms are presented, one for forward scattering in \S\ref{sec-forward-computation}, and one for inverse scattering in \S\ref{sec-inverse-computation}. These are illustrated in \S\ref{sec-computation-example} with a concrete example. 
Over the years, many
different computational schemes for inverse scattering have been proposed, typically involving the numerical solution of integral equations (see the survey \cite{Ne:1980}), but also based on other discretization schemes. Operation counts tied to the digital nature of recorded data seem to be largely absent however, precluding easy comparisons of efficiency. Both the forward and inverse formulas in Theorems~\ref{thm-hyperbolic} and \ref{thm-short-range-inversion} are amenable to direct numerical approximation. But while such direct approximation is feasible in principle, it is in practice very far from the most computationally efficient approach. What distinguishes the two algorithms presented in \S\ref{sec-fast-algorithms} (in particular from those proposed in \cite{WaAk:1969}, \cite{BuBu:1983} and \cite[III.5]{FoFr:1990}) is that they exploit the recursive structure inherent in OPUC for maximal efficiency. In particular, the computation of moments step in Algorithm~\ref{alg-inverse} reduces the operation count from $O(n^3)$ to $O(n^2)$. One would expect inverse scattering to be computationally more expensive than forward scattering, but, surprisingly, the opposite is true, at least at the data scale $n\leq 10000$. 

Standard scattering lore describes the scattering map as a nonlinear analogue of the Fourier transform, making its discrete version analogous to the DFT. The FFT famously reduces computational expense of the DFT to $O(n\log n)$ by exploiting Vandermonde symmetry. One cannot expect such results for the less symmetric nonlinear analogue (where, for example, the forward and inverse maps are completely different). But it seems reasonable to regard the $O(n^2)$ algorithms of \S\ref{sec-fast-algorithms} as nonlinear analogues of the FFT in that they take maximum advantage of the special structure available, namely that of OPUC.  
There is much more to be said concerning computation, however this will be consigned to separate work.

In \S\ref{sec-conclusion} the paper concludes with some technical and historical remarks, and a short discussion of open problems.

\section{Continuity of $g:\reg_+(X)\rightarrow K^\real$\label{sec-mapping}}

It is proved in this section that the system (\ref{equation},\ref{bc}) has a unique solution $u$ for almost every $\omega$, and that setting $g_\zeta^\omega(\xi)=u(x_0+)$ defines a continuous map $g:\reg_+(X)\rightarrow K^\real$.  

\subsection{Existence and uniqueness for $\zeta\in\step_+(X)$\label{sec-existence-step}}

Here the essential facts concerning the system (\ref{equation},\ref{bc}) in the case where $\zeta\in\step_+(X)$ are summarized.  Existence and uniqueness is established, as is the fact that $g_\zeta^\omega$ is a disk automorphism. The basic results are well known, but the present formulation is tailored to suit approximation of the general case. 

To begin, consider scattering at a discontinuity between two intervals on which $\zeta$ is constant.  Let $\zeta\in\reg_+(X)$  be constant on $(y_{j-1},y_j)$ and $(y_j,y_{j+1})$ for some $$x_0\leq y_{j-1}<y_j<y_{j+1}\leq x_1.$$ Assuming $\omega\in\real\setminus\{0\}$, the general solution to (\ref{equation}), restricted to these two intervals, may be expressed as 
\begin{equation}\label{general-solution-two-intervals}
u(x)=\left\{
\begin{array}{cc}
A_{j-1}e^{-i(x-y_{j-1})\omega}+B_{j-1}e^{i(x-y_{j-1})\omega}&\mbox{ if }y_{j-1}<x<y_{j}\\
A_{j}e^{-i(x-y_{j})\omega}+B_{j}e^{i(x-y_j)\omega}&\mbox{ if }y_{j}<x<y_{j+1}
\end{array}\right.
\end{equation}
where constants $A_{j-1},B_{j-1},A_j,B_j$ have the representation
\begin{equation}\label{constants-as-directional-components}
A_{j-1}=u^{\ell}(y_{j-1}+),\quad  B_{j-1}=u^{\rm r}(y_{j-1}+),\quad  A_j=u^{\ell}(y_j+),\quad  B_j=u^{\rm r}(y_j+). 
\end{equation}
Statement (1) of Proposition~\ref{prop-solution-interpretation} requires that the left and right-hand limits of $u$ and $\zeta u^\prime$ agree at $y_j$, leading directly to the relation
\begin{equation}\label{scattering-at-yj}
\binom{A_{j-1}}{B_{j-1}}=\frac{1}{\mu_j^{1/2}(1+r_j)}\begin{pmatrix}\mu_j&\mu_jr_j\\ r_j&1\end{pmatrix}
\binom{A_j}{B_j}
\end{equation}
where 
\begin{equation}\label{muj-rj}
\mu_j=e^{2i(y_j-y_{j-1})\omega}\quad\mbox{ and }\quad r_j=\frac{\zeta(y_j-)-\zeta(y_j+)}{\zeta(y_j-)+\zeta(y_j+)}.
\end{equation}
Here the intended interpretation of the square root is $\mu_j^{1/2}=e^{i(y_j-y_{j-1})\omega}$.  The invertible relation (\ref{scattering-at-yj}) is key, implying 
\begin{align}
\zeta(y_j-)|A_{j-1}|^2+\zeta(y_j+)|B_j|^2=&\zeta(y_j-)|B_{j-1}|^2+\zeta(y_j+)|A_j|^2\label{energy-1}\\
\mu_j^{-1/2}\zeta(y_j-)A_{j-1}+\zeta(y_j+)B_j=&\mu_j^{-1/2}\zeta(y_j-)B_{j-1}+\zeta(y_j+)A_j\label{momentum-1}
\end{align}
and, in terms of the notation (\ref{disk-automorphisms}) for automorphisms of the Poincar\'e disk, 
\begin{equation}\label{disk-automorphism-at-yj}
A_{j-1}/B_{j-1}=\varphi_{\mu_j,r_j}\bigl(A_j/B_j\bigr).
\end{equation}
Equations (\ref{energy-1}) and (\ref{momentum-1}) express conservation of energy and momentum, respectively. 
Together they imply (\ref{scattering-at-yj}), which was originally derived based on the purely technical requirement that (\ref{equation}) should have a meaningful interpretation in the sense of distributions, as discussed in \S\ref{sec-regularity-of-solutions}.  Thus the mathematical criteria of continuity of $u$ and $\zeta u^\prime$ at $y_j$ correspond to the physical principles of conservation of energy and momentum, in the sense that either starting point leads to (\ref{scattering-at-yj}).  

Consider now the system (\ref{equation},\ref{bc}) for an arbitrary $\zeta\in\step_+(X)$.  Let $y_1,\ldots,y_m$ denote the jump points of $\zeta$, with 
\begin{equation}\label{jump-points}
x_0=y_0<y_1<\cdots<y_m<y_{m+1}=x_1.
\end{equation}
The general solution to (\ref{equation}) has the representation (\ref{general-solution-two-intervals}) for $1\leq j\leq m$, subject to (\ref{scattering-at-yj},\ref{energy-1},\ref{momentum-1},\ref{disk-automorphism-at-yj}).  The boundary condition (\ref{bc}) translated in terms of the constants $A_j,B_j$ and the notation (\ref{muj-rj}) states
\begin{equation}\label{bc-1}
B_0=1\quad\mbox{ and }\quad A_m=\mu_{m+1}\xi B_m\quad\mbox{ where }\quad \xi\in\ddd.
\end{equation}
Invertibility of (\ref{scattering-at-yj}) implies that $u=0$ on the interval $(y_{j-1},y_j)$ if and only if $u=0$ on $(y_j,y_{j+1})$ $(1\leq j\leq m)$. 
The left boundary condition $B_0=1$ is incompatible with the trivial solution; therefore $B_m\neq 0$ in the right boundary condition, which may be equivalently stated as
\begin{equation}\label{right-bc}
A_m/B_m=e^{2i(x_1-y_m)\omega}\xi.
\end{equation}
Writing $g_\zeta^\omega(\xi)=u^{\ell}(x_0+)=A_0=A_0/B_0$, repeated application of (\ref{disk-automorphism-at-yj}) yields
\begin{prop}\label{prop-step-composition}
Fix $X=(x_0,x_1)$.  Let $\zeta\in\step(X)$ have $n\geq 1$ jump points $y_j$, indexed according to their natural order 
\[
x_0<y_1<\cdots<y_n<x_1, 
\]
and write $y_0=x_0$.  Set 
\[
\mu_j=e^{2i(y_j-y_{j-1})}\quad \mbox{ and }\quad
r_j=\frac{\zeta(y_j-)-\zeta(y_j+)}{\zeta(y_j-)+\zeta(y_j+)}\qquad(1\leq j\leq n).
\]
Then
\begin{equation}\label{g-step-composition}
g_\zeta^\omega(\xi)=\varphi_{\mu_1,r_1}\circ\cdots\circ\varphi_{\mu_n,r_n}\bigl(e^{2i(x_1-y_n)\omega}\xi\bigr).
\end{equation}
\end{prop}

Thus $\binom{A_0}{B_0}$ is uniquely determined by (\ref{equation},\ref{bc}), implying by (\ref{scattering-at-yj}) that the full solution $u$ as expressed by  (\ref{general-solution-two-intervals}) is uniquely determined. 

For constant $\zeta$ formula (\ref{muj-rj}) yields $r_j=0$ $(1\leq j\leq m)$, in which case (\ref{g-step-composition}) reduces to 
\begin{equation}\label{constant-zeta}
g_\zeta^\omega(\xi)=e^{2i(x_1-x_0)\omega}\xi.
\end{equation}
It will be seen that $g_\zeta^\omega(\xi)\sim e^{2i(x_1-x_0)\omega}\xi$ as $|\omega|\rightarrow\infty$ for any absolutely continuous $\zeta\in\reg_+(X)$.  Thus the constant case (\ref{constant-zeta}) represents the high frequency asymptotics of $g_\zeta^\omega$ for a general class of impedance functions. 

In summary,
\begin{prop}\label{prop-step-existence}
The system (\ref{equation},\ref{bc}) has a unique solution $u$ if 
 $\zeta\in\step_+(X)$, $\omega\in\real\setminus\{0\}$ and $\xi\in\ddd$.  Setting $g_\zeta^\omega(\xi)=u^{\ell}(x_0+)$ defines a mapping
\begin{equation}\label{g-step}
g:\step_+(X)\rightarrow \left(\aut\mathbb{P}\right)^\real.
\end{equation}
\end{prop}
\begin{pf}
That $g_\zeta^\omega$ is a disk automorphism, and so a non-constant member of $K$, follows from the explicit formulation (\ref{g-step-composition}) as a composition of disk automorphisms. (Note that $\varphi_{e^{2i(x_1-y_n)\omega},0}(\xi)=e^{2i(x_1-y_n)\omega}\xi$ is also a disk automorphism.)
\end{pf}

\subsection{Preliminary technical results\label{sec-technical}}

This section derives a series of technical results concerning $\reg_+(X)$, laying the groundwork for the main results formulated afterward in \S\ref{sec-general-existence}.  The goal of this and the next section is to extend the mapping $g$ of Proposition~\ref{prop-step-existence} to arbitrary regulated $\zeta$ and to prove that it is continuous with respect to the $L^\infty(X)$ norm on $\reg_+(X)$ and convergence pointwise almost everywhere on $K^\real$. 

The first step is to establish a weak equation implied by the original boundary-value problem (\ref{equation},\ref{bc}) based on part~2 of Proposition~\ref{prop-solution-interpretation}, as follows. The definitions (\ref{left-right-components}) of $u^{\ell}$ and $u^{\rm r}$ together with the boundary conditions (\ref{bc}) imply 
\begin{equation}\label{u-prime-at-endpoints}
u^\prime(x_0+)=i\omega\bigl(2-u(x_0+)\bigr)\quad\mbox{ and }\quad u^\prime(x_1-)=i\omega\frac{1-\xi}{1+\xi}u(x_1-).
\end{equation}
Since the linear fractional transformation $\xi\mapsto(1-\xi)/(1+\xi)$ is a bijection from $\ddd$ onto the open right half plane, one may write 
\begin{equation}\label{xi-transformed}
\frac{1-\xi}{1+\xi}=\nu+i\eta\quad\mbox{ where }\quad\nu>0\mbox{ and }\eta\in\real. 
\end{equation}
Given $v\in W^{1,2}(X)$, evaluation of the scalar product $\langle (\zeta u^\prime)^\prime+\omega^2\zeta u,v\rangle_{L^2}$ using integration by parts together with (\ref{u-prime-at-endpoints}) and (\ref{xi-transformed}) transforms (\ref{equation}) into the weak equation 
\begin{equation}\label{weak-equation}
a_1^\omega(u,v)+(1+\omega^2)a_2(u,v)=b_1^\omega(v)\qquad \bigl(v\in W^{1,2}(X)\bigr),
\end{equation}
where binary forms 
\[
a_1^\omega:W^{1,2}(X)\times W^{1,2}(X)\rightarrow\complex,\quad a_2:L^2(X)\times W^{1,2}(X)\rightarrow\complex
\]
and the conjugate linear functional $b_1^\omega:W^{1,2}(X)\rightarrow\complex$ are defined as
\begin{subequations}\label{forms}
\begin{align}
a_1^\omega(u,v)=&-\langle\zeta u,v\rangle_{L^2}-\langle\zeta u^\prime,v^\prime\rangle_{L^2}\label{a1}\\
&+i\omega\bigl(\zeta(x_0+)u(x_0+)\overline{v}(x_0+)+\zeta(x_1-)(\nu+i\eta)u(x_1-)\overline{v}(x_1-)\bigr)\nonumber\\
a_2(u,v)=&\langle\zeta u,v\rangle_{L^2}\label{a2}\\
b_1^\omega(v)=&2i\omega\zeta(x_0+)\overline{v}(x_0+).\label{b1}
\end{align}
\end{subequations}
To simplify notation, dependence of the forms $a_1^\omega,a_2$ and $b_1^\omega$ on $\zeta$ and $\xi$ is suppressed, but this dependency is to be understood.   

From a purely logical standpoint equation (\ref{weak-equation}) is weaker than the system (\ref{equation},\ref{bc}), since the boundary condition is rolled into its formulation.  If a solution $u$ to (\ref{weak-equation}) satisfies (\ref{bc}), then $u$ is also a solution to (\ref{equation}) in the sense of Proposition~\ref{prop-solution-interpretation}, but not necessarily otherwise: one has to prove separately that the boundary condition is satisfied before concluding that a solution to (\ref{weak-equation}) is a solution to (\ref{equation}) and hence to the system (\ref{equation},\ref{bc}).  This takes some work.  It will be proved over the course of this and the next section that (\ref{weak-equation}) and the system (\ref{equation},\ref{bc}) are indeed equivalent for any 
$\zeta\in\reg_+(X)$.

The following elementary estimate implies boundedness of $a_1^\omega$ and $b_1^\omega$.
\begin{prop}\label{prop-boundedness}
Let $X=(x_0,x_1)$ and set 
\[
C_X=\sqrt{2}\max\left\{\frac{1}{\sqrt{x_1-x_0}},2\sqrt{x_1-x_0}\right\}.
\]
If $f\in W^{1,2}(X)$ then for every $x\in X$
\begin{equation}\label{h-infinity}
\left|f(x)\right|\leq C_X\|f\|_{W^{1,2}}.
\end{equation}
\end{prop}

The next result asserts coercivity of $a_1^\omega$ for arbitrary $\omega$ in the closed upper half plane, which is slightly more general than what is actually needed. Subsequent applications only involve the real case $\omega\in\real$. 
\begin{prop}\label{prop-bounded-coercive}
For $\omega=\alpha+i\beta$, where $\alpha\in\real$ and $\beta\geq 0$, the sesquilinear form $a_1^\omega$ is bounded and coercive.  
\end{prop}
\begin{pf}
Denote the limiting values of $\zeta,u$ and $v$ at $x_0+$ by $\zeta_0,u_0$ and $v_0$ respectively; and similarly denote their limiting values at $x_1-$ by $\zeta_1,u_1$ and $v_1$.  
Boundedness of $a_1^\omega$ follows from Proposition~\ref{prop-boundedness}.  In detail,
\[
\begin{split}
|a_1^\omega(u,v)|&=\left|-\langle\zeta u,v\rangle_{L^2}-\langle\zeta u^\prime,v^\prime\rangle_{L^2}+i\omega\bigl(\zeta_0u_0\overline{v}_0+\zeta_1(\nu+i\eta)u_1\overline{v}_1\bigr)\right|\\
&\leq\|\zeta\|_{L^\infty} \|u\|_{\hilbert}\|v\|_{\hilbert}+|\omega|\|\zeta\|_{L^\infty} \left(1+\sqrt{\nu^2+\eta^2}\right)C_X^2\|u\|_{\hilbert}\|v\|_{\hilbert}\\
&=\left(1+|\omega|\left(1+\sqrt{\nu^2+\eta^2}\right)C_X^2\right)\|\zeta\|_{L^\infty} \|u\|_{\hilbert}\|v\|_{\hilbert}.
\end{split}
\]

To verify coercivity, first arrange $a_1^\omega(u,u)$ in terms of its real and imaginary parts as
\begin{equation}\label{real-and-imaginary}
\begin{split}
a_1^\omega(u,u)
=&-\langle\zeta u,u\rangle_{L^2}-\langle\zeta u^\prime,u^\prime\rangle_{L^2}
-\bigl(\beta\zeta_0|u_0|^2+\beta\nu\zeta_1|u_1|^2+\alpha\eta\zeta_1|u_1|^2\bigr)\\
&+
i\bigl(\alpha\zeta_0|u_0|^2+\alpha\nu\zeta_1|u_1|^2-\beta\eta\zeta_1\nu|u_1|^2\bigr).
\end{split}
\end{equation}
If $\left|\Re\,a_1^\omega(u,u)\right|\geq\frac{1}{2}\zeta_{\min}\|u\|^2_{\hilbert}$, then  
$a_1^\omega$ is coercive since $\zeta_{\min}>0$ by (\ref{positive-regulated}).  

Suppose $\left|\Re\,a_1^\omega(u,u)\right|<\frac{1}{2}\zeta_{\min}\|u\|^2_{\hilbert}$.  In this case  (\ref{real-and-imaginary}) implies that 
\begin{equation}\label{real-inequalities}
-\alpha\eta>\beta\nu\geq 0\quad\mbox{ and }\quad -\alpha\eta|u_1|^2>\frac{1}{2}\zeta_{\min}\|u\|^2_{\hilbert},
\end{equation}
since 
\[
\langle\zeta u,u\rangle_{L^2}+\langle\zeta u^\prime,u^\prime\rangle_{L^2}
+\beta\zeta_0|u_0|^2+\beta\nu\zeta_1|u_1|^2
\geq \langle\zeta u,u\rangle_{L^2}+\langle\zeta u^\prime,u^\prime\rangle_{L^2}
\geq\zeta_{\min}\|u\|^2_{\hilbert}.
\]
Therefore $\eta\neq 0$ and 
\begin{equation}\label{real-2-imaginary-inequality}
|\alpha|\,|u_1|^2>\frac{\zeta_{\min}\|u\|^2_{\hilbert}}{2|\eta|}. 
\end{equation}
But then, by (\ref{real-and-imaginary}), (\ref{real-inequalities}) and (\ref{real-2-imaginary-inequality}), the imaginary part of $a_1^\omega$ satisfies
\[
\left|\Im\,a_1^\omega(u,u)\right|\geq|\alpha|\left(\nu+\frac{\beta |\eta|}{|\alpha|}\right)|u_1|^2\geq
\frac{\nu\zeta_{\min}\|u\|^2_{\hilbert}}{2|\eta|}.
\]
Thus in any case, 
\[
|a_1^\omega(u,u)|\geq\textstyle\frac{1}{2}\min\bigl\{1,\frac{\nu}{|\eta|}\bigr\}\zeta_{\min}\|u\|^2_{\hilbert}.
\]

where $\zeta_{\min}>0$ by definition (\ref{positive-regulated}), and $\nu>0$ by (\ref{xi-transformed}).  This proves coercivity. \end{pf}

The functional $b_1^\omega:\hilbert(X)\rightarrow\complex$ is also bounded by Proposition~\ref{prop-boundedness},
\begin{equation}\label{b-bounded}
|b_1^\omega(v)|=|2i\omega\zeta(x_0+)\overline{v}(x_0+)|\leq 2|\omega|\zeta(x_0+)C_X\|v\|_{\hilbert}.
\end{equation}

Proposition~\ref{prop-bounded-coercive} asserts that $a_1^\omega$ satisfies the hypothesis of the Lax-Milgram lemma, and can thus represent any bounded linear or conjugate linear functional on $\hilbert(X)$.  The salient points, which will be applied repeatedly below, are as follows. Let $\ell\in\bigl(\hilbert(X)\bigr)^\ast$ be arbitrary, and let $C$ be a constant such that $|\ell(v)|\leq C\|v\|_{\hilbert}$ for every $v\in\hilbert(X)$.  Let 
\begin{equation}\label{mu}
\mu=\textstyle\frac{1}{2}\min\bigl\{1,\frac{\nu}{|\eta|}\bigr\}\zeta_{\min}
\end{equation}
denote the coercive estimate for $a_1^\omega$ arising in the proof of Proposition~\ref{prop-bounded-coercive}, and suppose $f\in\hilbert(X)$ represents $\overline{\ell}$ via $a_1^\omega$ in the sense that 
\begin{equation}\label{represents}
a_1^\omega(f,v)=\overline{\ell(v)}\qquad(v\in\hilbert(X)).  
\end{equation}
Then 
\begin{equation}\label{representor-bound}
\|f\|_{\hilbert}\leq C/\mu
\end{equation}
since
\[
\mu \|f\|^2_{\hilbert}\leq |a_1^\omega(f,f)|=|\ell(f)|\leq C\|f\|_{\hilbert}.
\]

Treating $a_2(u,v)=\langle\zeta u,v\rangle_{L^2}$ as a sesquilinear form on $L^2(X)\times\hilbert(X)$, note that for a given $\psi\in L^2(X)$, the conjugate linear functional on $\hilbert(X)$ defined by the equation  
\[
\ell(v)=a_2(\psi,v)\qquad(v\in\hilbert(X))
\]
satisfies the bound
\begin{equation}\label{a2-estimate}
|\ell(v)|\leq\|\zeta\|_{L^\infty} \|\psi\|_{L^2}\|v\|_{\hilbert}. 
\end{equation}

For each $\psi\in L^2(X)$, define $A_{\zeta,\omega}\psi\in\hilbert(X)$ using the Lax-Milgram lemma by the equation
\begin{equation}\label{A-defn}
a_1^\omega(A_{\zeta,\omega}\psi,v)=-a_2(\psi,v)\qquad(v\in\hilbert(X)).
\end{equation}
By (\ref{representor-bound}), the estimate (\ref{a2-estimate}) yields the bound 
\begin{equation}\label{A-bound}
\|A_{\zeta,\omega}\psi\|_{\hilbert}\leq\frac{\|\zeta\|_{L^\infty} }{\mu}\|\psi\|_{L^2}. 
\end{equation}
Thus $A_{\zeta,\omega}:L^2(X)\rightarrow\hilbert(X)$ is a bounded linear operator.  And  $A_{\zeta,\omega}:L^2(X)\rightarrow L^2(X)$ is compact, because $\hilbert(X)$ is compactly embedded in $L^2(X)$ by Rellich's theorem. 

The point of introducing $A_{\zeta,\omega}$ is to represent $a_1^\omega(u,v)+(1+\omega^2)a_2(u,v)$ in terms of a compact perturbation of the identity. By definition, 
\begin{equation}\label{compact-perturbation}
a_1^\omega\left((I-(1+\omega^2)A_{\zeta,\omega})u,v\right)=a_1^\omega(u,v)+(1+\omega^2)a_2(u,v). 
\end{equation}
Define $\rho\in\hilbert(X)$ using the Lax-Milgram lemma by the equation 
\begin{equation}\label{lax-milgram-rho}
a_1^\omega(\rho,v)=b_1^\omega(v)\qquad(v\in\hilbert(X)), 
\end{equation}
so that by (\ref{b-bounded}),  
\begin{equation}\label{rho-estimate}
\|\rho\|_{\hilbert}\leq\frac{2|\omega|\zeta_0C_X}{\mu},
\end{equation}
as in (\ref{representor-bound}).  Then, for $u\in\hilbert(X)$, 
\begin{subequations}\label{operator-reduction}
\begin{gather}
a_1^\omega(u,v)+(1+\omega^2)a_2(u,v)=b_1^\omega(v)\qquad(v\in\hilbert(X))\label{weak-for-all}\\
\Leftrightarrow\quad a_1^\omega\left((I-(1+\omega^2)A_{\zeta,\omega})u,v\right)=a_1^\omega(\rho,v)\qquad(v\in\hilbert(X))\label{a1-operator-equation}\\
\Leftrightarrow\quad (I-(1+\omega^2)A_{\zeta,\omega})u=\rho\label{operator-equation},
\end{gather}
\end{subequations}
the last equivalence by the Lax-Milgram lemma.  Observe that if (\ref{operator-equation}) has a solution $u\in L^2(X)$, then automatically
\begin{equation}\label{u-in-hilbert}
u=(1+\omega^2)A_{\zeta,\omega}u+\rho\in\hilbert(X). 
\end{equation}
This proves
\begin{prop}\label{prop-operator-equivalence}
Existence of a solution $u\in L^2(X)$ to (\ref{operator-equation}) is equivalent to the existence of a solution $u\in\hilbert(X)$ to (\ref{weak-for-all}).  
\end{prop} 
The question of injectivity of $I-(1+\omega^2)A_{\zeta,\omega}$ turns out to involve the spectrum of $A_{\zeta,0}$, which is easily seen to be real.  Observe furthermore that $1\in\spec A_{\zeta,0}$, as follows.  For $\omega=0$, definitions (\ref{forms}) and (\ref{A-defn}) imply
\[
A_{\zeta,0}\psi=\psi\quad\Leftrightarrow\quad \langle\zeta\psi^\prime,v^\prime\rangle_{L^2}=0\qquad(v\in\hilbert(X)),
\]
which is satisfied by any constant function $\psi$.  Any non-zero constant $\psi$ is therefore an eigenvector of $A_{\zeta,0}$ corresponding to eigenvalue $1\in\spec A_{\zeta,0}$.  

Set 
\begin{equation}\label{sigma-zeta}
\Sigma_\zeta=\left\{\left. \pm\sqrt{\lambda^{-1}-1}\;\right|\;\lambda\in\spec A_{\zeta,0}\cap(0,1]\right\}.
\end{equation}
Since $1\in\spec A_{\zeta,0}$, it follows from this definition that $0\in\Sigma_{\zeta}$ for any $\zeta$.   
Compactness of $A_{\zeta,0}$ implies that the set $\Sigma_\zeta$ is at most countable, with $\pm\infty$ the only possible accumulation points. 
It can be shown using Wirtinger's inequality that if $\omega\neq0$ and $\omega\in\Sigma_\zeta$, then
\begin{equation}\label{wirtinger}
|\omega|\geq\frac{\pi\sqrt{\zeta_{\min}}}{(x_1-x_0)\sqrt{\|\zeta\|_{L^\infty} }}. 
\end{equation}

\begin{prop}\label{prop-injectivity}
Let $\zeta\in\reg_+(X)$ and $\omega\in\real\setminus\Sigma_\zeta$.  
If $\psi\in\hilbert(X)$ and 
\begin{equation}\label{lhs-0}
a_1^\omega(\psi,v)+(1+\omega^2)a_2(\psi,v)=0
\end{equation}
for every $v\in\hilbert(X)$, then $\psi=0$.  
\end{prop}
\begin{pf}
Abbreviate $\zeta(x_0+)$ and $\psi(x_0+)$ to $\zeta_0$ and $\psi_0$ respectively, and similarly write $\zeta_1$ and $\psi_1$ for $\zeta(x_1-)$ and $\psi(x_1-)$. 
Setting $v=\psi$, the left-hand side of (\ref{lhs-0}) becomes
\begin{equation}\label{lhs-1}
\omega^2\langle \zeta\psi,\psi\rangle_{L^2}-\langle\zeta \psi^\prime,\psi^\prime\rangle_{L^2}-\omega\zeta_1\eta|\psi_1|^2+i\omega\bigl(\zeta_0|\psi_0|^2+\zeta_1\nu|\psi_1|^2\bigr).
\end{equation}
Recall from (\ref{xi-transformed}) that $\nu>0$; also, $\omega\neq0\in\Sigma_\zeta$.  
Thus the imaginary part of (\ref{lhs-1}) being zero implies that $\psi_0=\psi_1=0$, in which case the left-hand side of (\ref{lhs-0}) is
\[
\omega^2\langle \zeta v,\psi\rangle_{L^2}-\langle\zeta\psi^\prime,v^\prime\rangle_{L^2}=a_1^0\left((I-(1+\omega^2)A_{\zeta,0})\psi,v\right),
\]
where $a_1^0$ in the last expression is defined according to (\ref{a1}).  Since (\ref{lhs-0}) holds for every $v\in\hilbert(X)$, it follows from Proposition~\ref{prop-bounded-coercive} and the uniqueness part of the Lax-Milgram lemma that 
\[
(I-(1+\omega^2)A_{\zeta,0})\psi=0,
\]
which implies $\omega\in\Sigma_\zeta$, unless $\psi=0$. 
\end{pf}

\begin{prop}\label{prop-weak-estimate}
Let $\zeta\in\reg_+(X)$ and $\omega\in\real\setminus\Sigma_\zeta$.
Then the operator
\[
I-(1+\omega^2)A_{\zeta,\omega}:L^2(X)\rightarrow L^2(X)
\]
is a bijection, and the equation (\ref{weak-equation}),
\begin{equation}\label{the-weak-equation}
a_1^\omega(\psi,v)+(1+\omega^2)a_2(\psi,v)=b_1^\omega(v)\qquad(v\in\hilbert(X)),
\end{equation}
has a unique solution $\psi\in\hilbert(X)$. 
\end{prop}

\begin{pf}  By Proposition~\ref{prop-operator-equivalence} it suffices to consider the equation 
\[
(I-(1+\omega^2)A_{\zeta,\omega})\psi=\rho.
\]
Fredholm theory implies that the compact perturbation of the identity 
\[
I-(1+\omega^2)A_{\zeta,\omega}:L^2(X)\rightarrow L^2(X)
\]
is surjective if and only if it is injective.  To see that it is injective, suppose 
\[
(I-(1+\omega^2)A_{\zeta,\omega})\psi=0
\]
for some $\psi\in L^2(X)$.  Then 
\[
\psi=(I+\omega^2)A_{\zeta,\omega}\psi\in\hilbert(X),
\]
and for every $v\in\hilbert(X)$, 
\[
a_1^\omega((I-(1+\omega^2)A_{\zeta,\omega})\psi,v)=0,
\]
which implies $\psi=0$ by Proposition~\ref{prop-injectivity}. Thus 
$
I-(1+\omega^2)A_{\zeta,\omega}
$
is a bijection on $L^2(X)$, and $\psi=(I-(1+\omega^2)A_{\zeta,\omega})^{-1}\rho$ is uniquely determined.  Moreover, 
\[
\psi=(1+\omega^2)A_{\zeta,\omega}\psi+\rho\in\hilbert(X),
\]
completing the proof. \end{pf}

The parameters $\xi$ and $\omega\neq0$ being fixed, Proposition~\ref{prop-weak-estimate} determines an operator
\begin{equation}\label{solution-operator}
\left\{\zeta\in\reg_+(X)\,\left|\,\omega\not\in\Sigma_\zeta\right.\right\}\rightarrow\hilbert(X),\qquad\zeta\mapsto\psi_\zeta,
\end{equation}
where $\psi_\zeta$ denotes the unique solution to (\ref{the-weak-equation}).
For the remainder of the present section $\psi_\zeta$ will denote this unique solution.   

\begin{prop}\label{prop-proto-stability}
Let $\zeta,\widetilde\zeta\in\reg_+(X)$ and $\omega\in\real\setminus\bigl(\Sigma_\zeta\cup\Sigma_{\widetilde\zeta}\bigr)$. There exists a positive constant $F_{\omega,\zeta}$, independent of $\widetilde\zeta$, such that 
\begin{equation}
\|\psi_\zeta-\psi_{\widetilde\zeta}\|_{\hilbert}\leq F_{\omega,\zeta}\|\psi_{\widetilde\zeta} \|_{\hilbert}\, \|\zeta-\widetilde\zeta\|_{L^\infty}.
\end{equation}
\end{prop}
\begin{pf}
To reduce notational clutter, write $\psi=\psi_\zeta$ and $\widetilde\psi=\psi_{\widetilde\zeta}$.  
Let $\widetilde{a}_1^\omega, \widetilde{a}_2$ denote the forms defined according to (\ref{forms}) with $\widetilde\zeta$ in place of $\zeta$. Set
\[
a^\omega=a_1^\omega+(1+\omega^2)a_2,\qquad  \widetilde{a}^\omega=\widetilde{a}_1^\omega+(1+\omega^2)\widetilde{a}_2
\]
 and
\[
a_3^\omega(u,v)=\omega^2\langle(\zeta-\widetilde\zeta)u,v\rangle_{L^2}-\langle(\zeta-\widetilde\zeta)u^\prime,v^\prime\rangle_{L^2}.
\]
Set $\theta=\psi-\widetilde\psi$ and observe that for $v\in\hilbert(X)$, 
\begin{equation}\label{b2-derivation}
\begin{split}
a^\omega(\theta,v)&=a^\omega(\psi,v)-a^\omega(\widetilde\psi,v)\\
&=b_1^\omega(v)-\left(\widetilde{a}^\omega(\widetilde\psi,v)+a_3^\omega(\widetilde\psi,v)\right)\\
&=b_1^\omega(v)-b_1^\omega(v)-a_3^\omega(\widetilde\psi,v)\\
&=-\omega^2\langle(\zeta-\widetilde\zeta)\widetilde\psi,v\rangle_{L^2}+\langle(\zeta-\widetilde\zeta)\widetilde\psi^\prime,v^\prime\rangle_{L^2}. 
\end{split}
\end{equation}
Define $b_2^\omega(v)$ accordingly as 
\begin{equation}\label{b2-definition}
b_2^\omega(v)=-a_3^\omega(\widetilde\psi,v)=-\omega^2\langle(\zeta-\widetilde\zeta)\widetilde\psi,v\rangle_{L^2}+\langle(\zeta-\widetilde\zeta)\widetilde\psi^\prime,v^\prime\rangle_{L^2},
\end{equation}
so that $\theta$ satisfies the equation 
\begin{equation}\label{theta-equation}
a^\omega(\theta,v)=b_2^\omega(v)\qquad(v\in\hilbert(X)),
\end{equation}
and observe that 
\begin{equation}\label{b2-bound}
|b_2^\omega(v)|\leq (1+\omega^2)\|\widetilde\psi\|_{\hilbert}\, \|\zeta-\widetilde\zeta\|_{L^\infty}\|v\|_{\hilbert}.
\end{equation}
Define $\sigma$ using the Lax-Milgram lemma by the equation 
\begin{equation}\label{sigma-definition}
a_1^\omega(\sigma,v)=b_2^\omega(v)\qquad(v\in\hilbert(X)),
\end{equation}
so that, as per (\ref{representor-bound}), 
\begin{equation}\label{sigma-bound}
\|\sigma\|_{\hilbert}\leq \frac{(1+\omega^2)\|\widetilde\psi\|_{\hilbert}\, \|\zeta-\widetilde\zeta\|_{L^\infty}}{\mu}.
\end{equation}
By Proposition~\ref{prop-weak-estimate}, the operator
\[
I-(1+\omega^2)A_{\zeta,\omega}:L^2(X)\rightarrow L^2(X)
\]
is a bijection.  Any bounded bijection on a Banach space has a bounded inverse, so there exists a positive constant $C_{\omega,\zeta}>0$  such that 
\[
\|(I-(1+\omega^2)A_{\zeta,\omega})^{-1}f\|_{L^2}\leq C_{\omega,\zeta} \|f\|_{L^2}\qquad(f\in L^2(X)). 
\]
Since $\theta$ satisfies the equation (\ref{theta-equation}), it follows by definition of $\sigma$ that 
\begin{equation}\label{sigma-satisfies}
(I-(1+\omega^2)A_{\zeta,\omega})\theta=\sigma,
\end{equation}
whereby
\begin{equation}\label{sigma-estimate}
\|\theta\|_{L^2}\leq C_{\omega,\zeta}\|\sigma\|_{L^2}\leq C_{\omega,\zeta}\|\sigma\|_{\hilbert}. 
\end{equation}
Rewriting (\ref{sigma-satisfies}) as $\theta=\sigma+(1+\omega^2)A_{\zeta,\omega}\theta$, and using the estimate (\ref{A-bound}), yields
\begin{equation}\label{theta-estimate}
\begin{split}
\|\theta\|_{\hilbert}&\leq \|\sigma\|_{\hilbert}+\frac{(1+\omega^2)\|\zeta\|_{L^\infty} }{\mu}\|\theta\|_{L^2}\\
&\leq \|\sigma\|_{\hilbert}+\frac{(1+\omega^2)\|\zeta\|_{L^\infty} }{\mu}C_{\omega,\zeta}\|\sigma\|_{\hilbert}\\
&\leq \left(1+\frac{(1+\omega^2)\|\zeta\|_{L^\infty} }{\mu}C_{\omega,\zeta}\right)\frac{(1+\omega^2)\|\widetilde\psi\|_{\hilbert}\, \|\zeta-\widetilde\zeta\|_{L^\infty}}{\mu}.
\end{split}
\end{equation}
Setting 
\[
F_{\omega,\zeta}= \left(1+\frac{(1+\omega^2)\|\zeta\|_{L^\infty} }{\mu}C_{\omega,\zeta}\right)\frac{(1+\omega^2)}{\mu},
\]
one therefore has
\[
\|\psi-\widetilde\psi\|_{\hilbert}\leq F_{\omega,\zeta}\|\widetilde\psi\|_{\hilbert}\, \|\zeta-\widetilde\zeta\|_{L^\infty}
\]
as desired. 
\end{pf}

\begin{lem}[First singular approximation lemma.]\label{lem-forward-stability}
Let $\zeta,\zeta_n\in\reg_+(X)$ $(n\geq 1)$, set 
$
\Sigma=\Sigma_\zeta\cup\bigcup_{n=1}^\infty\Sigma_{\zeta_n},
$
and fix $\omega\in\real\setminus\Sigma$.  
If $\zeta_n\rightarrow\zeta$ uniformly, then 
$
\psi_{\zeta_n}\rightarrow\psi_\zeta
$
in $\hilbert(X)$. 
If in addition each $\psi_{\zeta_n}$ satisfies the boundary conditions (\ref{bc}), then $\psi_\zeta$ also satisfies (\ref{bc}), and moreover
\begin{equation}\label{component-convergence}
\psi_{\zeta_n}^\ell\rightarrow\psi_{\zeta}^\ell\quad\mbox{ and }\quad \psi_{\zeta_n}^{\it{r}}\rightarrow\psi_{\zeta}^{\it{r}}
\end{equation}
uniformly. 
\end{lem}
\begin{pf}
Write $\psi=\psi_\zeta$ and $\psi_n=\psi_{\zeta_n}$.  
First note that the sequence $\|\psi_n\|_{\hilbert}$ $(n\geq1)$ is bounded.  Otherwise there is a subsequence $\|\psi_{n_j}\|_{\hilbert}\rightarrow\infty$, in which case 
\begin{equation}\label{psi-blowup}
\frac{\|\psi-\psi_{n_j}\|_{\hilbert}}{\|\psi_{n_j}\|_{\hilbert}}\rightarrow 1\quad\mbox{ as }\quad j\rightarrow\infty. 
\end{equation}
But Proposition~\ref{prop-proto-stability} implies that 
\[
\frac{\|\psi-\psi_{n_j}\|_{\hilbert}}{\|\psi_{n_j}\|_{\hilbert}}\leq F_{\omega,\zeta}\|\zeta-\zeta_{n_j}\|_{L^\infty}\rightarrow 0\quad\mbox{ as }\quad j\rightarrow\infty,
\]
contradicting (\ref{psi-blowup}).  Letting $B=\sup_{n\geq 1}\|\psi_n\|_{\hilbert}<\infty$, Proposition~\ref{prop-proto-stability} implies in turn that 
\[
\|\psi-\psi_n\|_{\hilbert}\leq F_{\omega,\zeta}B\|\zeta-\zeta_n\|_{L^\infty}\rightarrow 0\quad\mbox{ as }\quad n\rightarrow\infty.
\]
It follows from Proposition~\ref{prop-boundedness} that $\psi_n\rightarrow\psi$ uniformly. Therefore $\zeta_n\psi_n\rightarrow\zeta\psi$ uniformly, since $\zeta_n\rightarrow\zeta$ uniformly. 

Suppose that each $\psi_n$ is a solution to (\ref{the-weak-equation},\ref{bc}).  Then, since $\zeta_n\psi_n$ converges uniformly, 
\[
M=\sup_{x\in X}\omega^2\zeta_n(x)|\psi_n(x)|<\infty.
\]
Recall that every solution to (\ref{the-weak-equation},\ref{bc}) satisfies (\ref{equation}). 
It follows from (\ref{equation}) that for every $x,y\in X$,
\[
|\zeta_n(y)\psi_{n\,}^\prime(y)-\zeta_n(x)\psi_{n\,}^\prime(x)|=\left|\int_x^y\left(\zeta_n(s)\psi_{n\,}^\prime(s)\right)^\prime\,ds\right|=\left|\int_x^y\omega^2\zeta_n(s)\psi_n(s)\,ds\right|\leq |y-x|M,
\]
whereby the functions $\zeta_n\psi_{n\,}^\prime$ are equicontinuous. By the Arzel\`a-Ascoli theorem, a subsequence $\zeta_{n_j}\psi_{n_j}^\prime$ of these functions converges uniformly.  Since the coefficients $\zeta_n$ converge uniformly, it follows that $\psi_{n_j}^\prime$ converges uniformly.  Since the functions $\psi_n$ themselves converge to $\psi$, it follows in turn that 
$
\psi_{n_j}^\prime\rightarrow\psi^\prime
$
uniformly. By similiar reasoning, uniform convergence $\psi_n\rightarrow\psi$ further implies uniform convergence of the full sequence
\begin{equation}\label{uniform-derivative-convergence}
\psi_{n\,}^\prime\rightarrow\psi^\prime.
\end{equation}
In particular, the boundary conditions (\ref{bc}) hold for $\psi$, since they hold for each $\psi_n$, from which it follows that $\psi$ is a solution to (\ref{equation},\ref{bc}).   Uniform convergence of the left- and right-moving components (\ref{component-convergence}) follows from that of $\psi_n$ together with (\ref{uniform-derivative-convergence}).  
\end{pf}

\subsection{Existence and uniqueness for $\zeta\in\reg_+(X)$\label{sec-general-existence}}

The main results concerning existence and uniqueness of solutions to the system (\ref{equation},\ref{bc}) for $\zeta\in\reg_+(X)$ follow from the technical results established in the previous section in combination with \S\ref{sec-existence-step}, as does the fact that, viewed as a function of the boundary parameter $\xi$, $u^{\ell}(x_0+)$ is a disk automorphism. 

Recall the definition (\ref{sigma-zeta}) of the countable set $\Sigma_\zeta$ of exceptional frequencies $\omega$.  
\begin{thm}\label{thm-general-existence}
Let $\zeta\in\reg_+(X)$, $\omega\in\real\setminus\Sigma_\zeta$ and $\xi\in\ddd$.  Then the system (\ref{equation},\ref{bc}) has a unique solution $u$ in the sense of Proposition~\ref{prop-solution-interpretation}.  
\end{thm}
\begin{pf}
Because $\zeta\in\reg_+(X)$ there exists a sequence of step functions $\zeta_n\in\step_+(X)$ $(n\geq 1)$ uniformly convergent to $\zeta$.  Since $0\in\Sigma_\zeta$, Proposition~\ref{prop-step-existence} implies that the system (\ref{equation},\ref{bc}), with $\zeta_n$ in place of $\zeta$, has a unique solution $u_n$.  Moreover, $u_n=\psi_{\zeta_n}$ by Proposition~\ref{prop-weak-estimate} and the fact that (\ref{the-weak-equation}) is implied by (\ref{equation},\ref{bc}).  Therefore $u_n\rightarrow\psi_\zeta$ by Lemma~\ref{lem-forward-stability}, and $u:=\psi_\zeta$ satisfies the boundary condition (\ref{bc}), since each $u_n$ does. It follows that $u$ satisfies the system (\ref{equation},\ref{bc}) in the sense of Proposition~\ref{prop-solution-interpretation}. It is the unique such solution by uniqueness of $\psi_\zeta$ as a solution to (\ref{the-weak-equation}).  
\end{pf}

Theorem~\ref{thm-general-existence} legitimizes the definition of $g$ in \S\ref{sec-generalized-reflection-coefficient-2}: for any $\zeta\in\reg_+(X)$, $\omega\in\real\setminus\Sigma_\zeta$ and $\xi\in\ddd$, let $u$ denote the corresponding unique solution to (\ref{equation},\ref{bc}), and set
\begin{equation}\label{general-g-definition}
g_\zeta^\omega(\xi)=u^{\ell }(x_0+).  
\end{equation}
Using Lemma~\ref{lem-forward-stability} the nature of the dependency of $g_\zeta^\omega(\xi)$ on $\xi$ can be inferred from that of $g_{\zeta_n}^\omega(\xi)$ for step functions $\zeta_n$ approximating $\zeta$, as follows. 
\begin{thm}\label{thm-g-in-K}
If $\zeta\in\reg_+(X)$ and $\omega\in\real\setminus\Sigma_\zeta$ then $g_\zeta^\omega\in K$. If in addition $\zeta$ has bounded variation then $g_\zeta^\omega\in\aut\mathbb{P}$ is non-constant.  
\end{thm}
\begin{pf}
Since $\zeta\in\reg_+(X)$ there exists a sequence $\zeta_n\in\step_+(X)$ such that $\zeta_n\rightarrow\zeta$, and with the following property for each $n\geq 1$.  Letting $y_1,\ldots,y_m$ denote the jump points of $\zeta_n$ in their natural order,
\[
x_0=y_0<y_1<\cdots<y_m<y_{m+1}=x_1,
\]
and writing $Y_j$ for the closure of the interval $(y_j,y_{j+1})$ in $X$ $(0\leq j\leq m)$, it may be assumed without loss of generality that 
\begin{equation}\label{bounded-variation-property}
\zeta_n(y_j+)=\zeta_n(y_{j+1}-)\in\zeta(Y_j)\qquad(0\leq j\leq m).
\end{equation}
For, given $\tilde{\zeta}_n\in\step_+(X)$, there exists $\zeta_n\in\step_+(X)$ satisfying (\ref{bounded-variation-property}) such that 
\[
\|\zeta_n-\zeta\|_{L^\infty}\leq 2\|\tilde{\zeta}_n-\zeta\|_{L^\infty}.
\]
By the results in \S\ref{sec-existence-step}, $g_{\zeta_n}^\omega\in\aut\mathbb{P}$ for each $n\geq1$.  The sequence $g_{\zeta_n}^\omega$ has a limit point in $f\in K$ since $K$ is compact. Pointwise convergence of $g_{\zeta_n}^\omega$ to $g_\zeta^\omega$ with respect to $\xi\in\ddd$ as implied by Lemma~\ref{lem-forward-stability} forces $f=g_\zeta^\omega$.  

Suppose $\zeta\in\bv(X)$ and denote the total variation of $\zeta$ by $\mathscr{V}$.  By definition of $\reg_+(X)$ there exists $\varepsilon>0$ such that $\zeta(x)>\varepsilon$ for every $x\in X$.  Set 
\begin{equation}\label{M-bound-on-approximants}
M=\tanh\left(\frac{\mathscr{V}}{2\varepsilon}\tanh^{-1}\frac{\|\zeta\|_{L^\infty}-\varepsilon}{\|\zeta\|_{L^\infty}+\varepsilon}\right),
\end{equation}
so that $0\leq M<1$.  It is claimed that $\bigl|g_{\zeta_n}^\omega(0)\bigr|\leq M$ for every $n\geq 1$.  Fix $n$, and let $y_1<\cdots<y_{m}$ denote the jump points of $\zeta_n$ as above. Then $g_{\zeta_n}^\omega(0)$ has the form
\begin{equation}\label{g-sequence-form}
g_{\zeta_n}^\omega(0)=\varphi_{\mu_1,r_1}\circ\cdots\circ\varphi_{\mu_{m},r_{m}}(0)
\end{equation}
where 
\begin{equation}\label{mu-r-g-sequence-form}
\mu_k=e^{2i(y_j-y_{j-1})\omega}\quad\mbox{ and }\quad r_j=\frac{\zeta(y_j-)-\zeta(y_j+)}{\zeta(y_j-)+\zeta(y_j+)}\qquad(1\leq j\leq m).
\end{equation}
Note that  
\begin{equation}\label{rj-bound-00}
|r_j|\leq\left|\zeta(y_j-)-\zeta(y_j+)\right|/(2\varepsilon)\quad (1\leq j\leq m)\quad\mbox{ and }\quad\sum_{j=1}^m\left|\zeta(y_j-)-\zeta(y_j+)\right|\leq\mathscr{V},
\end{equation}
the latter inequality by (\ref{bounded-variation-property}).  Also,
\begin{equation}\label{rj-bound-01}
|r_j|\leq\frac{\|\zeta\|_{L^\infty}-\varepsilon}{\|\zeta\|_{L^\infty}+\varepsilon}\quad\mbox{ and }\quad\tanh^{-1}|r_j|\leq|r_j|\tanh^{-1}\frac{\|\zeta\|_{L^\infty}-\varepsilon}{\|\zeta\|_{L^\infty}+\varepsilon}\qquad (1\leq j\leq m).
\end{equation}
Observe that for $\rho\in\ddd$ and $0\leq \delta<1$, 
\[
\max_{|v|\leq\delta}\left|\frac{v+\rho}{1+\overline{\rho}v}\right|=\frac{\delta+|\rho|}{1+\delta|\rho|}=\tanh\left(\tanh^{-1}\delta+\tanh^{-1}|\rho|\right).
\]
Applying this observation iteratively to the right-hand side of (\ref{g-sequence-form}) from right to left yields 
\[
\begin{split}
\left|g_{\zeta_n}^\omega(0)\right|\leq&\tanh\left(\sum_{j=1}^{m}\tanh^{-1}|r_j|\right)\\
\leq&\tanh\left(\sum_{j=1}^{m}|r_j|\tanh^{-1}\frac{\|\zeta\|_{L^\infty}-\varepsilon}{\|\zeta\|_{L^\infty}+\varepsilon}\right)\quad\mbox{by (\ref{rj-bound-01})}\\
\leq&\tanh\left(\sum_{j=1}^{m}\frac{\left|\zeta(y_j-)-\zeta(y_j+)\right|}{2\varepsilon}\tanh^{-1}\frac{\|\zeta\|_{L^\infty}-\varepsilon}{\|\zeta\|_{L^\infty}+\varepsilon}\right)\quad\mbox{by (\ref{rj-bound-00})}\\
\leq&\tanh\left(\frac{\mathscr{V}}{2\varepsilon}\tanh^{-1}\frac{\|\zeta\|_{L^\infty}-\varepsilon}{\|\zeta\|_{L^\infty}+\varepsilon}\right)= M,
\end{split}
\]
proving the claim.  

It follows that $g_{\zeta_n}^\omega$ and its inverse have the form 
\[
g_{\zeta_n}^\omega=\varphi_{\nu_n,\rho_n}\quad\mbox{ and }\quad (g_{\zeta_n}^\omega)^{-1}=\varphi_{\bar\nu_n,-\nu_n\rho_n},
\]
where $|\rho_n|=\bigl|g_{\zeta_n}^\omega(0)\bigr|=\bigl|\bigl( g_{\zeta_n}^\omega\bigr)^{-1}(0)\bigr|\leq M$ for each $n\geq 1$.  Each $g_{\zeta_n}^\omega$ and its inverse extend to holomorphic bijections of the closed disk $\overline{\ddd}$.  The bound $|\rho_n|\leq M<1$ implies that the sequence $g_{\zeta_n}^\omega$ $(n\geq 1)$ is equicontinuous, as is the corresponding sequence of inverses.  The Arzel\'a-Ascoli theorem thus implies the existence of subsequences $g_{\zeta_{n_j}}^\omega$ and $(g_{\zeta_{n_j}}^\omega)^{-1}$ $(j\geq 1)$ uniformly convergent on $\overline{\ddd}$.  Setting 
\[
f=\lim_{j\rightarrow\infty}g_{\zeta_{n_j}}^\omega\quad\mbox{ and }\quad k=\lim_{j\rightarrow\infty}(g_{\zeta_{n_j}}^\omega)^{-1},
\]
it follows that $f$ is holomorphic on $\ddd$ and $k=f^{-1}$, whereby $f$ is bijective.  As before, pointwise convergence of $g_{\zeta_n}^\omega$ to $g_\zeta^\omega$ with respect to $\xi\in\ddd$ forces $f=g_\zeta^\omega$.  Thus $g_\zeta^\omega$ is a disk automorphism.  
\end{pf}

Theorem~\ref{thm-g-in-K} completes the proof that the mapping $g:\reg_+(X)\rightarrow K^\real$ is well-defined. There exist regulated functions of unbounded variation, such as $x^2\sin(1/x^2)$ on $X=(0,1)$; the possibility of $g_\zeta^\omega\in K\setminus\aut\mathbb{P}$ being constant is restricted to such functions. 
Lemma~\ref{lem-forward-stability} also implies: 
\begin{thm}\label{thm-continuous}
The mapping $g:\reg_+(X)\rightarrow K^\real$ is continuous.  
\end{thm}
\begin{pf}  
Fix $\zeta, \zeta_n\in\reg_+(X)$ $(n\geq 1)$ such that $\zeta_n\rightarrow\zeta$ uniformly, and set 
\[
\Sigma=\Sigma_\zeta\cup\bigcup_{n=1}^\infty\Sigma_{\zeta_n}.
\]
Note that $\Sigma$ is countable, and let $\omega\in\real\setminus\Sigma$.  It suffices to show that $g_{\zeta_n}^\omega\rightarrow g_\zeta^\omega$ uniformly on compact subsets of $\ddd$.  Lemma~\ref{lem-forward-stability} implies that $g_{\zeta_n}^\omega(\xi)\rightarrow g_\zeta^\omega(\xi)$ for any given $\xi\in\ddd$, and in particular for $\xi=0$ and $\xi=1/2$.  Since $g_\zeta^\omega, g_{\zeta_n}^\omega\in K$ by Theorem~\ref{thm-g-in-K}, the desired result follows from the fact that a disk automorphism is determined by its value at these two points.  In detail, set 
\[
a_n=g_{\zeta_n}^\omega(0),\quad a=g_\zeta^\omega(0),\quad b_n=g_{\zeta_n}^\omega(1/2),\quad b=g_\zeta^\omega(1/2).
\]
For each $n$, either $g_{\zeta_n}^\omega=\mathbf{c}_{a_n}$ where $a_n\in S^1$ if $g_{\zeta_n}^\omega$ is constant, or 
\begin{equation}\label{two-values-n}
g_{\zeta_n}^\omega=\varphi_{\mu_n,\rho_n}\quad\mbox{ where }\quad \rho_n=\frac{a_n(1-\bar{a}_nb_n)}{2(b_n-a_n)}\mbox{ and }\mu_n=a_n/\rho_n,
\end{equation}
and similarly either $g_\zeta^\omega$ is constant or
\begin{equation}\label{two-values-limit}
g_{\zeta}^\omega=\varphi_{\mu,\rho}\quad\mbox{ where }\quad \rho=\frac{a(1-\bar{a}b)}{2(b-a)}\mbox{ and }\mu=a/\rho.
\end{equation}
If $g_\zeta^\omega=\mathbf{c}_\sigma$ has constant value $\sigma\in S^1$, then $a_n\rightarrow\sigma$ and $\frac{d}{d\xi}g_{\zeta_n}^\omega$ converges uniformly to 0 on compact subsets of $\ddd$, since either $\frac{d}{d\xi}g_{\zeta_n}^\omega(\xi)=0$ if $g_{\zeta_n}^\omega$ is constant, or
\[
\frac{d}{d\xi}g_{\zeta_n}^\omega(\xi)=\mu_n\frac{1-|\rho_n|^2}{(1+\bar{\rho}_n\xi)^2}
\]
and $|\rho_n|=|a_n|$. It follows that $g_{\zeta_n}^\omega\rightarrow \mathbf{c}_\sigma$ uniformly on compact subsets of $\ddd$.  On the other hand, if $g_\zeta^\omega=\varphi_{\mu,\rho}$ is non-constant, then, since $a_n\rightarrow a$ and $b_n\rightarrow b$, it follows from (\ref{two-values-n}) and (\ref{two-values-limit}) that $g_{\zeta_n}^\omega\rightarrow g_{\zeta}^\omega$ uniformly on the closed disk $\overline{\ddd}$. 
\end{pf}

Theorem~\ref{thm-continuous} in conjunction with density of $\step_+(X)$ in $\reg_+(X)$ yields the following result concerning concatenation. 
\begin{prop}\label{prop-concatenation}
Let $x_0<x_1<x_2$, set $X_1=(x_0,x_1)$, $X_2=(x_1,x_2)$ and $X=(x_0,x_2)$.  Let $\zeta\in\reg_+(X)$ be continuous at $x_1$, and define $\widetilde{\zeta}\in\step_+(X)$ to have a single jump point at $x_1$,
\[
\widetilde{\zeta}(x)=1+cH(x-x_1)\qquad (x\in X),
\]
where $c>-1$ and $H$ denotes the Heaviside function.  Set $\zeta_1=\zeta|_{X_1}$ and $\zeta_2=\zeta|_{X_2}$, so that $\zeta=\zeta_1^{\;\smallfrown}\zeta_2$.  Set 
\[
r=\frac{\widetilde{\zeta}(x_1-)-\widetilde{\zeta}(x_1+)}{\widetilde{\zeta}(x_1-)+\widetilde{\zeta}(x_1+)}=c/(c+2).
\]
Then 
\begin{equation}\label{concatenation-with-jump}
g_{\zeta}=g_{\zeta_1}\circ g_{\zeta_2}\quad\mbox{ and }\quad g_{\zeta\widetilde{\zeta}}=g_{\zeta_1}\circ\varphi_{1,r}\circ g_{\zeta_2}. 
\end{equation}
\end{prop}
\begin{pf}
The key formula is (\ref{g-step-composition}), which in combination with Theorem~\ref{thm-continuous} easily gives the desired result.  Fix $\omega\in\real\setminus\Sigma_{\zeta\widetilde{\zeta}}$.  Let $\zeta_1^n\in\step_+(X_1)$, $\zeta_2^n\in\step_+(X_2)$ $(n\geq 1)$ be sequences converging uniformly to $\zeta_1$ and $\zeta_2$ respectively.  For fixed $n$, list the jump points of $\bigl(\zeta_1^{n\;\smallfrown}\zeta_2^n\bigr)\widetilde{\zeta}$ as 
\[
y_1<\cdots<y_p<x_1<\eta_1<\cdots<\eta_q,
\]
$\zeta_1^n$ and $\zeta_2^n$ having $p$ and $q$ jump points respectively.  For $1\leq j\leq p$, define $\mu_j$ and $r_j$ in terms of $\zeta_1^n$ and the jump points $y_j$ according to (\ref{muj-rj}), and analogously set 
\[
\nu_j=e^{2i(\eta_j-\eta_{j-1})\omega}\quad\mbox{ and }\quad s_j=\frac{\zeta_2^n(\eta_j-)-\zeta_2^n(\eta_j+)}{\zeta_2^n(\eta_j-)+\zeta_2^n(\eta_j+)}\qquad(1\leq j\leq q).
\]
(Note that $s_j$ is unchanged if $\zeta_2^n$ is replaced by $\zeta_2^n\widetilde{\zeta}$.) Then according to (\ref{g-step-composition}),
\begin{equation}\label{concatenation-approximation}
\begin{split}
g_{(\zeta_1^{n\;\smallfrown}\zeta_2^n)\widetilde{\zeta}}^\omega(\xi)=&\varphi_{\mu_1,r_1}\circ\cdots\circ\varphi_{\mu_p,r_p}\circ\varphi_{e^{2i(x_1-y_p)\omega},r}\circ\varphi_{\nu_1,s_1}\circ\cdots\circ\varphi_{\nu_q,s_q}(e^{2i(x_2-\eta_q)\omega}\xi)\\
=&g_{\zeta_1^n}^\omega\circ\varphi_{1,r}\circ g_{\zeta_2^n}^\omega(\xi),
\end{split}
\end{equation}
since
\[
\begin{split}
\varphi_{\mu_1,r_1}\circ\cdots\circ\varphi_{\mu_p,r_p}\circ\varphi_{e^{2i(x_1-y_p)\omega},r}(\xi)
=&\varphi_{\mu_1,r_1}\circ\cdots\circ\varphi_{\mu_p,r_p}\bigl(e^{2i(x_1-y_p)\omega}\varphi_{1,r}(\xi)\bigr)\\
=&g_{\zeta_1^n}^\omega\circ\varphi_{1,r}(\xi).
\end{split}
\]
By Theorem~\ref{thm-continuous}, $g_{\zeta_1^n}^\omega\rightarrow g_{\zeta_1}^\omega$, $g_{\zeta_2^n}^\omega\rightarrow g_{\zeta_2}^\omega$ and $g_{(\zeta_1^{n\;\smallfrown}\zeta_2^n)\widetilde{\zeta}}^\omega\rightarrow g_{\zeta\widetilde{\zeta}}^\omega$.  Taking limits in (\ref{concatenation-approximation}) and applying the foregoing values yields $g_{\zeta\widetilde{\zeta}}=g_{\zeta_1}\circ\varphi_{1,r}\circ g_{\zeta_2}$.  In the case $c=0=r$ this reduces to $g_{\zeta}=g_{\zeta_1}\circ g_{\zeta_2}$. 
\end{pf}

Proposition~\ref{prop-concatenation} extends in an obvious way to the case where $\widetilde{\zeta}\in\step_+(X)$ has an arbitrary number $n$ of jump points, involving a representation of $\zeta$ as a concatenation of $n+1$ terms, and the interposition of $n$ corresponding terms $\varphi_{1,r_j}$.

\begin{thm}\label{thm-g-composition}
Given $X=(x_0,x_1)$, let $\zeta\in\pwacplus(X)$. Suppose $\zeta\not\in C(X)$, and denote by $y_j$ $(1\leq j\leq n)$ the points of discontinuity of $\zeta$, indexed according to their natural order,
\[
x_0<y_1<\cdots<y_n<x_1,
\]
with $y_0=x_0$ and $y_{n+1}=x_1$.  For each $1\leq j\leq n+1$, set $\zeta_j=\zeta|_{(y_{j-1},y_j)}$, so that 
\[
\zeta=\zeta_1^{\;\frown}\cdots^{\;\frown}\zeta_{n+1}.
\]
Write 
\[
r_j=\frac{\zeta(y_j-)-\zeta(y_j+)}{\zeta(y_j-)+\zeta(y_j+)}\qquad(1\leq j\leq n). 
\]
Then 
\[
g_\zeta=g_{\zeta_1}\circ\varphi_{1,r_1}\circ\cdots\circ g_{\zeta_n}\circ\varphi_{1,r_n}\circ g_{\zeta_{n+1}}.
\]
\end{thm}

The next step is to study the behaviour of $g_\zeta^\omega(\xi)$ as a function of $\omega$.   
For this an approach involving explicit formulas, completely different from that in \S\ref{sec-technical}, turns out to be fruitful.

A key to explicit representation of $R(\omega)=g_\zeta^\omega(0)$, when $\zeta\in\pwacplus$, is purely singular approximation via a uniformly convergent sequence of step functions $\zeta_n\rightarrow\zeta$.  Initially the continuous case $\zeta\in C_+(X)$ will be analyzed, in two steps: more regular $\zeta\in C^1_+(X)$ will be uniformly approximated by evenly spaced step functions; then general non-smoothly differentiable $\zeta\in C_+(X)$ will be approximated uniformly by more regular $\zeta_n\in C_+^1(X)$.  The generalization from results for $\zeta\in C_+(X)$ to discontinuous $\zeta\in\pwacplus$ is provided by Theorem~\ref{thm-g-composition}.  Before introducing the two flavours of the harmonic exponential operator---the crucial intermediary linking $\zeta$ to $R$---it will be useful first to establish a technical bound applicable to all $\zeta\in\pwacplus$. 

Each $\zeta\in\pwacplus$ has bounded variation by virtue of absolute continuity.  Fix $X=(x_0,x_1)$ and let $\zeta\in\reg_+(X)$ be of bounded variation. Note this implies $\log\zeta$ is also of bounded variation, since $\zeta$ is by definition bounded and bounded away from 0.  A step function $\tilde\zeta\in\step_+(X)$ of the form 
\begin{equation}\label{step-form}
\tilde\zeta=\sum_{j=1}^{n+1}c_j\chi_{(y_{j-1},y_j)}\quad\mbox{ where }\quad x_0=y_0<y_1<\cdots<y_n<y_{n+1}=x_1
\end{equation}
is said to be an \emph{interpolating approximant} to $\zeta$ if $c_j\in\ran\zeta_{(y_{j-1},y_j)}$ for each $1\leq j\leq n+1$.  
\begin{prop}\label{prop-bounded-variation} 
Fix $X=(x_0,x_1)$.  Let $\zeta\in\reg_+(X)$ have bounded variation, set $R=g_\zeta(0)$, and let $\mathscr{V}$ denote the total variation of $\frac{1}{2}\log\zeta$.  Let $\tilde\zeta\in Step_+(X)$ be an interpolating approximant to $\zeta$, and write $\tilde R=g_{\tilde\zeta}(0)$.  Then for almost every $\omega\in\real$,
\[
\bigl|\tilde R(\omega)\bigr|\leq\tanh \mathscr{V}\quad\mbox{ and }\quad \bigl|R(\omega)\bigr|\leq\tanh \mathscr{V}.
\]
\end{prop}
\begin{pf}
With $\tilde\zeta$ as in (\ref{step-form}), Proposition~\ref{prop-step-composition} implies $g_{\tilde\zeta}^\omega(\xi)$ has the form
\begin{equation}\label{g-form}
g_{\tilde\zeta}^\omega(\xi)=\varphi_{\mu_1,r_1}\circ\cdots\circ\varphi_{\mu_n,r_n}\bigl(e^{2i(x_1-y_n)\omega}\xi\bigr)
\end{equation}
where $\mu_j=e^{2i(y_j-y_{j-1})}$ and
\begin{equation}\label{r-form}
r_j=\frac{c_j-c_{j+1}}{c_j+c_{j+1}}=\tanh\textstyle\left(\frac{1}{2}\log c_j-\frac{1}{2}\log c_{j+1}\right)\qquad(1\leq j\leq n).
\end{equation}
Repeated application of the inequality 
\[
 \left|\varphi_{\mu,\rho}(\xi)\right|
\leq \frac{|\rho|+|\xi|}{1+|\rho\xi|}=\tanh\left(\tanh^{-1}|\rho|+\tanh^{-1}|\xi|\right)\qquad\bigl(\mu\in S^1, \rho, \xi\in\ddd\bigr)
\]
to the composition of disk automorphisms (\ref{g-form}) yields
\begin{equation}\label{g-bound}
\bigl|g_{\tilde\zeta}^\omega(\xi)\bigr|\leq \tanh\left(\tanh^{-1}|\xi|+\textstyle\sum_{j=1}^n\tanh^{-1}|r_j|\right),
\end{equation}
which, by (\ref{r-form}), in turn yields
\[
\bigl|g_{\tilde\zeta}^\omega(\xi)\bigr|\leq \tanh\bigl(\tanh^{-1}|\xi|+\sum_{j=1}^n\left|\textstyle\frac{1}{2}\log c_j-\frac{1}{2}\log c_{j+1}\right|\bigr)\leq \tanh\bigl(\tanh^{-1}|\xi|+\mathscr{V}\bigr),
\]
since $\tilde\zeta$ is interpolating. 
The desired inequality $\bigl|\tilde R(\omega)\bigr|\leq\tanh \mathscr{V}$ then follows from the fact that $\tilde R(\omega)=g_{\zeta_n}^\omega(0)$. 

By density of $\step_+(X)$ in $\reg_+(X)$ there exists a sequence $\zeta_n\in\step_+(X)$ such that $\zeta_n\rightarrow\zeta$ uniformly. Moreover, without loss of generality, $\zeta_n$ may be chosen to be interpolating, so that 
\begin{equation}\label{Rn-variation-bound}
\bigl| R_n(\omega)\bigr|\leq\tanh \mathscr{V}\qquad(n\geq 1,\omega\in\real),
\end{equation}
where $R_n=g_{\zeta_n}(0)$.  Theorem~\ref{thm-continuous} implies $\lim_{n\rightarrow\infty}R_n(\omega)=R(\omega)$ for almost every $\omega\in\real$; therefore $\bigl|R(\omega)\bigr|\leq\tanh \mathscr{V}$ almost everywhere by (\ref{Rn-variation-bound}). 
\end{pf}

\section{The harmonic exponential operator\label{sec-harmonic-exponential}}

\subsection{The singular harmonic exponential operator\label{sec-singular-harmonic}}

Let $X=(x_0,x_1)$ and $\zeta\in\step_+(X)$. Denote by $r:X\rightarrow\real$ the reflectivity function associated to $\zeta$ given by
\begin{equation}\label{reflectivity-function}
r(x)=\frac{\zeta(x-)-\zeta(x+)}{\zeta(x-)+\zeta(x+)}=\tanh\log\sqrt{\zeta(x-)/\zeta(x+)}\qquad(x\in X).
\end{equation}
For $\omega\in\real$, define
\begin{equation}\label{Z}
\somop:\step(X)\rightarrow\step(X),\qquad\somop f(y)=\sum_{x_0<x<y}r(x)e^{2i(x-x_0)\omega}\overline{f(x)}. 
\end{equation}
Note that $r(x)=0$ unless $x$ is a jump point of $\zeta$; so $\somop=0$ if $\zeta$ is constant. 
For $y\in (x_0,x_1]$, define the singular harmonic exponential operator
\[
\calE^{(x_0,y)}:\step_+(X)\rightarrow C^{\infty}(\real),\qquad\zeta\mapsto\calE^{(x_0,y)}_\zeta
\]
by the formula
\begin{equation}\label{singular-harmonic-exponential}
\calE^{(x_0,y)}_\zeta(\omega)=(1-\somop)^{-1}\one(y).
\end{equation}
In particular, $\calE^{(x_0,y)}_\zeta(\omega)=1$ is constant if $\zeta$ is constant.  
\begin{prop}\label{prop-singular-harmonic}
Fix $X=(x_0,x_1)$.  Let $\zeta\in\step(X)$ have $n\geq 1$ jump points $y_j$, indexed according to their natural order 
\[
x_0<y_1<\cdots<y_n<x_1, 
\]
and let $\nu$ or $\nu_j$ denote elements of the index set $\{1,\ldots, n\}$. Define
\begin{equation}\label{kappa-simplex-definition-0}
\kappa(s_1,\ldots,s_j)=2\sum_{\nu=1}^j(-1)^{j-\nu}(s_\nu-x_0)\qquad(1\leq j\leq n),
\end{equation}
and let $y\in(x_0,x_1]$. Then:
\begin{enumerate}[label={(\roman*)},itemindent=0em]
\item \ $\somop^k\one(y)=\displaystyle\sum\limits_{x_0<y_{\nu_1}<\cdots<y_{\nu_k}<y}\exp\bigl(i\omega\kappa(y_{\nu_1},\ldots,y_{\nu_k})\bigr)\prod\limits_{j=1}^kr(y_{\nu_j})\qquad(1\leq k\leq n)$;\label{somop-power}\\[5pt]
\item \ $\somop^k=0$ if $k>n$;\label{somop-zero}\\[5pt]
\item \ $\displaystyle
\calE^{(x_0,y)}_\zeta(\omega)=1+\sum_{j=1}^n\somop^j\one(y);$\label{singular-exponential-finite}\\[5pt]
\item \ $\calE^{(x_0,y)}_\zeta$ is almost periodic; \label{E-almost-periodic}\\[5pt]
\item \ $\calE^{(x_0,y)}_\zeta$ is periodic with period $p>0$ if and only if the quantities $y_j-x_0$ are all integer multiples of $\pi/p$ $(1\leq j\leq n)$.\label{periodicity-criterion}
\end{enumerate}
\end{prop}
\begin{pf}
According to (\ref{Z}), 
\[
\begin{split}
\somop f(y)&=\sum_{x_0<{\eta_1}<y}r({\eta_1})e^{2i({\eta_1-x_0})\omega}\overline{f({\eta_1})},\\
\somop^2 f(y)&=\sum_{x_0<{\eta_2}<y}r({\eta_2})e^{2i({\eta_2-x_0})\omega}\left(\sum_{x_0<{\eta_1}<{\eta_2}}r({\eta_1})e^{-2i({\eta_1-x_0})\omega}f({\eta_1})\right)\\
&=\sum_{x_0<{\eta_1}<{\eta_2}<y}f({\eta_1})r({\eta_1})r({\eta_2})e^{2i({\eta_2}-{\eta_1})\omega}.
\end{split}
\]
Repeated iteration yields the general form
\begin{equation}\label{Z-k}
\somop^k f(y)=\sum_{x_0<{\eta_1}<\cdots<{\eta_k}<y}f_k({\eta_1})e^{i\omega\kappa(\eta_1,\ldots,\eta_k)}\prod_{j=1}^kr(\eta_j),
\end{equation}
where $f_k=f$ if $k$ is even, and $f_k=\overline{f}$ if $k$ is odd.  
Now, the reflectivity function $r(x)$ is non-zero only if $x=y_j$ for some $1\leq j\leq n$.  Thus, if $k>n$ and 
\[
x_0<\eta_1<\cdots<\eta_k<y,
\]
then at least one $\eta_j$ is not a jump point of $\zeta$, and so $\prod_{j=1}^kr(\eta_j)=0$, which implies part~\ref{somop-zero} of the proposition.  If $1\leq k\leq n$, replacing $f$ with $\one$ in (\ref{Z-k}) yields part~\ref{somop-power}. Part~\ref{singular-exponential-finite} follows from \ref{somop-power} and \ref{somop-zero}.  

In light of \ref{somop-power}, the formulation \ref{singular-exponential-finite} expresses $\calE^{(x_0,y)}_\zeta(\omega)$ as a finite combination of exponentials $\exp\bigl(i\omega\kappa(y_{\nu_1},\ldots,y_{\nu_k})\bigr)$.  Since $\kappa$ is real-valued it follows that $\calE^{(x_0,y)}_\zeta$ is almost periodic.  Using (\ref{kappa-simplex-definition-0}),
periodicity of $\calE^{(x_0,y)}_\zeta$ is easily seen to be equivalent to co-rationality of the numbers $y_{j}-x_0$ $(1\leq j\leq n)$, the precise period (and its multiples) corresponding to part~\ref{periodicity-criterion}.
\end{pf}

\subsection{Representation of $g_\zeta$ for $\zeta\in\step_+(X)$\label{sec-singular-representation}}

$\aut\mathbb{P}$ consists of holomorphic bijections of the unit disk of the form
\begin{equation}\label{disk-automorphisms-2}
\varphi_{\mu,\rho}:\ddd\rightarrow\ddd,\quad \varphi_{\mu,\rho}(\xi)=\mu\frac{\xi+\rho}{1+\bar\rho\xi}\quad\mbox{ where }\quad\mu\in S^1\mbox{ and }\rho\in\ddd.
\end{equation}
Define 
\begin{equation}\label{pro}
\pro:\complex^2\setminus \{0\}\rightarrow\complex\cup\{\infty\}\qquad\pro\,\binom{w}{z}=w/z. 
\end{equation}
Two-by-two matrices of the form 
\begin{equation}\label{M-mu-rho}
M_{\mu,\rho}=\begin{pmatrix}\mu&\mu\rho\\ \bar\rho&1\end{pmatrix}\qquad\left(\mu\in S^1, \rho\in\ddd\right)
\end{equation}
serve as homogeneous coordinates representing $\aut\mathbb{P}$ in the sense that 
\begin{equation}\label{homogeneous-representation}
\varphi_{\mu,\rho}(\xi)=\pro\left(\alpha M_{\mu,\rho}\binom{\xi}{1}\right)\quad\mbox{ for any }0\neq\alpha\in\complex, 
\end{equation}
and, for any $(\mu_1,\rho_1), (\mu_2,\rho_2)\in S^1\times\ddd$,
\begin{equation}\label{multiplicative}
\pro\left(M_{\mu_1,\rho_1}M_{\mu_2,\rho_2}\binom{\xi}{1}\right)=\varphi_{\mu_1,\rho_1}\circ\varphi_{\mu_2,\rho_2}(\xi).
\end{equation}
(For present purposes the representation (\ref{M-mu-rho}) is preferable to the standard double cover of $\aut\mathbb{P}$ by $SU(1,1)$, 
\begin{equation}\label{su-one-one}
\begin{pmatrix}w&z\\ \bar{z}&\bar{w}\end{pmatrix}\mapsto\varphi_{\frac{w}{\bar{w}},\frac{z}{w}}\quad\mbox{ where }\quad (w,z)\in\complex^2\mbox{ such that }|w|^2-|z|^2=1, 
\end{equation}
because, unlike (\ref{su-one-one}), the mapping $M_{\mu,\rho}\mapsto\varphi_{\mu,\rho}$ is injective.)

\begin{prop}\label{prop-matrix-product}
Fix $X=(x_0,x_1)$.  Let $\zeta\in\step(X)$ have $n\geq 1$ jump points $y_j$, indexed according to their natural order 
\[
x_0<y_1<\cdots<y_n<x_1, 
\]
and write $y_0=x_0$.  Set 
\[
\mu_{n+1}=e^{2i(x_1-y_n)\omega},\quad r_{n+1}=0,\quad \mu_j=e^{2i(y_j-y_{j-1})\omega},\quad
r_j=\frac{\zeta(y_j-)-\zeta(y_j+)}{\zeta(y_j-)+\zeta(y_j+)}\qquad(1\leq j\leq n)
\]
and
\[
M_0=\begin{pmatrix}1&1\\ -1&1\end{pmatrix},\quad M_j=\begin{pmatrix}\mu_j&\mu_jr_j\\ r_j&1\end{pmatrix}\qquad(1\leq j\leq n+1).
\]
Then 
\begin{equation}\label{matrix-product-E}
M_0M_1\cdots M_{n+1}=
\begin{pmatrix}
e^{2i(x_1-x_0)\omega}\overline{\calE_\zeta^X(\omega)}& \calE_\zeta^X(\omega)\\[10pt]
-e^{2i(x_1-x_0)\omega}\overline{\calE_{1/\zeta}^X(\omega)}& \calE_{1/\zeta}^X(\omega)
\end{pmatrix}.
\end{equation}
\end{prop} 
\begin{pf} This is a matter of careful bookkeeping. In detail, consider first products of $n$ matrices of the form 
\begin{equation}\label{matrix-form-1}
\mathcal{M}_j=\begin{pmatrix}\mu_j &\mu_j\rho_j\\ \bar\rho_j&1\end{pmatrix}=
\begin{pmatrix}\mu_j&0\\ 0&1\end{pmatrix}\left(I+\begin{pmatrix}0&\rho_j\\ \bar\rho_j&0\end{pmatrix}\right)
\qquad\mu_j\in S^1, \rho_j\in\ddd
\end{equation}
corresponding to disk automorphisms.  The identity
\begin{equation}\label{basic}
\left(I+\begin{pmatrix}0&p\\ q&0\end{pmatrix}\right)\begin{pmatrix}s&0\\ 0&1\end{pmatrix}=
\begin{pmatrix}s&0\\ 0&1\end{pmatrix}\left(I+\begin{pmatrix}0&s^{-1}p\\ sq&0\end{pmatrix}\right)
\end{equation}
applied to the product $\mathcal{M}_1\cdots \mathcal{M}_n$ yields 
\begin{equation}\label{product-1}
\mathcal{M}_1\cdots \mathcal{M}_n=\begin{pmatrix}\mu_1\cdots\mu_n&0\\ 0&1\end{pmatrix}\left(I+\Lambda_1\right)\cdots\left(I+\Lambda_n\right),
\end{equation}
where 
\begin{equation}\label{lambda-formulas}
\lambda_j=\rho_j\overline{\mu_{j+1}\cdots\mu_n},\quad \Lambda_j=\begin{pmatrix}0&\lambda_j\\ \bar{\lambda}_j&0\end{pmatrix}\qquad(1\leq j\leq n).
\end{equation}
Expanding the right-hand side of (\ref{product-1}),
\begin{equation}\label{product-2}
\mathcal{M}_1\cdots \mathcal{M}_n=\begin{pmatrix}\mu_1\cdots\mu_n&0\\ 0&1\end{pmatrix}
\left(I+\sum_{k=1}^n\sum_{1\leq j_1<\cdots<j_k\leq n}\Lambda_{j_1}\cdots\Lambda_{j_k}\right).
\end{equation}
Noting that $\Lambda_{j_1}\Lambda_{j_2}=\begin{pmatrix}\lambda_{j_1}\bar{\lambda}_{j_2}&0\\ 0&\bar{\lambda}_{j_1}\lambda_{j_2}\end{pmatrix}$, one computes
\begin{equation}\label{A-B-matrix}
I+\sum_{k=1}^n\sum_{1\leq j_1<\cdots<j_k\leq n}\Lambda_{j_1}\cdots\Lambda_{j_k}=\begin{pmatrix}A&B\\ \bar{B}&\bar{A}\end{pmatrix},
\end{equation}
where 
\begin{equation}\label{A-B-definition}
\begin{split}
A&=1+\sum_{k=1}^{\lfloor n/2\rfloor}\sum_{1\leq j_1<\cdots<j_{2k}\leq n}\lambda_{j_1}\bar{\lambda}_{j_2}\cdots\lambda_{2k-1}\bar{\lambda}_{2k},\\
B&=\sum_{k=0}^{\lfloor(n-1)/2\rfloor}\sum_{1\leq j_1<\cdots<j_{2k+1}\leq n}\lambda_{j_1}\bar{\lambda}_{j_2}\cdots\lambda_{2k-1}\bar{\lambda}_{2k}\lambda_{2k+1}.
\end{split}
\end{equation}
In terms of the $\mu_j,\rho_j$, 
\begin{equation}\label{lambda-product-rewrite}
\begin{split}
\lambda_{j_1}\bar{\lambda}_{j_2}\cdots\lambda_{2k-1}\bar{\lambda}_{2k}=&\rho_{j_1}\bar{\rho}_{j_2}\cdots\rho_{j_{2k-1}}\bar{\rho}_{j_{2k}}\cdot\overline{\mu_{j_1+1}\cdots\mu_{j_2}}\cdots\overline{\mu_{j_{2k-1}+1}\cdots\mu_{j_{2k}}},\\
\lambda_{j_1}\bar{\lambda}_{j_2}\cdots\lambda_{2k-1}\bar{\lambda}_{2k}\lambda_{2k+1}=&
\rho_{j_1}\bar{\rho}_{j_2}\cdots\rho_{j_{2k-1}}\bar{\rho}_{j_{2k}}\rho_{j_{2k+1}}\cdot\\ 
&\overline{\mu_{j_1+1}\cdots\mu_{j_2}}\cdots\overline{\mu_{j_{2k-1}+1}\cdots\mu_{j_{2k}}}\cdot\overline{\mu_{j_{2k+1}+1}\cdots\mu_{n}}.
\end{split}
\end{equation}
Setting
\[
M_0=\begin{pmatrix}1&1\\ -1&1\end{pmatrix}\quad\mbox{ and }\quad \mathcal{M}_{n+1}=\begin{pmatrix}\mu_{n+1}&0\\ 0&1\end{pmatrix},
\]
this yields
\begin{equation}\label{expanded-AB-product}
M_0\mathcal{M}_1\cdots\mathcal{M}_{n+1}=\begin{pmatrix}\mu_{n+1}\bigl(\overline{B}+A\bigr)&\overline{A}+B\\ \mu_{n+1}\bigl(\overline{B}-A\bigr)&\overline{A}-B\end{pmatrix}.
\end{equation}

Define $\mu_j$ as in the hypothesis of Proposition~\ref{prop-matrix-product} for $1\leq j\leq n+1$, and set $\rho_j=r_j$ $(1\leq j\leq n)$, with the $r_j$ defined as in Proposition~\ref{prop-matrix-product}, so that $\mathcal{M}_j=M_j$ $(1\leq j\leq n+1)$.  Then formulas (\ref{A-B-definition}) and (\ref{lambda-product-rewrite}) in combination with Proposition~\ref{prop-singular-harmonic} yield
\[
M_0M_1\cdots M_{n+1}=\begin{pmatrix}\mu_{n+1}\bigl(\overline{B}+A\bigr)&\overline{A}+B\\ \mu_{n+1}\bigl(\overline{B}-A\bigr)&\overline{A}-B\end{pmatrix}=
\begin{pmatrix}
e^{2i(x_1-x_0)\omega}\overline{\calE_\zeta^X(\omega)}& \calE_\zeta^X(\omega)\\[10pt]
-e^{2i(x_1-x_0)\omega}\overline{\calE_{1/\zeta}^X(\omega)}& \calE_{1/\zeta}^X(\omega)
\end{pmatrix}.
\]
This completes the proof. \end{pf}

\begin{cor}\label{cor-determinantal-formula} For $\zeta$ as in Proposition~\ref{prop-matrix-product},
\[
\Re\left( \calE_{\zeta}^X(\omega)\,\overline{\calE_{1/\zeta}^X(\omega)}\right)=\prod_{j=1}^n(1-r_j^2)\qquad(\omega\in\real).
\]
\end{cor}
\begin{pf}Evaluate determinants on either side of (\ref{matrix-product-E}).\end{pf}

\begin{thm}\label{thm-singular-representation}
Let $\zeta\in\step_+(X)$ for some $X=(x_0,x_1)$.  Then: 
\begin{enumerate}[label={(\roman*)},itemindent=0em]
\item
\[
g_\zeta^\omega(\xi)=\mu\frac{\xi+\rho}{1+\overline{\rho}\xi}\qquad(\xi\in\ddd),
\]
where 
\[
\mu=e^{2i(x_1-x_0)\omega}\frac{\overline{\calE_{\zeta}^X(\omega)+\calE_{1/\zeta}^X(\omega)}}{\calE_{\zeta}^X(\omega)+\calE_{1/\zeta}^X(\omega)}\quad\mbox{ and }\quad 
\rho=e^{-2i(x_1-x_0)\omega}\frac{\calE_{\zeta}^X(\omega)-\calE_{1/\zeta}^X(\omega)}{\overline{\calE_{\zeta}^X(\omega)+\calE_{1/\zeta}^X(\omega)}};
\]\label{g-singular-representation}
\item
letting $R$ denote the reflection coefficient determined by $\zeta$,
\[
R=\frac{\calE_{\zeta}^X-\calE_{1/\zeta}^X}{\calE_{\zeta}^X+\calE_{1/\zeta}^X}.
\]
\label{singular-R-representation}
\end{enumerate}
\end{thm}
\begin{pf}
Using (\ref{multiplicative}), the representation 
\[
g_\zeta^\omega(\xi)=\varphi_{\mu_1,r_1}\circ\cdots\circ\varphi_{\mu_n,r_n}\bigl(e^{2i(x_1-y_n)\omega}\xi\bigr).
\]
established in Proposition~\ref{prop-step-composition}
may be rewritten in matrix form as 
\[
\begin{split}
g_\zeta^\omega(\xi)&=\pro\left(M_1\cdots M_{n+1}\binom{\xi}{1}\right)\\
&=\pro\left(M_0^{-1}\begin{pmatrix}
e^{2i(x_1-x_0)\omega}\overline{\calE_\zeta^X(\omega)}& \calE_\zeta^X(\omega)\\[10pt]
-e^{2i(x_1-x_0)\omega}\overline{\calE_{1/\zeta}^X(\omega)}& \calE_{1/\zeta}^X(\omega)
\end{pmatrix}\binom{\xi}{1}\right)\quad\mbox{ by Prop.~\ref{prop-matrix-product} }\\
&=\mu\frac{\xi+\rho}{1+\overline{\rho}\xi}.
\end{split}
\]
Setting $\xi=0$ then yields part~\ref{singular-R-representation} by the identity $R(\omega)=g^\omega_\zeta(0)$. 
\end{pf}

\begin{cor}\label{cor-singular-real-part}
Given $X=(x_0,x_1)$, let $\zeta\in\step_+(X)$ have $n\geq 1$ jump points $y_j$, with associated reflectivities $r_j=(\zeta(y_j-)-\zeta(y_j+))/(\zeta(y_j-)+\zeta(y_j+))$ $(1\leq j\leq n)$, and set $R(\omega)=g_\zeta^\omega(0)$. Then 
\[
2\Re\frac{R}{1-R}=\frac{\prod_{j=1}^n(1-r_j^2)}{\bigl|\calE_{1/\zeta}^X\bigr|^2}-1. 
\] 
\end{cor}
\begin{pf}
Note first that 
\begin{equation}\label{singular-fraction-real-part}
\begin{split}
\Re\frac{\calE_{\zeta}^X}{\calE_{1/\zeta}^X}&=\frac{1}{2}\frac{\overline{\calE_{1/\zeta}^X}\,\calE_{\zeta}^X+\calE_{1/\zeta}^X\,\overline{\calE_{\zeta}^X}}{\bigl|\calE_{1/\zeta}^X\bigr|^2}\\
&=\frac{\prod_{j=1}^n(1-r_j^2)}{\bigl|\calE_{1/\zeta}^X\bigr|^2}\quad\mbox{ by Cor.~\ref{cor-determinantal-formula}}.
\end{split}
\end{equation}
Part~\ref{singular-R-representation} of Theorem~\ref{thm-singular-representation} implies 
\[
\frac{2R}{1-R}=\frac{\calE_{\zeta}^X}{\calE_{1/\zeta}^X}-1,
\]
so that 
\[
2\Re\frac{R}{1-R}=\frac{\prod_{j=1}^n(1-r_j^2)}{\bigl|\calE_{1/\zeta}^X\bigr|^2}-1\]
by (\ref{singular-fraction-real-part}). 
\end{pf}

\subsection{Evenly-spaced jump points and OPUC\label{sec-OPUC}}

The following standard facts concerning finite sequences of orthogonal polynomials on the unit circle are needed in the sequel. Terminology conforms to that of Simon \cite{SiOPUC1:2005}. 
Given $n\geq 1$ points  $\rho_j\in\ddd$, let $z$ denote a complex variable, and set 
\begin{equation}\label{OPUC-M}
\mathcal{M}_0=\begin{pmatrix}1&1\\ -1&1\end{pmatrix}\quad\mbox{ and }\quad
\mathcal{M}_j=\begin{pmatrix}
z&z\rho_j\\
\overline{\rho_j}&1
\end{pmatrix}\qquad(1\leq j\leq n).
\end{equation}
Write 
\begin{equation}\label{OPUC-notation}
\begin{pmatrix}
\Psi_j(z)&\Psi_j^\ast(z)\\
-\Phi_j(z)&\Phi_j^\ast(z)
\end{pmatrix}
=\mathcal{M}_0\mathcal{M}_1\cdots\mathcal{M}_j\qquad(0\leq j\leq n)
\end{equation}
so that for $0\leq j\leq n-1$,
\begin{equation}\label{OPUC-recursion}
\begin{pmatrix}
\Psi_{j+1}(z)&\Psi_{j+1}^\ast(z)\\
-\Phi_{j+1}(z)&\Phi_{j+1}^\ast(z)
\end{pmatrix}
=
\begin{pmatrix}
z\Psi_j(z)+\overline{\rho_{j+1}}\Psi_j^\ast(z)&\Psi_j^\ast(z)+\rho_{j+1}z\Psi_j(z)\\
-\left(z\Phi_j(z)-\overline{\rho_{j+1}}\Phi_j^\ast(z)\right)&\Phi_j^\ast(z)-\rho_{j+1}z\Phi_j(z)
\end{pmatrix}. 
\end{equation}
The recurrence relations (\ref{OPUC-recursion}) characterize the sequences of monic orthogonal polynomials $\Phi_j$ and $\Psi_j$, and their duals $\Phi_j^\ast$ and $\Psi_j^\ast$.  More precisely, $\{\Phi_j\}_{j=0}^n$ is said to be the sequence of monic orthogonal polynomials on the unit circle determined by Verblunsky coefficients $\rho_1,\ldots,\rho_n$. And $\{\Psi_j\}_{j=0}^n$ is the associated sequence of monic orthogonal polynomials on the unit circle determined by $-\rho_1,\ldots,-\rho_n$.  The dual polynomials $\Phi_j^\ast$ and $\Psi_j^\ast$ are related to $\Phi_j$ and $\Psi_j$, respectively, by the formulas
\begin{equation}\label{dual-polynomial-formulas}
\Phi_j^\ast(z)=z^j\overline{\Phi_j(1/\bar{z})}\quad\mbox{ and }\quad\Psi_j^\ast(z)=z^j\overline{\Psi_j(1/\bar{z})}\qquad(0\leq j\leq n). 
\end{equation}
The sequences $\{\Phi_j\}_{j=0}^n$ and $\{\Psi_j\}_{j=0}^n$ are orthogonal in $L^2(S^1, d\mu)$ and $L^2(S^1, d\nu)$, respectively, for various probability measures $d\mu$ and $d\nu$ on the unit circle, including
\begin{equation}\label{unit-circle-measures}
d\mu(e^{i\theta})=\frac{\prod_{j=1}^n(1-|\rho_j|^2)}{\bigl|\Phi_n^\ast(e^{i\theta})\bigr|^2}\,\frac{d\theta}{2\pi}\qquad\mbox{ and }\qquad
d\nu(e^{i\theta})=\frac{\prod_{j=1}^n(1-|\rho_j|^2)}{\bigl|\Psi_n^\ast(e^{i\theta})\bigr|^2}\,\frac{d\theta}{2\pi}.
\end{equation}
(The theory of OPUC is generally concerned with the correspondence among non-degenerate probability measures on the circle,
infinite sequences of Verblunksy coefficients, and infinite sequences of polynomials. The above measures $d\mu$ and $d\nu$ correspond to the infinite sequence of coefficients obtained by extending $\rho_1,\ldots,\rho_n$ by $\rho_{n+j}=0$, with corresponding polynomials $\Phi_{n+j}(z)=z^j\Phi_n(z)$ and $\Psi_{n+j}(z)=z^j\Psi_{n+j}(z)$ $(j\geq 1)$. However, in the present paper only finite sequences come into play.)
A crucial fact for present purposes is that the dual polynomials $\Phi_j^\ast$ and $\Psi_j^\ast$ are zero free on the closed unit disk $\overline{\ddd}$, which is easily proved using the recurrence relations (\ref{OPUC-recursion}).  

In detail, let $0
\leq j\leq n-1$ and note by (\ref{dual-polynomial-formulas}) that $\bigl|\Phi_j(z)/\Phi_j^\ast(z)\bigr|=1$ for all $z\in S^1$.  If $\Phi_j^\ast$ is zero free on $\overline{\ddd}$, then the maximum principle implies $\bigl|\Phi_j(z)/\Phi_j^\ast(z)\bigr|\leq 1$ for all $z\in\overline{\ddd}$.  Now, if $\Phi_{j+1}^\ast(z_0)=0$ for some $z_0\in\overline{\ddd}$, the recurrence (\ref{OPUC-recursion}) shows that $z_0r_{j+1}\Phi_j(z_0)/\Phi_j^\ast(z_0)=1$.  But this is impossible, since $|r_{j+1}|<1$, while $|z_0|\leq 1$ and $\bigl|\Phi_j(z_0)/\Phi_j^\ast(z_0)\bigr|\leq 1$.  So $\Phi_{j+1}^\ast$ is also zero free on $\overline{\ddd}$.  Given that $\Phi_0^\ast=1$ is zero free, the claimed result follows by induction. Exactly the same argument applies to $\Psi_j^\ast$ as well. 

This concludes the summary of needed properties of orthogonal polynomials on the unit circle.
The formulation of Proposition~\ref{prop-matrix-product} makes obvious a connection to OPUC that turns out to be very useful from a technical standpoint, as follows.
\begin{prop}\label{prop-opuc}
Fix $X=(x_0,x_1)$. Suppose $\zeta\in\step_+(X)$ has $n\geq 1$ equally-spaced jump points 
\[
y_j=x_0+j\frac{x_1-x_0}{n+1}\qquad(1\leq j\leq n).
\]
Write $\Delta=(x_1-x_0)/(n+1)$, and set 
\[
r_j=\frac{\zeta(y_j-)-\zeta(y_j+)}{\zeta(y_j-)+\zeta(y_j+)}\qquad(1\leq j\leq n).
\]
Let $\Psi_n(z)$ and $\Phi_n(z)$ denote the degree $n$ monic, orthogonal polynomials on the unit circle determined by Verblunsky coefficients $-r_1,\ldots,-r_n$  and $r_1,\ldots,r_n$, respectively, with respective dual polynomials $\Psi^\ast_n(z)$ and $\Phi_n^\ast(z)$. 
Then:\\[0pt]
\begin{enumerate}[label={(\roman*)},itemindent=0pt]
\item
$
\calE_\zeta^X(\omega)=\Psi_n^\ast(e^{2i\Delta\omega})\quad\mbox{ and }\quad\calE_{1/\zeta}^X(\omega)=\Phi_n^\ast(e^{2i\Delta\omega})\qquad (\omega\in\real)
$;\label{singular-harmonic-OPUC}\\[5pt]
\item
the tempered distribution $\displaystyle\left({1}/{\overline{\calE_{1/\zeta}^X}}\right)^{\rule{-4pt}{0pt}\rwc{\rule{12pt}{0pt}}}$ 
is supported on a subset of $(-\infty, 0]$;  \label{reciprocal-transform-support}
\item
\[
g_\zeta^\omega(\xi)=\mu\frac{\xi+\rho}{1+\overline{\rho}\xi}\qquad(\xi\in\ddd),
\]
where 
\[
\mu
=e^{2i\Delta\omega}\,\frac{\Psi_n(e^{2i\Delta\omega})+\Phi_n(e^{2i\Delta\omega})}{\Psi_n^\ast(e^{2i\Delta\omega})+\Phi_n^\ast(e^{2i\Delta\omega})}\quad
\mbox{ and }\quad
\rho
=e^{-2i\Delta\omega}\,\frac{\Psi_n^\ast(e^{2i\Delta\omega})-\Phi_n^\ast(e^{2i\Delta\omega})}{\Psi_n(e^{2i\Delta\omega})+\Phi_n(e^{2i\Delta\omega})};
\]\label{g-equally-spaced-representation}
\item
letting $R$ denote the reflection coefficient determined by $\zeta$,
\[
R(\omega)=\frac{\Psi_n^\ast(e^{2i\Delta\omega})-\Phi_n^\ast(e^{2i\Delta\omega})}{\Psi_n^\ast(e^{2i\Delta\omega})+\Phi_n^\ast(e^{2i\Delta\omega})}\qquad(\omega\in\real);
\]
\label{equally-spaced-R-representation}
\item \label{Szego-theorem}
\[
\frac{\Delta}{\pi}\int\limits_{-\frac{\pi}{2\Delta}}^{\frac{\pi}{2\Delta}}-\log\left(1-\bigl|R(\omega)\bigr|^2\right)\,d\omega
=\sum_{j=1}^n-\log(1-|r_j|^2).
\]
\end{enumerate}
\end{prop}

\begin{pf}
Denote matrices $M_j$ as in Proposition~\ref{prop-matrix-product} for $1\leq j\leq n+1$. According to (\ref{OPUC-M}), if $z=e^{2i\Delta\omega}$ and $\rho_j=r_j$, then $\mathcal{M}_j=M_j$ $(0\leq j\leq n)$, and so
\[
\begin{pmatrix}\Psi_n(z)&\Psi^\ast_n(z)\\
-\Phi_n(z)& \Phi_n^\ast(z)\end{pmatrix}
= M_0\cdots M_n.
\]
It follows by Proposition~\ref{prop-matrix-product} that
\[
\begin{pmatrix}
e^{2i(x_1-x_0)\omega}\overline{\calE_\zeta^X(\omega)}& \calE_\zeta^X(\omega)\\[10pt]
-e^{2i(x_1-x_0)\omega}\overline{\calE_{1/\zeta}^X(\omega)}& \calE_{1/\zeta}^X(\omega)
\end{pmatrix}
=
\begin{pmatrix}\Psi_n(z)&\Psi^\ast_n(z)\\
-\Phi_n(z)& \Phi_n^\ast(z)\end{pmatrix}M_{n+1}
=
\begin{pmatrix}z\Psi_n(z)&\Psi^\ast_n(z)\\
-z\Phi_n(z)& \Phi_n^\ast(z)\end{pmatrix},
\]
giving the equations in part~\ref{singular-harmonic-OPUC}.  

The function $1/\Phi_n^\ast$ is holomorphic on the closed unit disk $\overline{\ddd}$ because $\Phi_n^\ast$ is holomorphic and zero-free on $\overline{\ddd}$. Its complex conjugate $1/\overline{\Phi_n^\ast}$ is therefore anti-holomorphic, having a Taylor series expansion in $\overline{z}$ of the form
\[
1/\overline{\Phi_n^\ast(z)}=\sum_{j=0}^\infty c_j\overline{z}^j
\]
that converges absolutely and uniformly on $\overline{\ddd}$.  Restricting the above expansion to the unit circle, and applying part~\ref{singular-harmonic-OPUC}, yields
\[
\displaystyle{1}/{\overline{\calE_{1/\zeta}^X(\omega)}}=\sum_{j=0}^\infty c_je^{-2ij\Delta\omega}\qquad(\omega\in\real).
\]
The latter is bounded (and $C^\infty$) and hence a tempered distribution. Its inverse Fourier transform
\[
\displaystyle\left({1}/{\overline{\calE_{1/\zeta}^X}}\right)^{\rule{-4pt}{0pt}\rwc{\rule{12pt}{0pt}}}\rule{-7pt}{0pt}(t)=
\sum_{j=0}^\infty c_j\delta(t+2j\Delta)
\]
is supported only at points of the form $t=-2j\Delta$ where $j\geq 0$, proving part~\ref{reciprocal-transform-support}.

Given~\ref{singular-harmonic-OPUC}, Theorem~\ref{thm-singular-representation} and formulas (\ref{dual-polynomial-formulas}) yield
\[
g_\zeta^\omega(\xi)=\mu\frac{\xi+\rho}{1+\overline{\rho}\xi}\qquad(\xi\in\ddd),
\]
where 
\[
\begin{split}
\mu&=e^{2i(x_1-x_0)\omega}\frac{\overline{\Psi_n^\ast(e^{2i\Delta\omega})+\Phi_n^\ast(e^{2i\Delta\omega})}}{\Psi_n^\ast(e^{2i\Delta\omega})+\Phi_n^\ast(e^{2i\Delta\omega})}
=e^{2i\Delta\omega}\,\frac{\Psi_n(e^{2i\Delta\omega})+\Phi_n(e^{2i\Delta\omega})}{\Psi_n^\ast(e^{2i\Delta\omega})+\Phi_n^\ast(e^{2i\Delta\omega})}\\[10pt] 
\mbox{ and }\quad\rho&=e^{-2i(x_1-x_0)\omega}\frac{\Psi_n^\ast(e^{2i\Delta\omega})-\Phi_n^\ast(e^{2i\Delta\omega})}{\overline{\Psi_n^\ast(e^{2i\Delta\omega})+\Phi_n^\ast(e^{2i\Delta\omega})}}
=e^{-2i\Delta\omega}\,\frac{\Psi_n^\ast(e^{2i\Delta\omega})-\Phi_n^\ast(e^{2i\Delta\omega})}{\Psi_n(e^{2i\Delta\omega})+\Phi_n(e^{2i\Delta\omega})},
\end{split}
\]
proving part~\ref{g-equally-spaced-representation}. Part \ref{equally-spaced-R-representation} is obtained by setting $\xi=0$ in \ref{g-equally-spaced-representation}.

A subsidiary claim to proving \ref{Szego-theorem} is that for every $0\leq j\leq n$, 
\begin{equation}\label{transmission-outer}
\bigl| \Psi_j(z)-\Phi_j(z)\bigr|<\bigl| \Psi_j^\ast(z)+\Phi_j^\ast(z)\bigr|\quad(z\in\overline{\ddd}).
\end{equation}
To prove (\ref{transmission-outer}), set 
\[
f_j(z)=\frac{ \Psi_j(z)-\Phi_j(z)}{\Psi_j^\ast(z)+\Phi_j^\ast(z)}\qquad(0\leq j\leq n),
\]
and note by (\ref{dual-polynomial-formulas}) that 
\begin{equation}\label{circle-dual-equality}
\left| \Psi_j(z)-\Phi_j(z) \right|=\left| \Psi^\ast_j(z)-\Phi^\ast_j(z) \right|\quad\mbox{ if }\quad |z|=1.
\end{equation}
Let $0\leq j\leq n-1$ and suppose $\Psi^\ast_j+\Phi^\ast_j$ is zero free on $\overline{\ddd}$, whence $f_j$ is holomorphic $\overline{\ddd}$.  The aim is to show $|f_j|<1$ and $\Psi^\ast_{j+1}+\Phi^\ast_{j+1}$ is zero free on $\overline{\ddd}$.  Observe that if $|z|=1$ then
\[
\left|f_j(z)\right|=\left|\frac{\Psi_j^\ast(z)-\Phi_j^\ast(z)}{\Psi_j^\ast(z)+\Phi_j^\ast(z)}\right|<1,
\]
by (\ref{circle-dual-equality}) and the fact that if $|z|=1$, then 
$
\bigl(\Psi_j^\ast(z)-\Phi_j^\ast(z)\bigr)/\bigl(\Psi_j^\ast(z)+\Phi_j^\ast(z)\bigr)
$
is the value at 0 of a disk automorphism. The latter assertion follows from setting $\rho_\nu=r_\nu$ in (\ref{OPUC-notation}) and writing
\[
\bigl(\Psi_j^\ast(z)-\Phi_j^\ast(z)\bigr)/\bigl(\Psi_j^\ast(z)+\Phi_j^\ast(z)\bigr)=\pro\left(\mathcal{M}_1\cdots\mathcal{M}_j\binom{0}{1}\right).
\]
The maximum principle therefore implies $\bigl|f_j(z)\bigr|<1$ for $z\in\overline{\ddd}$. To see that $\Psi^\ast_{j+1}+\Phi^\ast_{j+1}$ is zero free on $\overline{\ddd}$, note that (\ref{OPUC-recursion}) implies
\[
\Psi_{j+1}^\ast(z)+\Phi_{j+1}^\ast(z)=
\Psi_j^\ast(z)+\Phi_j^\ast(z)+r_{j+1}z\bigl(\Psi_j(z)-\Phi_j(z)\bigr).
\]
If $\Psi_{j+1}^\ast(z_0)+\Phi_{j+1}^\ast(z_0)=0$ and $|z_0|\leq 1$ then $\bigl|r_{j+1}z_0f_j(z_0)\bigr|=1$.
But $|r_{j+1}|<1$, so this contradicts the fact that $|f_j(z_0)|<1$.  Thus $\Psi^\ast_{j+1}+\Phi^\ast_{j+1}$ is zero free on $\overline{\ddd}$.  Now, $\Psi_0^\ast=\Phi_0^\ast=1$ and so $\Psi_0^\ast+\Phi_0^\ast =2$ is zero free.  Therefore $|f_j(z)|<1$ on $\overline{\ddd}$ for every $0\leq j\leq n$ by induction, proving (\ref{transmission-outer}). 

To prove part~\ref{Szego-theorem}, set 
\[
f(z)=\frac{4\prod_{j=1}^n\bigl(1-r_j^2\bigr)}{\bigl(\Psi_n^\ast(z)+\Phi_n^\ast(z)\bigr)^2},
\]
which is holomorphic on $\overline{\ddd}$, because $\Psi_n^\ast+\Phi_n^\ast$ is zero free there.   Evidently $f$ is itself zero free on $\overline{\ddd}$, and therefore $f=e^h$ for a function $h$ holomorphic on $\overline{\ddd}$.  Since $\Re h=\log|f|$ is harmonic, the mean value property for harmonic functions yields
\begin{equation}\label{outer-consequence}
\frac{1}{2\pi}\int_{-\pi}^\pi\log \bigl|f(e^{i\theta})\bigr|\,d\theta =\log |f(0)|=\sum_{j=1}^n\log\bigl(1-r_j^2\bigr),
\end{equation}
noting that $\Psi_n^\ast(0)=\Phi_n^\ast(0)=1$ (a fact equivalent to $\Psi_n$ and $\Phi_n$ being monic).  (Alternatively, in the terminology of Hardy space on the unit disk, $f$ is an outer function, and such functions are characterized by (\ref{outer-consequence}); see \cite[Thm.~2.7.10]{MaRo:2007}.) 
By parts~\ref{singular-harmonic-OPUC},\ref{g-equally-spaced-representation} and Corollary~\ref{cor-determinantal-formula},
\[
1-|R(\omega)|^2=\frac{4\Re\Psi_n^\ast\bigl(e^{2i\Delta\omega}\bigr)\overline{\Phi_n^\ast\bigl(e^{2i\Delta\omega}\bigr)}}{\bigl|\Psi_n^\ast\bigl(e^{2i\Delta\omega}\bigr)+\Phi_n^\ast\bigl(e^{2i\Delta\omega}\bigr)\bigr|^2}=\bigl|f\bigl(e^{2i\Delta\omega}\bigr)\bigr|. 
\]
Part \ref{Szego-theorem} then follows from (\ref{outer-consequence}) via the change variables $\theta=2\Delta\omega$. 
\end{pf}

\subsection{Arbitrarily-spaced jump points and OPUD\label{sec-OPUD}}

The representation of the generalized reflection coefficient $g_\zeta^\omega$ in terms of OPUC
in Proposition~\ref{prop-opuc} is valid only if the jump points of $\zeta\in\step_+(X)$ are equally spaced.  Two impedance functions with differing sequences of reflectivities lead to different polynomials $\Psi_n^\ast$ and $\Phi_n^\ast$.  Another, completely different, representation of $g_\zeta^\omega(0)$ in terms of orthogonal polynomials is presented in \cite{Gi:JFAA2017}, and the result extends easily to $g_\zeta^\omega$ itself.  The latter representation involves a special class of orthogonal polynomials on the unit disk (OPUD) called scattering polynomials, without any restriction on the spacing of the jump points of $\zeta$.  This single, universal family of bivariate orthogonal polynomials represents $g_\zeta^\omega$ for every $\zeta\in\step_+(X)$, for any $X$.  

Scattering polynomials are intimately related to the Riemannian structure on $\aut\mathbb{P}$ described in the introduction.  
Recall that $\aut\mathbb{P}$ consists of linear fractional transformations $\varphi_{\mu,\rho}$, where $\mu\in S^1$ and $\rho\in\ddd$, and so is represented by the manifold $S^1\times\ddd$, which carries a product metric comprised of the standard Euclidean metric on $S^1$ crossed with the metric 
\[
ds^2=\frac{4}{1-x^2-y^2}\left(dx^2+dy^2\right)
\]
on $\ddd$ (not to be confused with the hyperbolic metric); see \cite{EmSa:2012}.  The remarkable structure of the manifold $(\ddd,ds^2)$ is described in \cite{BaGi:Proc2017}.  Scattering polynomials are eigenfunctions of the Laplace-Beltrami operator associated to $ds^2$,
\begin{equation}\label{laplace-beltrami}
\Delta=-\frac{1-x^2-y^2}{4}\left(\frac{\partial^2}{\partial x^2}+\frac{\partial^2}{\partial y^2}\right).
\end{equation}
Here is a summary of the results needed from \cite{Gi:JFAA2017}.  

For each $(p,q)\in\integer^2$ define $\pi^{(p,q)}:\complex\rightarrow\complex$ as follows. If $\min\{p,q\}\geq 1$ set
\begin{equation}\label{scattering-polynomials}
\pi^{(p,q)}({z})=\textstyle\frac{(-1)^p}{q(p+q-1)!}\,\displaystyle(1-{z}\bar{{z}})\frac{\partial^{\,p+q}}{\partial\bar{{z}}^p\partial{z}^q}(1-{z}\bar{{z}})^{p+q-1}.
\end{equation}
If $\min\{p,q\}<0$ or $p=0<q$ set $\pi^{(p,q)}=0$; and if $p\geq 0$ set $\pi^{(p,0)}({z})={z}^p$.  Note in particular that $\pi^{(0,0)}=1$. The functions $\pi^{(p,q)}$ are polynomials in $z$ and $\overline{z}$, or equivalently in Euclidean coordinates $x$ and $y$, where $z=x+iy$. Furthermore, referring to (\ref{laplace-beltrami}),
\begin{equation}\label{eigenfunctions}
\Delta\pi^{(p,q)}=pq\,\pi^{(p,q)}\quad\mbox{ if }\quad \min\{p,q\}\geq 1\quad\mbox{ or if }\quad p\geq q=0.
\end{equation}
The operator $\Delta$ has no other eigenfunctions; see \cite{Gi:JFAA2017} and \cite{BaGi:Proc2017} for further details.  Let $\integer_+$ denote the nonnegative integers. 
\begin{thm}[{adapted from \cite[Thm.1, p.1497]{Gi:JFAA2017}}]\label{thm-fourier-series}
Fix $X=(x_0,x_1)$. Let $\zeta\in\step_+(X)$ have $n\geq 1$ jump points $y_j$, indexed according to their natural order, with corresponding reflectivities 
\[
r_j=\frac{\zeta(y_j-)-\zeta(y_j+)}{\zeta(y_j-)+\zeta(y_j+)}\qquad(1\leq j\leq n),
\]
and set
\[
\mathbf{y}=(y_1-x_0,y_2-y_1,\ldots,y_n-y_{n-1},x_1-y_n)\in\real^{n+1}_+\quad\mbox{ and }\quad r=(r_1,\ldots,r_n). 
\]
Then
\[
g_\zeta^\omega(\xi)=\sum_{k\in\{1\}\times\integer^n_+}\mathfrak{a}_k(r,\xi)e^{2i\langle k,\mathbf{y}\rangle\omega}\quad\mbox{ where }\quad\mathfrak{a}_k(r,\xi)=\left(\prod_{j=1}^{n}\pi^{(k_j,k_{j+1})}(r_j)\right)\xi^{k_{n+1}}.
\]
\end{thm} 
The interpretation $(r_j)^0=\xi^0=1$ is intended above, even if $r_j=0$ or $\xi=0$.  It follows from the definition of $\pi^{(p,q)}$ that the coefficient of $e^{2i\langle k,\mathbf{y}\rangle\omega}$ is non-zero for a given $k\in\{1\}\times\integer_+^n$ only if for every $1\leq j\leq n$, $k_{j+1}\neq0\Rightarrow k_j\neq0$.  In other words, $k$ consists of a block of positive entries followed by a block of zeros. 

Theorem~\ref{thm-fourier-series}, together with the fact that $g_\zeta^\omega$ is bounded, shows, independently of Proposition~\ref{prop-singular-harmonic}\ref{E-almost-periodic}, that with respect to $\omega$, $g_\zeta^\omega(\xi)$ is an almost periodic function in the sense of Besicovitch \cite{Be:1926}, having almost periods $\pi/\langle k,\mathbf{y}\rangle$, where $k\in\{1\}\times\integer_+^n$ ranges over the set $\mathfrak{K}$ of tuples consisting of a block of positive integers followed by a (possibly empty) block of zeros.  The almost periodic norm is
\begin{equation}\label{almost-periodic-norm}
\|g_\zeta(\xi)\|_{\rm ap}^2=\lim_{L\rightarrow\infty}\frac{1}{2L}\int_{-L}^L\left|g_\zeta^\omega(\xi)\right|^2\,d\omega\leq 1.
\end{equation}
In the generic case where the entries of $\mathbf{y}$ are linearly independent over the integers, 
\[
\langle k,\mathbf{y}\rangle=\langle \tilde{k},\mathbf{y}\rangle\quad\mbox{ only if }\quad k=\tilde{k},
\]
and so, by the Plancherel theorem for almost periodic function \cite[Ch.II,\S9]{Be:1955}, 
\begin{equation}\label{almost-periodic-norm-2}
\|g_\zeta(\xi)\|_{\rm ap}^2=
\sum_{k\in\mathfrak{K}}|\mathfrak{a}_k(r,\xi)|^2.
\end{equation}
In any case, (\ref{almost-periodic-norm}) shows that for $\zeta\in\step_+(X)$ the almost periodic structure of $g_\zeta^\omega(\xi)$ is determined by its high-frequency asymptotics as $|\omega|\rightarrow\infty$.  It will be shown later that this almost periodic structure is invariant under continuous perturbations of $\zeta$. 

Note the following immediate consequence of Theorem~\ref{thm-fourier-series}, under the same hypothesis. 
\begin{cor}\label{cor-lowest-frequency}
\[
g_\zeta^\omega(\xi)=r_1e^{2i(y_1-x_0)\omega}+
\sum_{\stackrel{k\in\{1\}\times\integer^n_+}{\langle k,\mathbf{y}\rangle>y_1-x_0}}
\mathfrak{a}_k(r,\xi)e^{2i\langle k,\mathbf{y}\rangle\omega}
\]
\end{cor}
\begin{pf}
The least possible value of $\langle k,\mathbf{y}\rangle$ for which $\mathfrak{a}_k(r,\xi)\neq 0$ is 
\[
\langle(1,0,\ldots,0),\mathbf{y}\rangle=y_1-x_0,
\]
in which case $\mathfrak{a}_k(r,\xi)=r_1$. And $\langle k,\mathbf{y}\rangle>y_1-x_0$ if $k\neq (1,0,\ldots,0)$ and $\mathfrak{a}_k(r,\xi)\neq 0$. 
\end{pf}

An obvious layer-stripping procedure results by combining Proposition~\ref{prop-step-composition} with the above corollary, since, for any $\xi\in\ddd$, 
\begin{equation}\label{layer-strip-frequency}
y_1-x_0=\frac{1}{2}\min\left\{\lambda\,\left|\,\lim_{L\rightarrow\infty}\frac{1}{2L}\int_{-L}^Lg_\zeta^\omega(\xi)e^{-i\lambda\omega}\,d\omega\neq 0\right.\right\}
\end{equation}
and
\begin{equation}\label{layer-strip-amplitude}
r_1=\lim_{L\rightarrow\infty}\frac{1}{2L}\int_{-L}^Lg_\zeta^\omega(\xi)e^{-2i(y_1-x_0)\omega}\,d\omega.
\end{equation}
(The above equations are valid in particular with $\xi=0$, so one can work with just the reflection coefficient $R(\omega)=g_{\zeta}^\omega(0)$.)  In detail, given 
\[
g_\zeta^\omega(\xi)=\varphi_{\mu_1,r_1}\circ\cdots\circ\varphi_{\mu_n,r_n}\bigl(e^{2i(x_1-y_n)\omega}\xi\bigr)
\]
as in Proposition~\ref{prop-step-composition}, one can determine $\mu_1=e^{2i(y_1-x_0)\omega}$ by (\ref{layer-strip-frequency}) and $r_1$ by (\ref{layer-strip-amplitude}) to yield $\varphi_{\mu_1,r_1}$.  The latter determines 
\[
\varphi_{\mu_1,r_1}^{-1}\circ g_\zeta^\omega(\xi)=\varphi_{\mu_2,r_2}\circ\cdots\circ\varphi_{\mu_n,r_n}\bigl(e^{2i(x_1-y_n)\omega}\xi\bigr).
\]
Iteration of the above procedure produces the sequences $r_1,\ldots,r_n$ and $y_1-x_0,\ldots,y_n-y_{n-1}$, allowing one to reconstruct $\zeta$ from $R$, provided $x_0$ and $\zeta(x_0+)$ are known.

\subsection{The regular harmonic exponential operator\label{sec-regular-harmonic}}

Recall the notation $\one$ for the function taking constant value 1, and let $L^1_\real(X)\subset L^1(X)
$ denote real-valued functions.  Let $C^1(\complex)$ denote the class of entire functions with its standard topology of uniform convergence on compact sets. 
Fix $X=(x_0,x_1)$ and $\alpha\in L^1_\real(X)
$. For each $\omega\in\complex$, define
\begin{equation}\label{A}
\omop:C(\overline{X})\rightarrow C(\overline{X}),\qquad \omop f(y)=\int_{x_0}^y\alpha(x)e^{2i(x-x_0)\omega}\overline{f(x)}\,dx.
\end{equation}
Note that $\omop$ is a conjugate linear Volterra operator, the standard analysis of which shows $1-\omop$ to be invertible. 
Define the harmonic exponential operator 
\[
E^{(x_0,y)}:L^1_\real(X)
\rightarrow C^1(\complex),\qquad \alpha\mapsto E^{(x_0,y)}_\alpha\qquad(x_0\leq y\leq x_1)
\]
by the formula
\begin{equation}\label{harmonic-exponential}
E^{(x_0,y)}_\alpha(\omega)=(1-\omop)^{-1}\one(y).
\end{equation}
Of course the assertion $E^{(x_0,y)}_\alpha\in C^1(\complex)$ has to be proved; this and other needed facts concerning the harmonic exponential are established in the next several propositions. 
\begin{prop}\label{prop-harmonic}
Fix $X=(x_0,x_1)$, $x_0<y\leq x_1$ and $\alpha\in L^1_\real(X)
$, and set 
\begin{equation}\label{kappa-simplex-definition}
\kappa(s_1,\ldots,s_j)=2\sum_{\nu=1}^j(-1)^{j-\nu}(s_\nu-x_0)\qquad(j\geq 1).
\end{equation}
Denote by $\simp^{(x_0,y)}_j$ the simplex
\[
\simp^{(x_0,y)}_j=\left\{(s_1,\ldots,s_j)\in\real^j\,|\,x_0<s_1<\cdots<s_j<y\right\}\qquad(j\geq 1),
\]
and let $\|\alpha\|$ denote the $L^1$ norm of $\alpha$. 
Then:\\[0pt]
\begin{enumerate}[label={(\roman*)},itemindent=0em]
\item \ $\omop^j\one(y)=\displaystyle\int\limits_{\simp^{(x_0,y)}_j}e^{i\omega\kappa(s_1,\ldots,s_j)}\prod_{\nu=1}^j\alpha(s_\nu)\,ds_1\cdots ds_j\;;\label{A-power-j}$\\[5pt]
\item \ $\bigl|\omop^j\one(y)\bigr|\leq M_{\omega}
\|\alpha\|^j/j!$ where $M_{\omega}
=\sup\left\{\left.|e^{it\omega}|\,\right|\,0<t<2(x_1-x_0)\right\}$;\label{A-power-bound}\\[5pt]
\item \ $\displaystyle
E^{(x_0,y)}_\alpha(\omega)=1+\displaystyle\sum_{j=1}^\infty\int\limits_{\simp^{(x_0,y)}_j}e^{i\omega\kappa(s_1,\ldots,s_j)}\prod_{\nu=1}^j\alpha(s_\nu)\,ds_1\cdots ds_j\; ;\label{E-series}
$\\[5pt]
\item \ $E^{(x_0,y)}_\alpha\in C^1(\complex)$;\label{E-entire}\\[5pt]
\item \ for any $\omega\in\complex$ and $\alpha_1,\alpha_2\in L^1_\real(X)
$, 
\[
\left|E_{\alpha_1}^{(x_0,y)}(\omega)-E_{\alpha_2}^{(x_0,y)}(\omega)\right|\leq M_\omega\|\alpha_1-\alpha_2\|\exp\left(\max\left\{\|\alpha_1\|,\|\alpha_2\|\right\}\right),
\]
where $M_{\omega}
=\sup\left\{\left.|e^{it\omega}|\,\right|\,0<t<2(x_1-x_0)\right\}$ as in \ref{A-power-bound};\label{E-difference-estimate}\\[5pt] 
\item the mapping $E^{(x_0,y)}:L^1_\real(X)
\rightarrow C^1(\complex)$ is continuous. \label{E-complex-continuous}
\end{enumerate}
\end{prop}
\begin{pf}
For $j\geq 1$, the $j$-fold iteration of (\ref{A}) yields part~\ref{A-power-j}.
Write $s=(s_1,\ldots,s_j)$, and observe that for $s\in\simp^{(x_0,y)}_j$, $0\leq\kappa(s)\leq 2(y-x_0)$, since 
\[
\textstyle\frac{1}{2}\kappa(s)=\displaystyle\left\{
\begin{array}{cc}
(s_j-s_{j-1})+(s_{j-2}-s_{j-3})+\cdots+(s_2-s_1)&\mbox{ if $j$ is even}\\
(s_j-s_{j-1})+(s_{j-2}-s_{j-3})+\cdots+(s_3-s_2)+(s_1-x_0)&\mbox{ if $j$ is odd}
\end{array}\right..
\]
It therefore follows from \ref{A-power-j} and the definition of $M_{\omega}
$ that
\[
\begin{split}
\left|\omop^j\one(y)\right|&\leq M_{\omega}
\negthickspace\negthickspace\int\limits_{\simp^{(x_0,y)}_j}\prod_{\nu=1}^j\left|\alpha(s_\nu)\right|\,ds_1\cdots ds_j\\
&=\frac{M_{\omega}
}{j!}\int\limits_{(x_0,y)^j}\prod_{\nu=1}^j\left|\alpha(s_\nu)\right|\,ds_1\cdots ds_j\\
&\leq\frac{M_{\omega}
\|\alpha\|^j}{j!}
\end{split}
\]
proving \ref{A-power-bound}.  Note that $M_{\omega}\geq 1$ varies continuously with respect to $\omega\in\complex$. Thus the right-hand series in part~\ref{E-series} converges absolutely to $(1-\omop)^{-1}\one(y)$, uniformly on compact sets in $\complex$, and 
\begin{equation}\label{E-bound}
\left|E^{(x_0,y)}_\alpha(\omega)\right|\leq M_{\omega}\exp\|\alpha\|\qquad(\omega\in\real).
\end{equation} 
To prove \ref{E-entire}, it suffices to observe that each component $\omop^j\one(y)$ is entire with respect to $\omega$ as a consequence of the integral formulation \ref{A-power-j}.  

Concerning \ref{E-difference-estimate}, let $\alpha_1,\alpha_2\in L^1_\real(X)
$ and $\omega\in\complex$, write $\omopone$ and $\omoptwo$ for the respective operators (\ref{A}), and set $\varepsilon=\alpha_1-\alpha_2$.  For $j\geq 1$ and $\omega\in\complex$,
\begin{equation}\label{alpha1-alpha2}
\begin{split}
\left|\omopone^j-\omoptwo^j\right|&=\left|\,\displaystyle\int\limits_{\simp^{(x_0,y)}_j}e^{i\omega\kappa(s_1,\ldots,s_j)}\left(\prod_{\nu=1}^j\alpha_1(s_\nu)-\prod_{\nu=1}^j\alpha_2(s_\nu)\right)\,ds_1\cdots ds_j\right|\\
&\leq M_\omega\int\limits_{\simp^{(x_0,y)}_j}\left|\prod_{\nu=1}^j\alpha_1(s_\nu)-\prod_{\nu=1}^j\alpha_2(s_\nu)\right|\,ds_1\cdots ds_j\\
&=\frac{M_\omega}{j!}\int\limits_{(x_0,y)^j}\left|\prod_{\nu=1}^j(\alpha_2+\varepsilon)(s_\nu)-\prod_{\nu=1}^j\alpha_2(s_\nu)\right|\,ds_1\cdots ds_j\\
&\leq\frac{M_\omega}{j!}j\int_{x_0}^y|\varepsilon|\left(\max\left\{\int_{x_0}^y|\alpha_1|,\int_{x_0}^y|\alpha_2|\right\}\right)^{j-1}\\
&=\frac{M_\omega}{(j-1)!}\|\alpha_1-\alpha_2\|\left(\max\left\{\|\alpha_1\|,\|\alpha_2\|\right\}\right)^{j-1}.
\end{split}
\end{equation}
Therefore
\[
\left|E_{\alpha_1}^{(x_0,y)}(\omega)-E_{\alpha_2}^{(x_0,y)}(\omega)\right|\leq M_\omega\|\alpha_1-\alpha_2\|\exp\left(\max\left\{\|\alpha_1\|,\|\alpha_2\|\right\}\right) 
\]
as desired, proving \ref{E-difference-estimate}.  

Lastly, given a convergent sequence $\alpha_n\rightarrow\alpha$ in $L^1_\real(X)$, the foregoing inequality implies $E^{(x_0,y)}_{\alpha_n}(\omega)\rightarrow E^{(x_0,y)}_\alpha(\omega)$ uniformly on compact subsets of $\complex$, since $M_\omega$ is continuous with respect to $\omega$. This proves \ref{E-complex-continuous}.
\end{pf}

Continuity of the harmonic exponential operator into the space of entire functions potentially brings complex analysis to bear on questions involving its values. However, applications to scattering only involve directly the restriction of $E^{(x_0,y)}_\alpha$ to the real line.  Thus in the next result, a given value $E^{(x_0,y)}_\alpha$ of the harmonic exponential will be regarded as a function of a real variable only.  Of course, being the restriction of an entire function, $E^{(x_0,y)}_\alpha$ necessarily belongs to $C^{\infty}(\real)$. 

\begin{prop}\label{prop-harmonic-2}
Fix $X=(x_0,x_1)$ and $x_0<y\leq x_1$.
\begin{enumerate}[label={(\roman*)},itemindent=0em]
\item Then $E^{(x_0,y)}_\alpha\in L^\infty(\real)\cap C^\infty(\real)$ for every $\alpha\in L^1_\real(X)
$, and the mapping 
\[
E^{(x_0,y)}:L^1_\real(X)
\rightarrow L^{\infty}(\real)
\]
is continuous;\label{E-continuous}\\[0pt]
\item for every $\alpha\in L^1_\real(X)$ and $\omega\in\real$, $\bigl|E^{(x_0,y)}_\alpha(\omega)\bigr|\leq e^{\|\alpha\|}$.\label{E-upper-bound}
\end{enumerate}
Let $\alpha\in L^1_\real(X)
$, and set $h=E^{(x_0,y)}_\alpha-1$.  Then:\\[0pt]
\begin{enumerate}[label={(\roman*)},itemindent=0em]
\setcounter{enumi}{2}
\item \  $\displaystyle\supp\rwc{h}\subset [0,2(y-x_0)]$ and $\rwc{h}\in L^1(\real);$\label{h-check-properties}\\[5pt]
\item \ $h(\omega)\rightarrow 0$ as $|\omega|\rightarrow\infty$;\label{h-decay}\\[5pt]
\item \ $h\in L^2(\real)$ if $\alpha\in L^2(X)$.\label{h-properties}\\[0pt]
\end{enumerate}
\end{prop}
\begin{pf}
First consider part \ref{h-check-properties}. Fix $j\geq 1$, and recall from the proof of Proposition~\ref{prop-harmonic} that $0\leq\kappa(s)\leq 2(y-x_0)$ for $s\in\simp^{(x_0,y)}_j$. For $0\leq t\leq 2(y-x_0)$, denote by $\mathcal{H}_t$ the hyperplane
\[
\mathcal{H}_t=\left\{\left.s\in\real^j\,\right|\,\kappa(s)=t\right\},
\]
and let $\mu(s)$ denote ordinary Lebesgue measure on $\mathcal{H}_t$.  Re-write the right-hand side of Proposition~\ref{prop-harmonic}\ref{A-power-j} as
\[
\begin{split}
\omop^j\one(y)&=\frac{1}{2\sqrt{j}}\int_{t=0}^{2(y-x_0)}\int_{\simp^{(x_0,y)}_j\cap\mathcal{H}_t}e^{i\omega\kappa(s_1,\ldots,s_j)}\prod_{\nu=1}^j\alpha(s_\nu)\,d\mu(s)dt \\
&=\int_{t=-\infty}^{\infty}\left(\frac{1}{2\sqrt{j}}\chi_{(0,2(y-x_0))}(t)\int_{\simp^{(x_0,y)}_j\cap\mathcal{H}_t}\prod_{\nu=1}^j\alpha(s_\nu)\,d\mu(s)\right)e^{i\omega t}dt \\
\end{split}
\]
and take the inverse Fourier transform with respect to $\omega\in\real$ to yield
\begin{equation}\label{A-j-Fourier}
\rwc{\left(\omop^j\one(y)\right)}(t)=\frac{1}{2\sqrt{j}}\chi_{(0,2(y-x_0))}(t)\int_{\simp^{(x_0,y)}_j\cap\mathcal{H}_t}\prod_{\nu=1}^j\alpha(s_\nu)\,d\mu(s),
\end{equation}
the interpretation of which in the case $j=1$ is 
\begin{equation}\label{A-Fourier}
\rwc{\left(\omop\one(y)\right)}(t)=\frac{1}{2}\chi_{(0,2(y-x_0))}(t)\alpha\bigl(\textstyle\frac{1}{2}t+x_0\bigr).
\end{equation}
The formulation (\ref{A-j-Fourier}) makes clear that $\supp\rwc{\left(\omop^j\one(y)\right)}\subset[0,2(y-x_0)]$
for each $j\geq 1$. 
By parts \ref{A-power-j} and \ref{E-series} of Proposition~\ref{prop-harmonic},
\[
h(\omega)=E^{(x_0,y)}_\alpha(\omega)-1=\sum_{j\geq 1}\omop^j\one(y).
\]
Therefore 
\[
\supp\rwc{h}\subset\bigcup_{j\geq 1}\supp\rwc{\left(\omop^j\one(y)\right)}\subset[0,2(y-x_0)],
\]
proving the first statement in part~\ref{h-check-properties}.  Next observe by (\ref{A-j-Fourier}) that 
\[
\begin{split}
\int_{\real}\left|\rwc{\omop^j\one(y)}(t)\right|dt&\leq \frac{1}{2\sqrt{j}}\int\limits_0^{2(y-x_0)}\int_{\simp^{(x_0,y)}_j\cap\mathcal{H}_t}\prod_{\nu=1}^j\left|\alpha(s_\nu)\right|\,d\mu(s)\,dt\\
&=\int_{\simp^{(x_0,y)}_j}\prod_{\nu=1}^j\left|\alpha(s_\nu)\right|\,ds\\
&=\frac{1}{j!}\int_{(x_0,y)^j}\prod_{\nu=1}^j\left|\alpha(s_\nu)\right|\,ds\\
&=\frac{1}{j!}\left(\int_{x_0}^y\left|\alpha\right|\right)^j\\
&\leq \|\alpha\|^j/j!
\end{split}
\]
from which it follows by Minkowski's inequality that 
\[
\int_\real\bigl|\rwc{h}\bigr|\leq\sum_{j\geq 1}\int_{\real}\left|\rwc{\omop^j\one(y)}\right|\leq\sum_{j\geq1}\|\alpha\|^j/j!=e^{\|\alpha\|}-1,
\]
proving $\rwc{h}\in L^1(\real)$.  This completes the proof of part \ref{h-check-properties}.  Part \ref{h-decay} follows from \ref{h-check-properties} by the Riemann-Lebesgue lemma.  It follows in turn that $E^{(x_0,y)}_\alpha=h+1$ is bounded, since $E^{(x_0,y)}_\alpha\in C^\infty(\real)$; thus $E^{(x_0,y)}_\alpha\in C^\infty(\real)\cap L^\infty(\real)$.  Continuity of the mapping $\alpha\mapsto E^{(x_0,y)}_\alpha$ with respect to the norms on $L^1_\real(X)
$ and $L^\infty(\real)$ then follows from Proposition~\ref{prop-harmonic}\ref{E-difference-estimate}, since the bound $M_\omega$ takes the constant value $M_\omega=1$ when $\omega\in\real$.   This proves \ref{E-continuous}; part~\ref{E-upper-bound} follows directly from (\ref{E-bound}). 

Lastly, to prove \ref{h-properties}, suppose $\alpha\in L^2(X)$, and note that 
\begin{multline}\nonumber
\left(\int_{\simp^{(x_0,y)}_j\cap\mathcal{H}_t}\prod_{\nu=1}^j\left|\alpha(s_\nu)\right|\,d\mu(s)\right)^2\leq\\
\vol\left(\simp^{(x_0,y)}_j\cap\mathcal{H}_t\right)\int_{\simp^{(x_0,y)}_j\cap\mathcal{H}_t}\prod_{\nu=1}^j\left|\alpha(s_\nu)\right|^2\,d\mu(s)
\end{multline}
by Jensen's inequality.
Geometric considerations give a bound on
$\vol\left(\simp^{(x_0,y)}_j\cap\mathcal{H}_t\right)$ as follows.  The volume of the intersection of a simplex with a hyperplane cannot exceed the volume of the largest face of the simplex.  The faces of the simplex $\simp^{(x_0,y)}_j$ either have volume $(y-x_0)^{j-1}/(j-1)!$ or $\sqrt{2}(y-x_0)^{j-1}/(j-1)!$.  Thus 
\begin{equation}\label{volume}
\vol\left(\simp^{(x_0,y)}_j\cap\mathcal{H}_t\right)\leq \sqrt{2}(y-x_0)^{j-1}/(j-1)!,
\end{equation}
and consequently
\begin{multline}\label{jensen}
\left(\int_{\simp^{(x_0,y)}_j\cap\mathcal{H}_t}\prod_{\nu=1}^j\left|\alpha(s_\nu)\right|\,d\mu(s)\right)^2\leq\\
\frac{\sqrt{2}(y-x_0)^{j-1}}{(j-1)!}\int_{\simp^{(x_0,y)}_j\cap\mathcal{H}_t}\prod_{\nu=1}^j\left|\alpha(s_\nu)\right|^2\,d\mu(s).
\end{multline}
By (\ref{A-j-Fourier}),
\[
\begin{split}
\int_{\real}\left|\rwc{\omop^j\one(y)}(t)\right|^2dt&\leq \frac{1}{4j}\int\limits_0^{2(y-x_0)}\left(\int_{\simp^{(x_0,y)}_j\cap\mathcal{H}_t}\prod_{\nu=1}^j\left|\alpha(s_\nu)\right|\,d\mu(s)\right)^2dt\\
&\leq\frac{\sqrt{2}(y-x_0)^{j-1}}{4j(j-1)!}\int\limits_0^{2(y-x_0)}\int_{\simp^{(x_0,y)}_j\cap\mathcal{H}_t}\prod_{\nu=1}^j\left|\alpha(s_\nu)\right|^2\,d\mu(s)\,dt\quad\mbox{by (\ref{jensen})}\\
&=\frac{(y-x_0)^{j-1}}{2\sqrt{j}j!(j-1)!}\int_{\simp^{(x_0,y)}_j}\prod_{\nu=1}^j\left|\alpha(s_\nu)\right|^2\,ds\\
&=\frac{(y-x_0)^{j-1}}{2\sqrt{j}j!(j-1)!}\left(\int_{x_0}^y\left|\alpha\right|^2\right)^j\\
&\leq \frac{\|\alpha\|^2_{L^2}}{2j\sqrt{j}}\;\frac{\left(\sqrt{x_1-x_0}\,\|\alpha\|_{L^2}\right)^{2(j-1)}}{(j-1)!^2}.
\end{split}
\]
Thus 
\begin{equation}\label{A-j-two-norm-estimate}
\|\rwc{\omop^j\one(y)}\|_{L^2}\leq \frac{\|\alpha\|_{L^2}}{\sqrt{2} j^{3/4}}\;\frac{\left(\sqrt{x_1-x_0}\,\|\alpha\|_{L^2}\right)^{j-1}}{(j-1)!},
\end{equation}
and hence Minkowski's inequality yields
\begin{equation}\label{minkowski-h-check}
\|\rwc{h}\|_{L^2}\leq \sum_{j\geq 1}\|\rwc{\omop^j\one(y)}\|_{L^2}\leq \frac{\|\alpha\|_{L^2}}{\sqrt{2}}\exp\left(\sqrt{x_1-x_0}\,\|\alpha\|_{L^2}\right),
\end{equation}
independently of $y\in(x_0,x_1]$.  Thus $\rwc{h}$, and hence $h$, belongs to $L^2(\real)$, proving \ref{h-properties}. 
\end{pf}

Of course since $h=E_\alpha^{(x_0,y)}-1$, part~\ref{h-properties} of Proposition~\ref{prop-harmonic} is equivalent to
\begin{cor}\label{cor-harmonic-exponential-bounded}
For every $\alpha\in L^1_\real(X)$, $E_\alpha^{(x_0,y)}(\omega)\rightarrow 1$ as $|\omega|\rightarrow\infty$. 
\end{cor}
Thus the harmonic exponential of a real-valued, integrable $\alpha$, $E^{(x_0,y)}_\alpha=1+h$, is a smooth, localized perturbation of the constant function $\one$.  
Note that (\ref{A-j-Fourier}) implies continuity with respect to $t\in\bigl(0,2(y-x_0)\bigr)$ of $\rwc{\left(\omop^j\one(y)\right)}(t)$ if $j\geq 2$, since the integral varies continuously in $t$.  By contrast, there is no integral if $j=1$; it then follows from (\ref{A-Fourier}) that $\rwc{\left(\omop\one(y)\right)}(t)$ has precisely the same discontinuities in $\bigl(0,2(y-x_0)\bigr)$ as $\alpha\bigl(\textstyle\frac{1}{2}t+x_0\bigr)$, with the same jumps. Thus $\rwc{h}(t)$ has precisely the same discontinuities as $\alpha\bigl(\textstyle\frac{1}{2}t+x_0\bigr)$, and it follows by the (continuity part of the) Riemann-Lebesgue lemma that $h\in L^1(\real)$ only if $\alpha$ is continuous.  So in general $h$ need not belong to $L^1(\real)$, a fact that limits the rate of decay of $h(\omega)$ as $|\omega|\rightarrow\infty$ when $\alpha$ has discontinuities.  The next proposition records the fact that, given $\zeta(x_0+)$, $E^{(x_0,y)}_\alpha$ encodes the value $\zeta(y)$ when $\alpha=-\zeta/(2\zeta)$.   
\begin{prop}\label{prop-harmonic-at-zero}
Let $X=(x_0,x_1)$, $\zeta\in C_+(X)$ and set $\alpha=-\zeta^\prime/(2\zeta)$.  If $0<y\leq x_1$,
\[
E^{(x_0,y)}_\alpha(0)=\exp\int_{x_0}^y\alpha=\sqrt{\zeta(x_0+)/\zeta(y)}\quad\mbox{ and }\quad
E^{(x_0,y)}_{-\alpha}(0)=\sqrt{\zeta(y)/\zeta(x_0+)}.
\]
\end{prop}
\begin{pf}
Setting $\omega=0$ in the right-hand side of part \ref{A-power-j} of Proposition~\ref{prop-harmonic} yields
\[
\int_{\simp^{(x_0,y)}_j}\prod_{\nu=1}^j\alpha(s_\nu)\,ds_1\cdots ds_j=\frac{1}{j!}\int_{x_0}^y\cdots\int_{x_0}^y\prod_{\nu=1}^j\alpha(s_\nu)\,ds_1\cdots ds_j=\frac{\left(\int_{x_0}^y\alpha\right)^j}{j!}.
\]
It then follows from the expansion Proposition~\ref{prop-harmonic}\ref{E-series} that 
\[
\begin{split}
E^{(x_0,y)}_\alpha(0)&=\exp\int_{x_0}^y\alpha\\
&=\exp\left(-\frac{1}{2}\left(\log\zeta(y)-\log\zeta(x_0+)\right)\right)\\
&=\sqrt{\zeta(x_0+)/\zeta(y)}.
\end{split}
\]
Replacing $\zeta$ by its reciprocal in the formula $\alpha=-\zeta^\prime/(2\zeta)$ corresponds to negation of $\alpha$, 
so $E^{(x_0,y)}_{-\alpha}(0)=\sqrt{\zeta(y)/\zeta(x_0+)}$. 
\end{pf}

\subsection{Singular approximation of the regular harmonic exponential\label{sec-approximation}}
Many different sequences of step functions converge uniformly to a given $\zeta\in\pwacplus(X)$. Here a standard choice of approximating sequence is fixed in order to reduce notational overhead and simplify presentation.  

Fix $X=(x_0,x_1)$ and $\zeta\in C_+(X)$.  For each $n\geq 1$, define the \emph{$n$th standard approximant} of $\zeta$ to be the step function $\zeta_n$ with $n$ equally-spaced jump points that interpolates $\zeta$ at the midpoints of successive jump points, as follows. Set 
\begin{equation}\label{zeta-n-spacing}
\Delta_n=\frac{x_1-x_0}{n+1}
\end{equation}
and 
\begin{equation}\label{zeta-n-jump-points}
y_{n,j}=x_0+j\Delta_n\qquad (0\leq j\leq n+1),
\end{equation}
so that, $y_{n,0}=x_0$ and $y_{n,n+1}=x_1$.  For each $x\in X$ and $0\leq j\leq n$, set 
\begin{equation}\label{zeta-n-formula}
\zeta_n(x)=\zeta(y_{n,j}+\Delta_n/2)\quad\mbox{ if }\quad y_{n,j}\leq x<y_{n,j+1}.
\end{equation}
Lastly, set 
\begin{subequations}\label{rnj}
\begin{align}
r_{n,j}&=\frac{\zeta_n(y_{n,j}-)-\zeta_n(y_{n,j}+)}{\zeta_n(y_{n,j}-)+\zeta_n(y_{n,j}+)\label{rnj-n}}\\[5pt]
&=\frac{\zeta(y_{n,j}-\Delta_n/2)-\zeta(y_{n,j}+\Delta_n/2)}{\zeta(y_{n,j}-\Delta_n/2)+\zeta(y_{n,j}+\Delta_n/2)}
\qquad(1\leq j\leq n).
\label{rnj-zeta}
\end{align}
\end{subequations}
The \emph{standard approximation} to $\zeta$ is defined to be the sequence $\zeta_n$ $(n\geq 1)$ of standard approximants. Since $\zeta$ is uniformly continuous, the standard approximation converges uniformly to $\zeta$.  

More generally, for any $\zeta\in\pwacplus(X)$, if $\zeta\not\in C_+(X)$, let $y_j$ denote the $k\geq 1$ points of discontinuity of $\zeta$ indexed according to their natural order,
\[
x_0<y_1<\cdots<y_k<x_1,
\]
and set $y_0=x_0$ and $y_{k+1}=x_1$.  
Let $\zeta_{n,j}$ denote the $n$th standard approximant of 
\[
\zeta|_{(y_{j-1},y_{j})}\in C_+\bigl((y_{j-1},y_{j})\bigr)\qquad(1\leq j\leq k+1),
\] 
and define the $n$th standard approximant of $\zeta$ to be the concatenation 
\begin{equation}\label{standard-approximant-concatenation}
\zeta_n=\zeta_{n,1}^{\rule{10pt}{0pt}\frown}\cdots^{\;\frown}\zeta_{n,k+1}.
\end{equation}
Considering $\zeta_n$ as an element of $L^\infty(X)$, its values at the $k$ points of discontinuity of $\zeta$ are immaterial; but for concreteness one can set $\zeta_n(y_j)=\zeta(y_j)$ $(1\leq j\leq k)$.
Define the standard approximation to $\zeta$ to be the sequence of standard approximants.  As before, the standard approximation converges uniformly to $\zeta$, by construction. 

The next results assemble needed estimates concerning the standard approximation.  Recall from \S\ref{sec-regulated} that $\alpha=-\zeta^\prime/(2\zeta)$ is by definition regulated, and hence bounded, if $\zeta\in C^1_+(X)$. 
\begin{prop}\label{prop-standard-approximation}
Let $X=(x_0,x_1)$ and $\zeta\in C^1_+(X)$. Set $\alpha=-\zeta^\prime/(2\zeta)$ and
\[
\alpha_{\max}=\sup\limits_{x\in X}|\alpha(x)|.
\]
Fix notation as in (\ref{zeta-n-spacing},\ref{zeta-n-jump-points},\ref{zeta-n-formula},\ref{rnj}), define $\kappa$ as in 
Proposition~\ref{prop-harmonic}, and denote by 
$\somopn$ the operator (\ref{Z}) 
corresponding to $\zeta_n$ $(n\geq 1)$. 
For every $n\geq 1$, there exist $n$ points 
\begin{equation}\nonumber
y_{n,j}^\ast\in(y_{n,j}-\Delta_n/2,y_{n,j}+\Delta_n/2)
\quad
\mbox{ with }
\quad
r_{n,j}=\tanh\left(\Delta_n\alpha(y_{n,j}^\ast)\right)\qquad(1\leq j\leq n)
\end{equation}
such that the following estimates hold. 
\begin{enumerate}[label={(\roman*)},itemindent=0pt]
\item \label{eta-estimate}  
Setting $
\eta_{n,j}=r_{n,j}-\Delta_n\alpha(y_{n,j}^\ast)
$,
\[
\bigl|r_{n,j}\bigr|\leq \Delta_n\alpha_{\max}\quad\mbox{ and }\quad
\bigl|\eta_{n,j}\bigr|\leq \frac{1}{3}\bigl(\Delta_n\alpha_{\max}\bigr)^3\qquad(1\leq j\leq n).
\] 
\item \label{Z-power-bound}For every $\omega\in\real$ and $k\geq 1$, $\bigl|\somopn^k\one(y)\bigr|\leq \bigl(\alpha_{\max}(x_1-x_0)\bigr)^k/k!\qquad(x_0<y\leq x_1)$.\\[0pt]
\item \label{varepsilon-estimate} For every $\omega\in\real$, the quantity
\[
\varepsilon_{n,j}=r_{n,j}e^{2i\omega(y_{n,j}-x_0)}-\Delta_n\alpha(y_{n,j}^\ast)e^{2i\omega(y_{n,j}^\ast-x_0)}
\]
satisfies
\[
\bigl|\varepsilon_{n,j}\bigr|\leq \Delta_n^2\alpha_{\max}B_{n,\omega}\qquad(1\leq j\leq n),
\]
where $B_{n,\omega}=|\omega|+\Delta_n\alpha_{\max}^2(1+\Delta_n)/3$.\\[0pt]
\item \label{r-prod-estimate} For every $1\leq k\leq n$, and every increasing sequence of $k$ jump points 
\[
y_{n,\nu_1}<\cdots<y_{n,\nu_k},
\] 
\begin{multline}\nonumber
\left|\exp\bigl(i\omega\kappa(y_{n,\nu_1},\ldots,y_{n,\nu_k})\bigr)\prod_{j=1}^k r_{n,\nu_j} -
\Delta_n^k\exp\bigl(i\omega\kappa(y_{n,\nu_1}^\ast,\ldots,y_{n,\nu_k}^\ast)\bigr)\prod_{j=1}^k\alpha\bigl(y_{n,\nu_j}^\ast\bigr)\right|\\
\leq
kB_{n,\omega}\alpha_{\max}^k\Delta_n^{k+1}\bigl(1+\Delta_nB_{n,\omega}\bigr)^{k-1}.
\end{multline}
\item \label{Z-k-estimate} For every $1\leq k\leq n$ and $y\in(x_0,x_1]$,
\begin{multline}\nonumber
\left|\somopn^k\one(y) -
\sum\limits_{y_{n,\nu_1}<\cdots<y_{n,\nu_k}<y}\Delta_n^k\exp\bigl(i\omega\kappa(y_{n,\nu_1}^\ast,\ldots,y_{n,\nu_k}^\ast)\bigr)\prod_{j=1}^k\alpha\bigl(y_{n,\nu_j}^\ast\bigr)\right|\\
\leq
\alpha_{\max}(x_1-x_0)\Delta_nB_{n,\omega}\frac{\bigl(\alpha_{\max}(x_1-x_0)(1+\Delta_nB_{n,\omega})\bigr)^{k-1}}{(k-1)!}.
\end{multline}
\end{enumerate}
\end{prop}
\begin{pf}
By (\ref{rnj-zeta}), 
\begin{equation}\label{rnj-in-terms-of-log}
\begin{split}
r_{n,j}&=\frac{\zeta(y_{n,j}-\Delta_n/2)-\zeta(y_{n,j}+\Delta_n/2)}{\zeta(y_{n,j}-\Delta_n/2)+\zeta(y_{n,j}+\Delta_n/2)}\\
&=\tanh\left(-\frac{1}{2}\log\zeta(y_{n,j}+\Delta_n/2)+\frac{1}{2}\log\zeta(y_{n,j}-\Delta_n/2)\right)
\qquad(1\leq j\leq n).
\end{split}
\end{equation}
Since $\alpha$ is the derivative of $-\frac{1}{2}\log\zeta$, the mean value theorem produces points 
\[
y_{n,j}-\Delta_n/2<y_{n,j}^\ast<y_{n,j}+\Delta_n/2
\]
at which  
\[
\alpha\bigl(y_{n,j}^\ast\bigr)=\frac{-\frac{1}{2}\log\zeta(y_{n,j}+\Delta_n/2)+\frac{1}{2}\log\zeta(y_{n,j}-\Delta_n/2)}{\Delta_n}\qquad\qquad(1\leq j\leq n). 
\]
Combined with (\ref{rnj-in-terms-of-log}) this yields 
\begin{equation}\label{r-as-tanh}
r_{n,j}=\tanh\left(\Delta_n\alpha(y_{n,j}^\ast)\right)\qquad(1\leq j\leq n). 
\end{equation}
Thus $|r_{n,j}|\leq \Delta_n\alpha_{\max}$, since $|\tanh\eta|\leq |\eta|$ if $\eta\in\real$. 
On the other hand, the inequality
\[
|\tanh \eta-\eta|\leq \bigl|\textstyle\frac{1}{3}\eta^3\bigr|\qquad(\eta\in\real)
\]
implies by (\ref{r-as-tanh}) that
\[
\bigl|\eta_{n,j}\bigr|=\bigl|\tanh\left(\Delta_n\alpha(y_{n,j}^\ast)\right)-\Delta_n\alpha(y_{n,j}^\ast)\bigr|\leq \frac{1}{3}\bigl(\Delta_n\alpha_{\max}\bigr)^3\qquad(1\leq j\leq n),
\]
proving part \ref{eta-estimate}.  Given that $\bigl|r_{n,j}\bigr|\leq \Delta_n\alpha_{\max}$, parts \ref{somop-power} and \ref{somop-zero} of Proposition~\ref{prop-singular-harmonic} imply 
\[
\bigl|\somopn^k\one(y)\bigr|\leq \binom{n}{k}\Delta_n^k\alpha_{\max}^k\leq\frac{n^k}{k!}\Delta_n^k\alpha_{\max}^k\leq\bigl(\alpha_{\max}(x_1-x_0)\bigr)^k/k!\qquad(k\geq1),
\]
since there are at most $\binom{n}{k}$ tuples satisfying $y_{\nu_1}<\cdots<y_{\nu_k}<y$, which proves part \ref{Z-power-bound} of the present proposition. 

It follows from the inequality $|1-e^{i\theta}|\leq |\theta|$ $(\theta\in\real)$ that 
\[
\bigl|e^{2i\omega(y_{n,j}-x_0)}-e^{2i\omega(y_{n,j}^\ast-x_0)}\bigr|\leq |2\omega(y_{n,j}^\ast-y_{n,j})|\leq |\omega|\Delta_n.
\]
Therefore
\[
\begin{split}
\bigl|\varepsilon_{n,j}\bigr|&=\left|r_{n,j}\left(e^{2i\omega(y_{n,j}-x_0)}-e^{2i\omega(y_{n,j}^\ast-x_0)}\right)+\left(r_{n,j}-\Delta_n\alpha(y_{n,j}^\ast)\right)e^{2i\omega(y_{n,j}^\ast-x_0)}\right|\\
&\leq |r_{n,j}\omega\Delta_n|+|\eta_{n,j}|\\
&=\bigl|\left(\Delta_n\alpha\bigl(y_{n,j}^\ast\bigr)+\eta_{n,j}\right)\omega\Delta_n\bigr|+|\eta_{n,j}|\quad\mbox{ since $r_{n,j}=\Delta_n\alpha\bigl(y_{n,j}^\ast\bigr)+\eta_{n,j}$}\\
&\leq |\omega|\Delta_n^2\alpha_{\max}+(|\omega|\Delta_n+1)|\eta_{n,j}|\\
&\leq\Delta_n^2\alpha_{\max}B_{n,\omega}\quad\mbox{ by \ref{eta-estimate}},
\end{split}
\] 
proving \ref{varepsilon-estimate}.

Let $\mathscr{C}$ denote the operation of complex conjugation.  To prove \ref{r-prod-estimate}, note by definition of $\kappa$ (\ref{kappa-simplex-definition}) that 
\[
\begin{split}
\exp\bigl(i\omega\kappa(y_{n,\nu_1},\ldots,y_{n,\nu_k})\bigr)\prod_{j=1}^k r_{n,\nu_j}
&=
\prod_{j=1}^k\mathscr{C}^{k-j}\bigl(r_{n,\nu_j}e^{2i\omega(y_{n,\nu_j}-x_0)}\bigr)\\
&=
\prod_{j=1}^k\mathscr{C}^{k-j}\bigl(\Delta_n\alpha(y_{n,\nu_j}^\ast)e^{2i\omega(y_{n,\nu_j}^\ast-x_0)}+\varepsilon_{n,\nu_j}\bigr),
\end{split}
\] 
and
\[
\Delta_n^k\exp\bigl(i\omega\kappa(y_{n,\nu_1}^\ast,\ldots,y_{n,\nu_k}^\ast)\bigr)\prod_{j=1}^k\alpha\bigl(y_{n,\nu_j}^\ast\bigr)
=
\prod_{j=1}^k\mathscr{C}^{k-j}\bigl(\Delta_n\alpha(y_{n,\nu_j}^\ast)e^{2i\omega(y_{n,\nu_j}^\ast-x_0)}\bigr).
\]
The algebraic identity
\[
\begin{split}
&\left|\prod_{j=1}^k(A_j+\epsilon_j)-\prod_{j=1}^kA_j\right|\leq k\epsilon_{\max}\left(A_{\max}+\epsilon_{\max}\right)^{k-1}\qquad(A_j,\epsilon_j\in\complex),\\[5pt]
&\quad\mbox{ where }\quad
A_{\max}=\sup\limits_{1\leq j\leq k}|A_j|\quad
\mbox{ and }\quad
\epsilon_{\max}=\sup\limits_{1\leq j\leq k}|\epsilon_j|,
\end{split}
\]
therefore implies 
\[
\begin{split}
\left|\exp\bigl(i\omega\kappa(y_{n,\nu_1},\ldots,y_{n,\nu_k})\bigr)\prod_{j=1}^k r_{n,\nu_j}\right. & -\left.
\Delta_n^k\exp\bigl(i\omega\kappa(y_{n,\nu_1}^\ast,\ldots,y_{n,\nu_k}^\ast)\bigr)\prod_{j=1}^k\alpha\bigl(y_{n,\nu_j}^\ast\bigr)\right|\\
&\leq k\sup_{1\leq j\leq k}|\varepsilon_{n,\nu_j}|\left(\Delta_n\alpha_{\max}+\sup_{1\leq j\leq k}|\varepsilon_{n,\nu_j}|\right)^{k-1}\\[5pt]
&\leq 
k\Delta_n^2\alpha_{\max}B_{n,\omega}\left(\Delta_n\alpha_{\max}+\Delta_n^2\alpha_{\max}B_{n,\omega}\right)^{k-1}
\quad\mbox{ by \ref{varepsilon-estimate}}.
\end{split}
\]
Part \ref{r-prod-estimate} follows directly.  

Applying Proposition~\ref{prop-singular-harmonic}\ref{somop-power},
\begin{equation}\label{expansion-of-Z-in-terms-of-r}
\begin{split}
&\somopn^k\one(y) -
\sum\limits_{y_{n,\nu_1}<\cdots<y_{n,\nu_k}<y}\Delta_n^k\exp\bigl(i\omega\kappa(y_{n,\nu_1}^\ast,\ldots,y_{n,\nu_k}^\ast)\bigr)\prod_{j=1}^k\alpha\bigl(y_{n,\nu_j}^\ast\bigr)\\
&=
\sum\limits_{y_{n,\nu_1}<\cdots<y_{n,\nu_k}<y}
\left(\exp\bigl(i\omega\kappa(y_{n,\nu_1},\ldots,y_{n,\nu_k})\bigr)\prod_{j=1}^k r_{n,\nu_j}-\right.\\
&\rule{1.5in}{0pt}
\left.\Delta_n^k\exp\bigl(i\omega\kappa(y_{n,\nu_1}^\ast,\ldots,y_{n,\nu_k}^\ast)\bigr)\prod_{j=1}^k\alpha\bigl(y_{n,\nu_j}^\ast\bigr)\right).
\end{split}
\end{equation}
For any $y\in(x_0,x_1]$, there are at most $\binom{n}{k}$ tuples such that $y_{n,\nu_1}<\cdots<y_{n,\nu_k}<y$.  Applying part~\ref{r-prod-estimate} to (\ref{expansion-of-Z-in-terms-of-r}), it therefore follows that 
\[
\begin{split}
\left|\rule{0pt}{22pt}\somopn^k\one(y)\right.&\left. -
\sum\limits_{y_{n,\nu_1}<\cdots<y_{n,\nu_k}<y}\Delta_n^k\exp\bigl(i\omega\kappa(y_{n,\nu_1}^\ast,\ldots,y_{n,\nu_k}^\ast)\bigr)\prod_{j=1}^k\alpha\bigl(y_{n,\nu_j}^\ast\bigr)\right|\\
&\leq\binom{n}{k}kB_{n,\omega}\alpha_{\max}^k\Delta_n^{k+1}\bigl(1+\Delta_nB_{n,\omega}\bigr)^{k-1}\\
&\leq \frac{n^k}{(k-1)!}B_{n,\omega}\alpha_{\max}^k\Delta_n^{k+1}\bigl(1+\Delta_nB_{n,\omega}\bigr)^{k-1}\\
&\leq \alpha_{\max}(x_1-x_0)\Delta_nB_{n,\omega}\frac{\bigl(\alpha_{\max}(x_1-x_0)(1+\Delta_nB_{n,\omega})\bigr)^{k-1}}{(k-1)!},
\end{split}
\]
proving \ref{Z-k-estimate}.
\end{pf}

The following consequence of Proposition~\ref{prop-standard-approximation} is noteworthy in view of Proposition~\ref{prop-opuc}\ref{Szego-theorem}.
\begin{cor}\label{cor-log-rnj-limit} 
Let $X=(x_0,x_1)$, $\zeta\in C^1_+(X)$, and set $\alpha=-\zeta^\prime/(2\zeta)$. Let $\zeta_n$ $(n\geq 1)$ be the regular approximation to $\zeta$, and denote by $r_{n,j}$ $(1\leq j\leq n)$ the reflectivities associated to $\zeta_n$.  Then: 
\begin{enumerate}[label={(\roman*)},itemindent=0pt]
\item\label{alpha-L2-convergence}\ 
$
\lim_{n\rightarrow\infty}\frac{1}{\Delta_n}\sum_{j=1}^n-\log\bigl(1-r_{n,j}^2\bigr)=\int\limits_X|\alpha|^2
$; and\\[5pt]
\item\label{rnj-square-limit}\  $\lim_{n\rightarrow\infty}\prod_{j=1}^n(1-r_{n,j}^2\bigr)=1.$
\end{enumerate}
\end{cor}
\begin{pf}
The expansion
\[
-\log(1-x)=x+\frac{x^2}{2}+\frac{x^3}{3}+\cdots
\]
yields that 
\[
r_{n,j}^2\leq -\log(1-r_{n,j}^2)\leq r_{n,j}^2+\frac{r_{n,j}^4}{2(1-r_{n,j})}.
\]
Proposition~\ref{prop-standard-approximation}\ref{eta-estimate} implies
\[
-\frac{1}{\Delta_n}\log(1-r_{n,j}^2)=
\Delta_n\alpha(y_{n,j}^\ast)^2+O(\Delta_n^3).
\]
Therefore 
\[
\frac{1}{\Delta_n}\sum_{j=1}^n-\log\bigl(1-r_{n,j}^2\bigr)=\left(\sum_{j=1}^n\Delta_n\alpha(y_{n,j}^\ast)^2\right) + O(\Delta_n^2),
\]
and thus
\[
\lim_{n\rightarrow\infty}\frac{1}{\Delta_n}\sum_{j=1}^n-\log\bigl(1-r_{n,j}^2\bigr)=\lim_{n\rightarrow\infty}\sum_{j=1}^n\Delta_n\alpha(y_{n,j}^\ast)^2=\int_X|\alpha|^2.
\]
It follows that 
\[
\lim_{n\rightarrow\infty}\sum_{j=1}^n-\log\bigl(1-r_{n,j}^2\bigr)=0,
\]
and hence, by exponentiating, that $\lim_{n\rightarrow\infty}\prod_{j=1}^n(1-r_{n,j}^2)=1$. 
\end{pf}

\begin{prop}\label{prop-step-exponential-bound}
Let $X=(x_0,x_1)$ and $\zeta\in C^1_+(X)$, with $\zeta_n\in\step_+(X)$ $(n\geq 1)$ the standard approximation to $\zeta$. 
Let $\mathscr{V}$ denote the total variation of $\frac{1}{2}\log\zeta$, and let $r_{n,j}$ $(1\leq j\leq n)$ denote the reflectivities corresponding to $\zeta_n$. Set
\[
\nu_n=\prod_{j=1}^n\sqrt{1-r_{n,j}^2}\quad(n\geq 1)\quad\mbox{ and }\quad \nu=\inf_{n\geq 1}\nu_n.
\]
Then $\nu>0$, and for every $n\geq 1$ and $\omega\in\real$,
\[
0<\nu e^{-\mathscr{V}}\leq \nu_n e^{-\mathscr{V}}\leq 
\min\left\{
\bigl|\calE_{\zeta_n}^X(\omega)\bigr|,\bigl|\calE_{1/\zeta_n}^X(\omega)\bigr|,\bigl|\calE_{\zeta_n}^X(\omega)+\calE_{1/\zeta_n}^X(\omega)\bigr|\right\}.
\]
\end{prop}
Thus each of $\calE_{\zeta_n}^X$, $\calE_{1/\zeta_n}^X$ and $\calE_{\zeta_n}^X+\calE_{1/\zeta_n}^X$ is bounded away from 0, independently of $n$. 
\begin{pf}
Set $R_n=g_{\zeta_n}(0)$ $(n\geq 1)$. Since the standard approximation is interpolating in the sense of Proposition~\ref{prop-bounded-variation}, the latter implies
\begin{equation}\nonumber
\left|2\Re\frac{R_n}{1-R_n}\right|\leq \frac{2\tanh\mathscr{V}}{1-\tanh\mathscr{V}}=e^{2\mathscr{V}}-1,
\end{equation}
while Corollary~\ref{cor-singular-real-part} implies
\begin{equation}\nonumber
2\Re\frac{R_n}{1-R_n}=\frac{\prod_{j=1}^n(1-r_{n,j}^2)}{\bigl|\calE_{1/\zeta_n}^X\bigr|^2}-1\quad(n\geq 1),
\end{equation}
where $r_{n,j}$ $(1\leq j\leq n)$ denote the reflectivities corresponding to $\zeta_n$.  Thus
\begin{equation}\nonumber
\left|\frac{\nu_n^2}{\bigl|\calE_{1/\zeta_n}^X\bigr|^2}-1\right|\leq e^{2\mathscr{V}}-1,
\end{equation}
and, since $\nu_n^2\rightarrow 1$ as $n\rightarrow\infty$ by Corollary~\ref{cor-log-rnj-limit}, 
it follows that $\nu>0$
and
\begin{equation}\label{singular-E-bound-1}
\bigl|\calE_{1/\zeta_n}^X\bigr|\geq\nu_n e^{-\mathscr{V}}\geq\nu e^{-\mathscr{V}}\qquad(n\geq 1). 
\end{equation}
The same argument goes through with $1/\zeta$ in place of $\zeta$, since the total variation of $\frac{1}{2}\log(1/\zeta)$ is the same as that of $\frac{1}{2}\log\zeta$, and the square reflectivities $r_{n,j}^2$ are the same for $1/\zeta_n$ as for $\zeta_n$.  Therefore 
\begin{equation}\label{singular-E-bound-2}
\bigl|\calE_{\zeta_n}^X\bigr|\geq\nu_n e^{-\mathscr{V}}\geq\nu e^{-\mathscr{V}}\qquad(n\geq 1)
\end{equation}
also.
Part \ref{singular-R-representation} of Theorem~\ref{thm-singular-representation} combined with Proposition~\ref{prop-bounded-variation} yields
\begin{equation}\nonumber
\left|\frac{\calE_{\zeta_n}^X-\calE_{1/\zeta_n}^X}{\calE_{\zeta_n}^X+\calE_{1/\zeta_n}^X}\right|=|R_n|\leq\tanh\mathscr{V}<1\qquad(n\geq 1).
\end{equation}
Thus
\[
2\bigl|\calE_{\zeta_n}^X\bigr|\leq \bigl|\calE_{\zeta_n}^X-\calE_{1/\zeta_n}^X\bigr|+\bigl|\calE_{\zeta_n}^X+\calE_{1/\zeta_n}^X\bigr|\leq \bigl(1+\tanh\mathscr{V}\bigr)\bigl|\calE_{\zeta_n}^X+\calE_{1/\zeta_n}^X\bigr|,
\]
whereby 
\begin{equation}\label{E-sum-lower-bound}
\nu_n e^{-\mathscr{V}}<\frac{2\nu_n e^{-\mathscr{V}}}{1+\tanh\mathscr{V}}\leq
\frac{2\bigl|\calE_{\zeta_n}^X\bigr|}{1+\tanh\mathscr{V}}
\leq\bigl|\calE_{\zeta_n}^X+\calE_{1/\zeta_n}^X\bigr|.
\end{equation}
This completes the proof.
\end{pf}

So much for technical estimates. Now comes a key singular approximation lemma. 
\begin{lem}[Second singular approximation lemma]\label{lem-singular-approximation}
Let $X=(x_0,x_1)$ and $\zeta\in C^1_+(X)$. Let $\zeta_n$ $(n\geq 1)$ be the standard approximation to $\zeta$, and set $\alpha=-\zeta^\prime/(2\zeta)$.  For every $y\in(x_0,x_1]$, 
\[
\calE_{\zeta_n}^{(x_0,y)}\rightarrow E_\alpha^{(x_0,y)}\quad\mbox{ as }\quad n\rightarrow\infty
\]
uniformly on compact sets. 
\end{lem}
\begin{pf}
Fix $\Omega>0$.  Given $\varepsilon>0$, choose $K\geq 1$ sufficiently large that 
\begin{equation}\label{epsilon-over-3-1}
2\sum_{k=K+1}^\infty\frac{\bigl(\alpha_{\max}(x_1-x_0)\bigr)^k}{k!}<\varepsilon/3. 
\end{equation}
Set 
\begin{multline}\label{cn-bound}
c_n=\alpha_{\max}(x_1-x_0)\Delta_nB_{n,\Omega}\exp\bigl(\alpha_{\max}(x_1-x_0)(1+\Delta_nB_{n,\Omega})\bigr)
\qquad(n\geq K),
\end{multline}
and note that $c_n\rightarrow 0$ as $n\rightarrow\infty$.  Observe by Proposition~\ref{prop-standard-approximation}\ref{Z-k-estimate} that if $-\Omega\leq\omega\leq\Omega$, then
\begin{equation}\label{epsilon-over-3-2a}
\sum_{k=1}^K
\left|\somopn^k\one(y) -
\sum\limits_{y_{n,\nu_1}<\cdots<y_{n,\nu_k}<y}\Delta_n^k\exp\bigl(i\omega\kappa(y_{n,\nu_1}^\ast,\ldots,y_{n,\nu_k}^\ast)\bigr)\prod_{j=1}^k\alpha\bigl(y_{n,\nu_j}^\ast\bigr)\right|
\leq c_n.
\end{equation}
Choose $N_1\geq K$ sufficiently large that $c_n\leq\varepsilon/3$ for every $n\geq N_1$.

Uniform continuity of $\alpha$ on $(x_0,x_1)$ implies that for each $k\geq1$ the family of uniformly continuous functions 
$
e^{i\omega\kappa(y_1,\ldots,y_k)}\prod_{j=1}^k\alpha(y_j), 
$
indexed by $\omega\in[-\Omega,\Omega]$, is equicontinuous on $\simp_k^{(x_0,y)}$.  Comparing Proposition~\ref{prop-harmonic}\ref{A-power-j} with Proposition~\ref{prop-standard-approximation}\ref{Z-k-estimate}, note that 
\begin{equation}\label{Riemann-approximant}
\sum\limits_{y_{n,\nu_1}<\cdots<y_{n,\nu_k}<y}\Delta_n^k\exp\bigl(i\omega\kappa(y_{n,\nu_1}^\ast,\ldots,y_{n,\nu_k}^\ast)\bigr)\prod_{j=1}^k\alpha\bigl(y_{n,\nu_j}^\ast\bigr)
\end{equation}
is a Riemann sum approximant to 
\begin{equation}\label{Riemann-approximand}
\omop^k\one(y)=\displaystyle\int\limits_{\simp^{(x_0,y)}_k}e^{i\omega\kappa(s_1,\ldots,s_k)}\prod_{\nu=1}^j\alpha(s_\nu)\,ds_1\cdots ds_k.
\end{equation}
By equicontinuity, there exists an index $M_k$, independent of $\omega\in[-\Omega,\Omega]$, such that for every $n\geq M_k$
\begin{equation}\label{epsilon-over-3-3a}
\left|\omop^k\one(y) -
\sum\limits_{y_{n,\nu_1}<\cdots<y_{n,\nu_k}<y}\Delta_n^k\exp\bigl(i\omega\kappa(y_{n,\nu_1}^\ast,\ldots,y_{n,\nu_k}^\ast)\bigr)\prod_{j=1}^k\alpha\bigl(y_{n,\nu_j}^\ast\bigr)\right|\leq \frac{\varepsilon}{3K}. 
\end{equation}

Set $N_2=\max\{N_1,M_1,\ldots,M_K\}$. Then, for every $n\geq N_2$,  
\begin{equation}\label{epsilon-over-three-two-thirds}
\begin{split}
&\left|\sum_{k=1}^K\somopn^k\one(y) -\omop^k\one(y)\right|\\
&\rule{18pt}{0pt}\leq 
\sum_{k=1}^K
\left|\somopn^k\one(y) -
\sum\limits_{y_{n,\nu_1}<\cdots<y_{n,\nu_k}<y}\Delta_n^k\exp\bigl(i\omega\kappa(y_{n,\nu_1}^\ast,\ldots,y_{n,\nu_k}^\ast)\bigr)\prod_{j=1}^k\alpha\bigl(y_{n,\nu_j}^\ast\bigr)\right|\\
&\rule{54pt}{0pt}+\left|\omop^k\one(y) -
\sum\limits_{y_{n,\nu_1}<\cdots<y_{n,\nu_k}<y}\Delta_n^k\exp\bigl(i\omega\kappa(y_{n,\nu_1}^\ast,\ldots,y_{n,\nu_k}^\ast)\bigr)\prod_{j=1}^k\alpha\bigl(y_{n,\nu_j}^\ast\bigr)\right|\\
&\rule{18pt}{0pt}\leq 2\varepsilon/3\quad\mbox{ by (\ref{epsilon-over-3-2a}), the choice of $N_1$, and (\ref{epsilon-over-3-3a})}.
\end{split}
\end{equation}
Therefore if $n\geq N_2$,
\[
\begin{split}
\left|\calE_{\zeta_n}^{(x_0,y)}(\omega)- E_\alpha^{(x_0,y)}(\omega)\right|&=
\left|\sum_{k=1}^\infty\somopn^k\one(y) -\omop^k\one(y)\right|\quad\mbox{ by Prop.~\ref{prop-singular-harmonic}\ref{singular-exponential-finite} and \ref{prop-harmonic}\ref{E-series}}\\
&\leq\left|\sum_{k=1}^K\somopn^k\one(y) -\omop^k\one(y)\right|
+\sum_{k=K+1}^\infty\left|\somopn^k\one(y)\right| +\left|\omop^k\one(y)\right|\\
&\leq 2\varepsilon/3+\varepsilon/3=\varepsilon\quad\mbox{ by (\ref{epsilon-over-three-two-thirds}), (\ref{epsilon-over-3-1}), Prop.~\ref{prop-standard-approximation}\ref{Z-power-bound} and Prop.~\ref{prop-harmonic}\ref{A-power-bound}}, 
\end{split}
\]
for every $\omega\in[-\Omega,\Omega]$.  Hence $\calE_{\zeta_n}^{(x_0,y)}\rightarrow E_\alpha^{(x_0,y)}$ uniformly on compact sets as $n\rightarrow\infty$. 
\end{pf}

\begin{prop}\label{prop-E-bound-from-zero}
Fix $X=(x_0,x_1)$ and let $\alpha\in L_\real^1(X)$ have norm $\|\alpha\|$. For every $\omega\in\real$,
\[
e^{-\|\alpha\|}\leq
\min\left\{
\bigl|E_\alpha^X(\omega)\bigr|,\bigl|E_{-\alpha}^X(\omega)\bigr|,\bigl|E_{\alpha}^X(\omega)+E_{-\alpha}^X(\omega)\bigr|\right\}.
\]
\end{prop}
\begin{pf}
Suppose first that $\alpha$ is continuous, and define $\zeta\in C^1_+(X)$ by 
\[
\zeta(x)=e^{-2\int_{x_0}^x\alpha}\qquad(x_0<x<x_1),
\]
so that the total variation $\mathscr{V}$ of $\frac{1}{2}\log\zeta$ is precisely $\|\alpha\|$.   
Let $\zeta_n$ $(n\geq 1)$ be the standard approximation to $\zeta$, with associated reflectitivities $r_{n,j}$ $(1\leq j\leq n)$. By Proposition~\ref{prop-step-exponential-bound},
\begin{equation}\label{En-bound}
\nu_n e^{-\mathscr{V}}\leq 
\min\left\{
\bigl|\calE_{\zeta_n}^X(\omega)\bigr|,\bigl|\calE_{1/\zeta_n}^X(\omega)\bigr|,\bigl|\calE_{\zeta_n}^X(\omega)+\calE_{1/\zeta_n}^X(\omega)\bigr|\right\}\qquad(n\geq1),
\end{equation}
where $\nu_n=\prod_{j=1}^n\sqrt{1-r_{n,j}^2}$. 
Letting $n\rightarrow\infty$, $\nu_n\rightarrow 1$ by Corollary~\ref{cor-log-rnj-limit}, and 
\[
\left\{
\bigl|\calE_{\zeta_n}^X(\omega)\bigr|,\bigl|\calE_{1/\zeta_n}^X(\omega)\bigr|,\bigl|\calE_{\zeta_n}^X(\omega)+\calE_{1/\zeta_n}^X(\omega)\bigr|\right\}\rightarrow
\left\{
\bigl|E_\alpha^X(\omega)\bigr|,\bigl|E_{-\alpha}^X(\omega)\bigr|,\bigl|E_{\alpha}^X(\omega)+E_{-\alpha}^X(\omega)\bigr|\right\}
\]
by Lemma~\ref{lem-singular-approximation}, whereby (\ref{En-bound}) yields
\begin{equation}\label{foregoing}
e^{-\|\alpha\|}\leq 
\min\left\{
\bigl|E_\alpha^X(\omega)\bigr|,\bigl|E_{-\alpha}^X(\omega)\bigr|,\bigl|E_{\alpha}^X(\omega)+E_{-\alpha}^X(\omega)\bigr|\right\}
\end{equation}
as desired. 

Suppose next that $\alpha\in L^1_\real(X)\setminus C_\real(X)$.  By density of $C_\real(X)$ in $L_\real^1(X)$ there exists a sequence $\alpha_n\in C_\real(X)$ $(n\geq 1)$ such that $\alpha_n\rightarrow\alpha$ in $L^1(X)$ as $n\rightarrow\infty$. 
By (\ref{foregoing}), 
\begin{equation}\label{succeeding}
e^{-\|\alpha_n\|}\leq 
\min\left\{
\bigl|E_{\alpha_n}^X(\omega)\bigr|,\bigl|E_{-\alpha_n}^X(\omega)\bigr|,\bigl|E_{\alpha_n}^X(\omega)+E_{-\alpha_n}^X(\omega)\bigr|\right\}\qquad(n\geq 1). 
\end{equation}
Part \ref{E-continuous} of Proposition~\ref{prop-harmonic-2} implies $E_{\alpha_n}^X\rightarrow E_{\alpha}^X$ and $E_{-\alpha_n}^X\rightarrow E_{-\alpha}^X$ uniformly.  Therefore letting $n\rightarrow\infty$ in (\ref{succeeding}) yields 
\[
e^{-\|\alpha\|}\leq 
\min\left\{
\bigl|E_\alpha^X(\omega)\bigr|,\bigl|E_{-\alpha}^X(\omega)\bigr|,\bigl|E_{\alpha}^X(\omega)+E_{-\alpha}^X(\omega)\bigr|\right\},
\]
completing the proof. 
\end{pf}

\section{Scattering for $\zeta\in\pwacplus(X)$\label{sec-scattering-for-piecewise-continuous}}

\subsection{The hyperbolic tangent operator and forward scattering\label{sec-forward-scattering}}

Given $X=(x_0,x_1)$, define the hyperbolic tangent operator 
\[
\Th^{(x_0,y)}:L^1_\real(X)\rightarrow C^\infty(\real)\cap L^\infty(\real),\qquad \alpha\mapsto \Th^{(x_0,y)}_\alpha\qquad(x_0< y\leq x_1)
\]
by the formula
\begin{equation}\label{harmonic-exponential}
\Th^{(x_0,y)}_\alpha=\frac{E^{(x_0,y)}_\alpha-E^{(x_0,y)}_{-\alpha}}{E^{(x_0,y)}_\alpha+E^{(x_0,y)}_{-\alpha}}.
\end{equation}
Boundedness and boundedness away from 0 of $E^{(x_0,y)}_\alpha,E^{(x_0,y)}_{-\alpha}$ and $E^{(x_0,y)}_\alpha+E^{(x_0,y)}_{-\alpha}$ (Propositions~\ref{prop-harmonic-2}\ref{E-upper-bound} and \ref{prop-E-bound-from-zero}) and continuity of $E^{(x_0,y)}$ (Proposition~\ref{prop-harmonic-2}\ref{E-continuous}) imply that $\Th^{(x_0,y)}$ is continuous with respect to the standard norms on $L^1_\real(X)$ and $L^\infty(\real)$. 

\begin{thm}\label{thm-hyperbolic}
Let $X=(x_0,x_1)$, $\zeta\in C_+(X)$, and set $\alpha=-\zeta^\prime/(2\zeta)$. Let $y\in(x_0,x_1]$, and set $R=g_{\zeta|_{(x_0,y)}}(0)$. Then:
\begin{enumerate}[label={(\roman*)},itemindent=0pt]
\item \label{R-Th}\  $R=\Th_\alpha^{(x_0,y)}$; \\[0pt]
\item \label{Th-bound} for every $\omega\in\real$, 
\[
\bigl|\Th_\alpha^{(x_0,y)}(\omega)\bigr|\leq \tanh\int_{x_0}^y|\alpha|\;;
\]
\item \label{g-regular-representation} $g_{\zeta|_{(x_0,y)}}^\omega(\xi)$ is $C^\infty$ with respect to $\omega$, and
\[
g_{\zeta|_{(x_0,y)}}^\omega(\xi)=\mu\frac{\xi+\rho}{1+\overline{\rho}\xi}\qquad(\xi\in\ddd),
\]
where 
\[
\mu=e^{2i(y-x_0)\omega}\frac{\overline{E_{\alpha}^{(x_0,y)}(\omega)+E_{-\alpha}^{(x_0,y)}(\omega)}}{E_{\alpha}^{(x_0,y)}(\omega)+E_{-\alpha}^{(x_0,y)}(\omega)}\quad\mbox{ and }\quad 
\rho=e^{-2i(y-x_0)\omega}\frac{E_{\alpha}^{(x_0,y)}(\omega)-E_{-\alpha}^{(x_0,y)}(\omega)}{\overline{E_{\alpha}^{(x_0,y)}(\omega)+E_{-\alpha}^{(x_0,y)}(\omega)}};
\]
\item \label{Bessel-inequality} if $\alpha\in L^2(X)$, then 
\[
\int_{-\infty}^\infty-\log\left(1-\bigl|\Th_{\alpha}^{(x_0,y)}(\omega)\bigr|^2\right)d\omega\leq \pi\int_{x_0}^y|\alpha|^2.
\]
\end{enumerate}
\end{thm}
\begin{pf}
The first step is to prove part~\ref{g-regular-representation}. Suppose $\zeta\in C^1_+(X)$ to begin with, so that $\alpha$ is continuous and its $L^1$ norm $\|\alpha\|$ coincides with the total variation $\mathscr{V}$ of $\frac{1}{2}\log\zeta$. Let $\zeta_n$ $(n\geq 1)$ be the standard approximation to $\zeta|_{(x_0,y)}$.  By Theorem~\ref{thm-singular-representation}, with $y$ replacing $x_1$ and $\zeta_n$ in place of $\zeta$, 
\[
g_{\zeta_n}^\omega(\xi)=\mu_n\frac{\xi+\rho_n}{1+\overline{\rho_n}\xi}\qquad(\xi\in\ddd),
\]
where 
\begin{equation}\label{mu-rho-step}
\mu_n=e^{2i(y-x_0)\omega}\frac{\overline{\calE_{\zeta_n}^{(x_0,y)}(\omega)+\calE_{1/\zeta_n}^{(x_0,y)}(\omega)}}{\calE_{\zeta_n}^{(x_0,y)}(\omega)+\calE_{1/\zeta_n}^{(x_0,y)}(\omega)}\quad\mbox{ and }\quad 
\rho_n=e^{-2i(y-x_0)\omega}\frac{\calE_{\zeta_n}^{(x_0,y)}(\omega)-\calE_{1/\zeta_n}^{(x_0,y)}(\omega)}{\overline{\calE_{\zeta_n}^{(x_0,y)}(\omega)+\calE_{1/\zeta_n}^{(x_0,y)}(\omega)}}.
\end{equation}
Lemma~\ref{lem-singular-approximation} implies $\calE_{\zeta_n}^{(x_0,y)}\rightarrow E^{(x_0,y)}_\alpha$ and $\calE_{1/\zeta_n}^{(x_0,y)}\rightarrow E^{(x_0,y)}_{-\alpha}$ uniformly on compact sets. 
Together with the fact that $E^{(x_0,y)}_\alpha+E^{(x_0,y)}_{-\alpha}$ is bounded away from 0, by Proposition~\ref{prop-E-bound-from-zero}, this implies
\begin{equation}\label{parameter-convergence}
\mu_n\rightarrow\mu\quad\mbox{ and }\quad \rho_n\rightarrow\rho\quad\mbox{ as }\quad n\rightarrow\infty
\end{equation}
uniformly on compact sets relative to $\omega$, whence
\[
\mu_n\frac{\xi+\rho_n}{1+\overline{\rho_n}\xi}\quad\rightarrow\quad\mu\frac{\xi+\rho}{1+\overline{\rho}\xi}\qquad(\xi\in\ddd). 
\]
On the other hand, $g_{\zeta_n}^\omega\rightarrow g_{\zeta|_{(x_0,y)}}^\omega$ by Theorem~\ref{thm-continuous}, since the standard approximation $\zeta_n$ to $\zeta|_{(x_0,y)}$ converges uniformly.  The representation 
\begin{equation}\label{g-representation}
g_{\zeta|_{(x_0,y)}}^\omega(\xi)=\mu\frac{\xi+\rho}{1+\bar{\rho}\xi}
\end{equation}
in part~\ref{g-regular-representation} follows.

Given the formula (\ref{g-representation}) for $\zeta\in C^1_+(X)$, Proposition~\ref{prop-harmonic-2}\ref{E-continuous} facilitates its extension to $\zeta\in C_+(X)$ as follows. Given $\zeta\in C_+(X)$, let $\alpha_n\in C_\real(X)$ $(n\geq 1)$ converge to $\alpha=-\zeta^\prime/(2\zeta)$ in $L^1_\real(X)$, and set
\[
\zeta_n(x)=\zeta(x_0+)e^{-2\int_{x_0}^x\alpha_n}\qquad (x_0<x<x_1).
\]
Observe that $\zeta_n\rightarrow\zeta$ uniformly as $n\rightarrow\infty$, so that $g_{\zeta_n|_{(x_0,y)}}^\omega\rightarrow g_{\zeta|_{(x_0,y)}}^\omega$ for almost every $\omega$, by Theorem~\ref{thm-continuous}.  On the other hand, with the notation
\[
\mu_n=e^{2i(y-x_0)\omega}\frac{\overline{E_{\alpha_n}^{(x_0,y)}(\omega)+E_{-\alpha_n}^{(x_0,y)}(\omega)}}{E_{\alpha_n}^{(x_0,y)}(\omega)+E_{-\alpha_n}^{(x_0,y)}(\omega)}\quad\mbox{ and }\quad 
\rho_n=e^{-2i(y-x_0)\omega}\frac{E_{\alpha_n}^{(x_0,y)}(\omega)-E_{-\alpha_n}^{(x_0,y)}(\omega)}{\overline{E_{\alpha_n}^{(x_0,y)}(\omega)+E_{-\alpha_n}^{(x_0,y)}(\omega)}},
\] 
(\ref{parameter-convergence}) follows from continuity of $E^{(x_0,y)}$ (Proposition~\ref{prop-harmonic-2}\ref{E-continuous}) together with the bound 
\begin{equation}\label{E-bound-from-zero-2}
\bigl|E_{\alpha}^{(x_0,y)}(\omega)+E_{-\alpha}^{(x_0,y)}(\omega)\bigr|\geq e^{-\|\alpha\|}
\end{equation}
(again, by Proposition~\ref{prop-E-bound-from-zero}).  Therefore $g_{\zeta_n|_{(x_0,y)}}^\omega$ converges to the formula on the right-hand side of (\ref{g-representation}), and thus (\ref{g-representation}) holds in the general case $\zeta\in C_+(X)$. That $g_{\zeta|_{(x_0,y)}}^\omega(\xi)$ is $C^\infty$ with respect to $\omega$ follows from the property $E_\alpha^{(x_0,y)}\in C^\infty(\real)$ in conjunction with (\ref{E-bound-from-zero-2}), completing the proof of part~\ref{g-regular-representation}.  Setting $\xi=0$ in \ref{g-regular-representation} yields \ref{R-Th}.  

To prove \ref{Th-bound}, suppose first that $\zeta\in C^1_+(X)$, so $\alpha\in C(X)$, and let $\zeta_n$ $(n\geq 1)$ be the standard approximation to $\zeta|_{(x_0,y)}$ as in the previous argument.  Consider the inequality (\ref{g-bound}) from the proof of Proposition~\ref{prop-bounded-variation} (interpreted with $x_1$ in the statement of the proposition replaced by $y$), 
\[
\bigl|g_{\zeta_n}^\omega(\xi)\bigr|\leq \tanh\left(\tanh^{-1}|\xi|+\textstyle\sum_{j=1}^n\tanh^{-1}|r_{n,j}|\right).
\]
Note that $\tanh|x|=|\tanh x|$ if $x\in\real$. By Proposition~\ref{prop-standard-approximation}, 
\[
|r_{n,j}|\leq\tanh\left(\Delta_n\bigl|\alpha(y_{n,j}^\ast)\bigr|\right)
\]
for some $y_{n,j}^\ast\in(y_{n,j}-\Delta_n/2,y_{n,j}+\Delta_n/2)$ $(1\leq j\leq n)$, where $y_{n,j}$ denote the evenly spaced jump points of $\zeta_n$.  Therefore 
\[
\bigl|g_{\zeta_n}^\omega(\xi)\bigr|\leq \tanh\left(\tanh^{-1}|\xi|+\textstyle\sum_{j=1}^n\Delta_n\bigl|\alpha(y_{n,j}^\ast)\bigr|\right).
\]
Since $|\alpha|$ is Riemann integrable, taking limits as $n\rightarrow\infty$ yields
\[
\left|g_{\zeta|_{(x_0,y)}}^\omega(\xi)\right|\leq \tanh\left(\tanh^{-1}|\xi|+\int_{x_0}^y|\alpha|\right).
\]
Given \ref{R-Th}, the inequality in part~\ref{Th-bound} follows upon setting $\xi=0$.  It extends to the case of arbitrary $\zeta\in C_+(X)$ by continuity of the operator $\Th^{(x_0,y)}$ on $L^1_\real(X)$.

To prove \ref{Bessel-inequality}, first suppose $\zeta\in C^1_+(X)$ and moreover that $\alpha=-\zeta^\prime/(2\zeta)$ is bounded.  Let $\zeta_n$ $(n\geq 1)$ be the standard approximation to $\zeta|_{(x_0,y)}$, with $r_{n,j}$ $(1\leq j\leq n)$ denoting the reflectivities associated to $\zeta_n$.  Set $R_n=g_{\zeta_n}(0)$, with $\mu_n$  and $\rho_n$ as in (\ref{mu-rho-step}), for $(n\geq 1)$.  The convergence result (\ref{parameter-convergence}) implies 
\begin{equation}\label{uniform-Rn-convergence}
R_n(\omega)=\mu_n\rho_n\rightarrow \mu\rho=\Th_{\alpha}^{(x_0,y)}(\omega)
\end{equation}
uniformly on compact sets.  
Let $L>0$ be arbitrary.  Since (\ref{uniform-Rn-convergence}) holds uniformly on compacts sets,
\[
\begin{split}
\int_{-L}^L-\log\bigl(1-\bigl|\Th_{\alpha}^{(x_0,y)}(\omega)\bigr|^2\bigr)d\omega
&=\lim_{n\rightarrow\infty}
\int_{-L}^L-\log\bigl(1-\bigl|R_n(\omega)\bigr|^2\bigr)d\omega\\
&\leq \lim_{n\rightarrow\infty}
\int\limits_{-\frac{\pi}{2\Delta_n}}^{\frac{\pi}{2\Delta_n}}-\log\left(1-\bigl|R_n(\omega)\bigr|^2\right)\,d\omega
\quad\mbox{ since $\Delta_n\rightarrow 0$ }\\
&=\lim_{n\rightarrow\infty}
\frac{\pi}{\Delta_n}\sum_{j=1}^n-\log(1-|r_{n,j}|^2)\quad\mbox{by Prop.~\ref{prop-opuc}\ref{Szego-theorem}}\\
&=\pi\int_{x_0}^y|\alpha|^2\quad\mbox{by Cor.~\ref{cor-log-rnj-limit}}.
\end{split}
\]
The inequality in part~\ref{Bessel-inequality} then follows from the fact that $L$ is arbitrary.  

Next suppose $\alpha\in L^2_\real(X)$ is arbitrary.  By density of $C(X)\cap L^\infty_\real(X)$ in $L^2_\real(X)$ there exist a sequence $\alpha_n\in C(X)\cap L^\infty_\real(X)$ $(n\geq 1)$ that converges to $\alpha$ in $L^2_\real(X)$, and hence also in $L^1_\real(X)$.  Again, let $L>0$ be arbitrary.  Letting $\mathscr{V}$ denote total variation of $\frac{1}{2}\log\zeta$, which is finite because $\zeta\in C_+(X)$, note that 
\[
\int_{-L}^L-\log\bigl(1-\bigl|\Th_{\alpha}^{(x_0,y)}(\omega)\bigr|^2\bigr)d\omega\leq\int_{-L}^L-\log\bigl(1-\bigl|\tanh\mathscr{V}\bigr|^2\bigr)d\omega<\infty
\]
since $R(\omega)=\Th_{\alpha}^{(x_0,y)}(\omega)\leq\tanh\mathscr{V}$, by Proposition~\ref{prop-bounded-variation}. Therefore 
\[
\int_{-L}^L-\log\bigl(1-\bigl|\Th_{\alpha_n}^{(x_0,y)}(\omega)\bigr|^2\bigr)d\omega\quad\rightarrow\quad
\int_{-L}^L-\log\bigl(1-\bigl|\Th_{\alpha}^{(x_0,y)}(\omega)\bigr|^2\bigr)d\omega,
\]
since $\Th_{\alpha_n}^{(x_0,y)}\rightarrow\Th_\alpha^{(x_0,y)}$ uniformly, by continuity of $\Th^{(x_0,y)}$. It was proved earlier that
\[
\int_{-L}^L-\log\bigl(1-\bigl|\Th_{\alpha_n}^{(x_0,y)}(\omega)\bigr|^2\bigr)d\omega\leq \pi\int_{x_0}^y|\alpha_n|^2
\]
for each $n\geq 1$.  Letting $n\rightarrow\infty$ yields
\[
\int_{-L}^L-\log\bigl(1-\bigl|\Th_{\alpha}^{(x_0,y)}(\omega)\bigr|^2\bigr)d\omega\leq \pi\int_{x_0}^y|\alpha|^2,
\]
and part \ref{Bessel-inequality} follows since $L>0$ was arbitrary. 
\end{pf}

By Corollary~\ref{cor-harmonic-exponential-bounded}, $E_{\pm\alpha}^{(x_0,y)}(\omega)\rightarrow 1$ as $|\omega|\rightarrow\infty$. Therefore $\mu\rightarrow e^{2i(y-x_0)\omega}$ and $\rho\rightarrow 0$ as $|\omega|\rightarrow\infty$ by part~\ref{g-regular-representation} of the above theorem, and so 
\[
\lim_{|\omega|\rightarrow\infty}e^{-2i(y-x_0)\omega}g_{\zeta|_{(x_0,y)}}^\omega(\xi)=\xi.
\]
In other words,
\begin{cor}\label{cor-g-regular-asymptotic}
Let $X=(x_0,x_1)$, $\zeta\in C_+(X)$ and $y\in(x_0,x_1]$.  Then
\[
g_{\zeta|_{(x_0,y)}}^\omega(\xi)\sim e^{2i(y-x_0)\omega}\xi\quad\mbox{ as }\quad |\omega|\rightarrow\infty.
\]
\end{cor}
The analysis of $\rwc{h}$ in Proposition~\ref{prop-harmonic} based on the structure of $\rwc{\left(\omop^j\one(y)\right)}(t)$ in equations (\ref{A-j-Fourier}) and (\ref{A-Fourier}) can be applied also to $\rwc{R}$ by Theorem~\ref{thm-hyperbolic}\ref{R-Th}. The upshot is that if $\alpha$ is bounded, discontinuities of $\rwc{R}$ correspond precisely to those of $\alpha$, as follows. 
\begin{prop}\label{prop-continuity-of-data}
Fix $X=(x_0,x_1)$, and let $\zeta\in C_+(X)$. Suppose $\alpha=-\zeta^\prime/(2\zeta)$ belongs to $L^\infty(\real)$, set $R=g_\zeta(0)$, and define $\widetilde{\alpha}:(0,2(x_1-x_0))\rightarrow\real$ by the formula
\[
\widetilde{\alpha}(t)=\textstyle\frac{1}{2}\alpha\bigl(\textstyle\frac{1}{2}t+x_0\bigr).
\]  
Then $\rwc{R}-\widetilde{\alpha}$ is continuous on $(0,2(x_1-x_0))$.   
\end{prop}
\begin{pf}
Note first by parts \ref{A-power-j} and \ref{E-series} of Proposition~\ref{prop-harmonic} that $E_\alpha^X-E_{-\alpha}^X=2(\omop\one(x_1)+H)$, where
\[
H=\sum_{j=2}^{\infty}\omop^{2j-1}\one(x_1),
\]
and $E_\alpha^X+E_{-\alpha}^X=2(1+J)$, where 
\[
J=\sum_{j=1}^\infty\omop^{2j}\one(x_1).
\]
It follows by Theorem~\ref{thm-hyperbolic}\ref{R-Th} that 
\[
R=\omop\one(x_1)+H-J\frac{\omop\one(x_1)+H}{1+J}, 
\]
and hence that 
\begin{equation}\label{R-check-decomposition}
\rwc{R}=\rwc{\omop\one(x_1)}+\rwc{H}-\rwc{J}\ast\left(\frac{\omop\one(x_1)+H}{1+J}\right)^{\rule{-4pt}{0pt}\rwc{\rule{12pt}{0pt}}}.
\end{equation}
Setting $\alpha_{\max}=\esssup|\alpha|$, it follows from formulas (\ref{A-j-Fourier}) and (\ref{volume}) that 
\[
\left|  \rwc{\left(\omop^j\one(y)\right)}(t)   \right|\leq\frac{\alpha_{\max}}{\sqrt{2j}}\,\frac{\bigl((y-x_0)\alpha_{\max}\bigr)^{j-1}}{(j-1)!}.
\]
Therefore each of $\rwc{H}$ and $\rwc{J}$ is a uniformly convergent series of continuous functions and hence continuous.  The term on the right-hand side of (\ref{R-check-decomposition}) involving convolution is automatically continuous. The difference $\rwc{R}-\rwc{\omop\one(x_1)}$ is thus continuous on the interval $(0,2(x_1-x_0)$, forcing the points of discontinuity (and values of the corresponding jumps) of $\rwc{R}$ and $\rwc{\omop\one(x_1)}$ to coincide.  Since 
\[
\rwc{\omop\one(x_1)}(t)=\textstyle\frac{1}{2}\alpha\bigl(\frac{1}{2}t+x_0\bigr)\qquad(0<t<2(x_1-x_0)),
\]
the conclusion of the proposition follows. 
\end{pf}
 
The assertion $\rwc{R}\approx\widetilde{\alpha}$ is the classical Born approximation, and the quantity $\rwc{R}-\widetilde{\alpha}$ in the statement of the proposition is the \emph{Born residual}.  While the Born residual for bounded $\alpha$ is guaranteed to be continuous, it is in general non-zero, and there are examples where it is large---see, for example, Figure~\ref{fig-d} in \S\ref{sec-computation-example}. 
A further application of Lemma~\ref{lem-singular-approximation} establishes some technical facts needed later.
\begin{prop}\label{prop-preliminary-approximation-application}
Fix $X=(x_0,x_1)$, $\zeta\in\pwacplus(X)$, and set $\alpha=-\zeta^\prime/(2\zeta)$. Suppose that $\zeta$ is continuous on $(x_0,y)$ for a given $y\in (x_0,x_1]$, and set $X_1=(x_0,y)$.  Set $R_1=g_{\zeta|_{X_1}}(0)$ and write $\alpha_1=\alpha|_{X_1}$.
Then:\\[0pt]
\begin{enumerate}[label={(\roman*)},itemindent=0pt]
\item\label{real-part-unity} \ $\Re E^{(x_0,y)}_{\alpha_1}\overline{E^{(x_0,y)}_{-\alpha_1}}=1$;\\[5pt]
\item\label{inverse-fourier-support} the inverse Fourier transform of 
$
\displaystyle\frac{1}{\overline{E_{-\alpha_1}^{(x_0,y)}}}-1
$
is supported on $(-\infty,0]$;
\item \label{R1-formula}
\[
2\Re\frac{R_1}{1-R_1}=\frac{1}{\bigl|E_{-\alpha_1}^{(x_0,y)}\bigr|^2}-1.
\] 
\end{enumerate}
\end{prop}
\begin{pf}
Suppose first that $\zeta|_{X_1}\in C^1(X_1)$, and let $\zeta_n$ $(n\geq 1)$ denote the standard approximation to $\zeta|_{X_1}$.  Denote by $r_{n,j}$ $(1\leq j\leq n)$ the reflectivities corresponding to $\zeta_n$.  

Then by Corollary~\ref{cor-singular-real-part} and Lemma~\ref{lem-singular-approximation},
\[
\Re\left( \calE_{\zeta_n}^{X_1}(\omega)\,\overline{\calE_{1/\zeta_n}^{X_1}(\omega)}\right)=\prod_{j=1}^n(1-r_{n,j}^2)\rightarrow\Re E^{(x_0,y)}_{\alpha_1}\overline{E^{(x_0,y)}_{-\alpha_1}}\qquad(\omega\in\real),
\]
whereby $\Re E^{(x_0,y)}_{\alpha_1}\overline{E^{(x_0,y)}_{-\alpha_1}}$ is constant. 
Corollary~\ref{cor-harmonic-exponential-bounded} then implies $\Re E^{(x_0,y)}_{\alpha_1}\overline{E^{(x_0,y)}_{-\alpha_1}}=1$.  

Proposition~\ref{prop-opuc}\ref{reciprocal-transform-support} asserts
\begin{equation}\label{support-assertion}
\supp\displaystyle\left({1}/{\overline{\calE_{1/\zeta_n}^{(x_0,y)}}}\right)^{\rule{-4pt}{0pt}\rwc{\rule{12pt}{0pt}}}\subset(-\infty,0]\qquad(n\geq 1).
\end{equation}
Since $\calE_{1/\zeta_n}^{(x_0,y)}$ and $E_{-\alpha_1}^{(x_0,y)}$ are bounded away from zero, Lemma~\ref{lem-singular-approximation} implies
\[
\frac{1}{\overline{\calE_{1/\zeta_n}^{(x_0,y)}}}\rightarrow \frac{1}{\overline{E_{-\alpha_1}^{(x_0,y)}}}\quad\mbox{ as }\quad n\rightarrow\infty
\]
uniformly on compact sets.  Uniform convergence on compact sets implies convergence in the sense of tempered distributions.  The inverse Fourier transform is a continuous operator on the space of tempered distributions, so 
\[
\left({1}/{\overline{\calE_{1/\zeta_n}^{(x_0,y)}}}\right)^{\rule{-4pt}{0pt}\rwc{\rule{12pt}{0pt}}}\rightarrow \left({1}/{\overline{E_{-\alpha_1}^{(x_0,y)}}}\right)^{\rule{-4pt}{0pt}\rwc{\rule{12pt}{0pt}}}\quad\mbox{ as }\quad n\rightarrow\infty
\]
in the sense of distributions.  Evaluation of the distributions $\left({1}/{\overline{\calE_{1/\zeta_n}^{(x_0,y)}}}\right)^{\rule{-4pt}{0pt}\rwc{\rule{12pt}{0pt}}}$ on a test functions supported on $(0,\infty)$ yields 0, by (\ref{support-assertion}).  Therefore the same is true of the limiting distribution $\left({1}/{\overline{E_{-\alpha_1}^{(x_0,y)}}}\right)^{\rule{-4pt}{0pt}\rwc{\rule{12pt}{0pt}}}$. In other words, the support of  $\left({1}/{\overline{E_{-\alpha_1}^{(x_0,y)}}}\right)^{\rule{-4pt}{0pt}\rwc{\rule{12pt}{0pt}}}$ is contained in $(-\infty,0]$, and hence so is that of 
\[
\left({1}/{\overline{E_{-\alpha_1}^{(x_0,y)}}}\right)^{\rule{-4pt}{0pt}\rwc{\rule{12pt}{0pt}}}-\delta=
 \left(\frac{1}{\overline{E_{-\alpha_1}^{(x_0,y)}}}-1\right)^{\rule{-4pt}{0pt}\rwc{\rule{12pt}{0pt}}}.
 \]
Both of the foregoing results extend to the more general case $\zeta|_{X_1}\in C(X_1)$ by continuity of the harmonic exponential with respect to $\alpha_1\in L^1_\real(\real)$ and density of $C_\real(X_1)\cap L^\infty(X_1)$ in $L_\real^1(X_1)$, proving \ref{real-part-unity} and \ref{inverse-fourier-support}. 

Part \ref{R-Th} of Theorem~\ref{thm-hyperbolic}, part \ref{real-part-unity} of the present proposition and the identity
\[
2\Re\frac{R_1}{1-R_1}=\frac{\Re E^{(x_0,y)}_{\alpha_1}\overline{E^{(x_0,y)}_{-\alpha_1}}}{\bigl|E_{-\alpha_1}^{(x_0,y)}\bigr|^2}-1
\]
imply part \ref{R1-formula}.
\end{pf}

Theorems~\ref{thm-g-composition} and \ref{thm-hyperbolic} combine to yield an explicit description of the generalized reflection coefficient corresponding to an arbitrary $\zeta\in\pwacplus(X)$ directly in terms of $\zeta$.  Of course if $\zeta\in C_+(X)$, then Theorem~\ref{thm-hyperbolic}  already gives such a description, so it suffices to consider the case where $\zeta$ has at least one point of discontinuity. 

For absolutely continuous $\zeta$, the distributional derivative $\zeta^\prime$ corresponds to a regular function, and so the formulas $\alpha=-\zeta^\prime/(2\zeta)$ and $\alpha=-\bigl(\frac{1}{2}\log\zeta\bigr)^\prime$ are equivalent.  On the other hand, given $X=(x_0,x_1)$, if $\zeta\in\pwacplus(X)$ has discontinuities, then only the latter formula makes sense, and it determines a distribution $\alpha=-\bigl(\frac{1}{2}\log\zeta\bigr)^\prime$ whose singular part is supported at the jump points of $\zeta$.  It will be useful later to refer to the decomposition 
\begin{equation}\label{reg-sing-decomposition}
\alpha=\alpha_{\regular}+\alpha_{\sing},
\end{equation}
in which $\alpha_{\regular}\in L_\real^1(X)$ and $\alpha_{\sing}$ is a linear combination of Dirac delta distributions.  

\begin{thm}\label{thm-forward-scattering}
Given $X=(x_0,x_1)$, let $\zeta\in\pwacplus(X)$ and set $\alpha=-\bigl(\textstyle\frac{1}{2}\log\zeta\bigr)^\prime$. Suppose $\zeta\not\in C_+(X)$, and denote by $y_j$ $(1\leq j\leq n)$ the points of discontinuity of $\zeta$, indexed according to their natural order,
\[
x_0<y_1<\cdots<y_n<x_1,
\]
with $y_0=x_0$ and $y_{n+1}=x_1$.  Write 
\[
r_j=\frac{\zeta(y_j-)-\zeta(y_j+)}{\zeta(y_j-)+\zeta(y_j+)}\qquad(1\leq j\leq n). 
\]
For each $1\leq j\leq n+1$, set $\alpha_j=\alpha|_{(y_{j-1},y_j)}$, $\mu_j=e^{2i(y_j-y_{j-1})\omega}$,  
\[
\tilde\mu_j=\mu_j\frac{\overline{E_{\alpha_j}^{(y_{j-1},y_j)}(\omega)+E_{-\alpha_j}^{(y_{j-1},y_j)}(\omega)}}{E_{\alpha_j}^{(y_{j-1},y_j)}(\omega)+E_{-\alpha_j}^{(y_{j-1},y_j)}(\omega)}\quad\mbox{ and }\quad
\rho_j=\overline{\mu_j}\frac{E_{\alpha_j}^{(y_{j-1},y_j)}(\omega)-E_{-\alpha_j}^{(y_{j-1},y_j)}(\omega)}{\overline{E_{\alpha_j}^{(y_{j-1},y_j)}(\omega)+E_{-\alpha_j}^{(y_{j-1},y_j)}(\omega)}}.
\]
Then $g_\zeta^\omega(\xi)$ is $C^\infty$ with respect to $\omega$, and:
\begin{enumerate}[label={(\roman*)},itemindent=0pt]
\item \label{forward-scattering-formula}
\[
g_\zeta^\omega=\varphi_{\tilde\mu_1,\rho_1}\circ\varphi_{1,r_1}\circ\cdots\circ\varphi_{\tilde\mu_n,\rho_n}\circ\varphi_{1,r_n}\circ\varphi_{\tilde\mu_{n+1},\rho_{n+1}};
\]
\item \label{forward-scattering-asymptotics}
\[
g_\zeta^\omega(\xi)\sim\varphi_{\mu_1,r_1}\circ\cdots\circ\varphi_{\mu_n,r_n}(\mu_{n+1}\xi)\quad\mbox{ as }\quad |\omega|\rightarrow\infty.
\]
\end{enumerate}
\end{thm}
\begin{pf}
Theorem~\ref{thm-hyperbolic}\ref{g-regular-representation} ensures $g_{\zeta_j}^\omega=\varphi_{\tilde\mu_j,\rho_j}$ is $C^\infty$ with respect to $\omega$, for each $1\leq j\leq n+1$.  Part~\ref{forward-scattering-formula} of the present theorem then follows by Theorem~\ref{thm-g-composition}.  Since each $\varphi_{\tilde\mu_j,\rho_j}(\xi)$ and $\varphi_{1,r_j}(\xi)$ is analytic with respect to $\xi$, it follows that the composition $g_\zeta^\omega$ is $C^\infty$ with respect to $\omega$. 

Corollary~\ref{cor-harmonic-exponential-bounded} implies $\tilde\mu_j\sim\mu_j$, $\rho_j\sim0$ and hence $\varphi_{\tilde\mu_j,\rho_j}\sim\varphi_{\mu_j,0}$ as $|\omega|\rightarrow\infty$ for each $1\leq j\leq n+1$.  The fact that 
\[
\varphi_{\mu_j,0}\circ\varphi_{1,r_j}=\varphi_{\mu_j,r_j}\qquad(1\leq j\leq n)
\]
then gives part~\ref{forward-scattering-asymptotics}.
\end{pf}

Note that in terms of the factorization $\zeta=\zeta_1\zeta_2$ of $\zeta\in\pwacplus(X)\setminus C_+(X)$, where $\zeta_1$ is continuous and $\zeta_2\in\step_+(X)$, part~\ref{forward-scattering-asymptotics} of Theorem~\ref{thm-forward-scattering} may be expressed as
\begin{equation}\label{forward-scattering-asymptotics-2}
g_\zeta^\omega(\xi)\sim g_{\zeta_2}^\omega(\xi)\quad\mbox{ as }\quad |\omega|\rightarrow\infty.
\end{equation}
(Cf.~Proposition~\ref{prop-step-composition}.)
\begin{cor}\label{cor-almost-periodic-computation}
Under the hypotheses of Theorem~\ref{thm-forward-scattering}: $r_1\neq 0$; for any $\xi\in\ddd$,
\begin{equation}\nonumber
y_1-x_0=\frac{1}{2}\min\left\{\lambda\,\left|\,\lim_{L\rightarrow\infty}\frac{1}{2L}\int_{-L}^Lg_\zeta^\omega(\xi)e^{-i\lambda\omega}\,d\omega\neq 0\right.\right\};
\end{equation}
and
\begin{equation}\nonumber
r_1=\lim_{L\rightarrow\infty}\frac{1}{2L}\int_{-L}^Lg_\zeta^\omega(\xi)e^{-2i(y_1-x_0)\omega}\,d\omega.
\end{equation}
\end{cor}
\begin{pf}
Note $r_1\neq 0$ because $y_1$ is a jump point of $\zeta$. 
Theorem~\ref{thm-forward-scattering}\ref{forward-scattering-asymptotics} implies that for any $\lambda\in\real$, 
\[
\lim_{L\rightarrow\infty}\frac{1}{2L}\int_{-L}^Lg_\zeta^\omega(\xi)e^{-i\lambda\omega}\,d\omega
=\lim_{L\rightarrow\infty}\frac{1}{2L}\int_{-L}^L\varphi_{\mu_1,r_1}\circ\cdots\circ\varphi_{\mu_n,r_n}(\mu_{n+1}\xi)e^{-i\lambda\omega}\,d\omega.
\]
The present result then follows from Corollary~\ref{cor-lowest-frequency} and the discussion following it. 
\end{pf}

\subsection{Short-range inversion and inverse scattering\label{sec-inverse-scattering}}

\begin{thm}[Short-range inversion formula]\label{thm-short-range-inversion}
Fix $X=(x_0,x_1)$, $\zeta\in\pwacplus(X)$, and set $\alpha=-\bigl(\textstyle\frac{1}{2}\log\zeta\bigr)^\prime$. Suppose for a given $y_1\in (x_0,x_1]$ that $\zeta$ is continuous on $(x_0,y_1)$ and $\int_{x_0}^{y_1}|\alpha|^2<\infty$.  Set 
\begin{equation}\label{gamma-constant}
\gamma=\frac{\left(1-\tanh\int_{x_0}^{y_1}|\alpha|\right)^2}{4\int_{x_0}^{y_1}|\alpha|^2}\quad\mbox{ and }\quad S(t)=-\left(\frac{R}{1-R}\right)^{\rule{-6pt}{0pt}\rwc{\rule{12pt}{0pt}}}\rule{-6pt}{0pt}(|t|)\qquad(t\in\real),
\end{equation}
where $R=g_\zeta(0)$.  If $x_0<y<\min\{x_0+\gamma,y_1\}$, then 
\begin{equation}
\zeta(y)=\zeta(x_0+)\left(1+\sum_{j=1}^\infty\;\int\limits_{(0,2(y-x_0))^j}\negthickspace S(t_j)S(t_j-t_{j-1})\cdots S(t_2-t_1)\,dt_1\cdots dt_j\right)^2.\label{magic-formula-2}
\end{equation}
\end{thm}
\begin{pf}
Recall the notation $\one$ for the function taking constant value 1.  
For $y\in(x_0,y_1]$, define
\[
\yop:L^2([0,2(y-x_0)])\rightarrow C([0,2(y-x_0)]),\quad \yop f(t)=\int\limits_{0}^{2(y-x_0)}S(t-s)f(s)\,ds.
\]
Iterating the latter yields  
\[
\yop^j\one(0)=\int\limits_{(0,2(y-x_0))^j}\negthickspace S(t_j)S(t_j-t_{j-1})\cdots S(t_2-t_1)\,dt_1\cdots dt_j\qquad(j\geq1).
\]
Thus, in order to prove (\ref{magic-formula-2}), it suffices to prove the equivalent formula
\begin{equation}
\zeta(y)=\zeta(x_0+)\left(\sum_{j=0}^\infty\yop^j\one(0)\right)^2\label{magic-formula-1},
\end{equation}
as follows. 

Given $y\in(x_0,y_1]$, set $\zeta_1=\zeta|_{(x_0,y)}$, $\alpha_1=\alpha|_{(x_0,y)}$ and $R_1=g_{\zeta_1}(0)$.  To begin, note that $\yop$ may be expressed in terms of $R_1$, as
\begin{equation}\label{yop-1}
\yop f(t)=-\int_{0}^{2(y-x_0)}\left(2\Re R_1/(1-R_1)\right)^{\rule{-4pt}{0pt}\rwc{\rule{12pt}{0pt}}}\rule{-2pt}{0pt}(t-s)f(s)\,ds\qquad 0\leq t\leq 2(y-x_0).
\end{equation}
Two observations justify this.  Firstly, $\rwc{R}_1$ is real-valued and supported on the positive half line $\real_+$.  Therefore the same is true of $(R_1/(1-R_1))^{\rule{-4pt}{0pt}\rwc{\rule{12pt}{0pt}}}\rule{-2pt}{0pt}$; the reflection of the latter about zero,
$(\overline{R}_1/(1-\overline{R}_1))^{\rule{-4pt}{0pt}\rwc{\rule{12pt}{0pt}}}\rule{-2pt}{0pt}$, is correspondingly supported on $\real_-$.  Thus, for $0<|\eta|< 2(y-x_0)$, 
\begin{equation}\label{K-kernel}
\begin{split}
K(\eta):=&\left(2\Re R_1/(1-R_1)\right)^{\rule{-4pt}{0pt}\rwc{\rule{12pt}{0pt}}}\rule{-2pt}{0pt}(\eta)\\
=&(R_1/(1-R_1))^{\rule{-4pt}{0pt}\rwc{\rule{12pt}{0pt}}}\rule{-2pt}{0pt}(\eta)+(\overline{R}_1/(1-\overline{R}_1))^{\rule{-4pt}{0pt}\rwc{\rule{12pt}{0pt}}}\rule{-2pt}{0pt}(\eta)\\
=&(R_1/(1-R_1))^{\rule{-4pt}{0pt}\rwc{\rule{12pt}{0pt}}}\rule{-2pt}{0pt}(|\eta|).
\end{split}
\end{equation}
The essential fact here is that the kernel $K$ is even. A second observation, which cinches (\ref{yop-1}), is that 
\[
R_1/(1-R_1))^{\rule{-4pt}{0pt}\rwc{\rule{12pt}{0pt}}}\rule{-2pt}{0pt}(|\eta|)=(R/(1-R))^{\rule{-4pt}{0pt}\rwc{\rule{12pt}{0pt}}}\rule{-2pt}{0pt}(|\eta|)
\]
if $|\eta|< 2(y-x_0)$, since by finite speed of propagation the right-hand expression depends only on $\zeta_1$. 

The formula for $\gamma>0$ in (\ref{gamma-constant}) is concocted to ensure that the $L^2$ operator norm of $\yop$ is less than $1$ if $x_0<y<x_0+\gamma$.  In detail,  
\[
\begin{split}
\|\yop\|^2\leq& \int_{0}^{2(y-x_0)}\int_{0}^{2(y-x_0)}K(t-s)^2\,ds dt\\
\leq& \int_{0}^{2(y-x_0)}\left(\int_{-\infty}^\infty K(t)^2\,dt\right)\,ds\\
=&2(y-x_0)\int_{-\infty}^\infty K(t)^2\,dt\\
=&\frac{(y-x_0)}{\pi}\int_{-\infty}^\infty \left(2\Re R_1(\omega)/(1-R_1(\omega))\right)^2\,d\omega\quad\mbox{ since $\|f\|_2^2=\textstyle\frac{1}{2\pi}\|\hat{f}\|_2^2$ }\\
\leq&\frac{4(y-x_0)}{\pi}\int_{-\infty}^\infty \left|R_1(\omega)/(1-R_1(\omega))\right|^2\,d\omega\\
\leq&\frac{4(y-x_0)}{\pi\left(1-\tanh\int_{x_0}^{y}|\alpha_1|\right)^2}\int_{-\infty}^\infty \left|R_1(\omega)\right|^2\,d\omega\quad\mbox{by Thm.~\ref{thm-hyperbolic}\ref{Th-bound}}\\
\leq&\frac{4(y-x_0)}{\pi\left(1-\tanh\int_{x_0}^{y}|\alpha_1|\right)^2}\int_{-\infty}^\infty-\log\left(1-|R_1(\omega)|^2\right)\,d\omega\\
\leq&\frac{4(y-x_0)}{\left(1-\tanh\int_{x_0}^{y}|\alpha_1|\right)^2}\int_{x_0}^{y}|\alpha_1|^2\quad\mbox{by Thm.~\ref{thm-hyperbolic}\ref{Bessel-inequality}}\\
\leq&\frac{4(y-x_0)}{\left(1-\tanh\int_{x_0}^{y_1}|\alpha|\right)^2}\int_{x_0}^{y_1}|\alpha|^2\quad\mbox{ since $\alpha_1=\alpha|_{(x_0,y)}$ }\\
=&(y-x_0)/\gamma. 
\end{split}
\]
Thus $\|\yop\|<1$ if $x_0<y<x_0+\gamma$, implying that the operator
\[
(1-\yop)^{-1}:L^2([0,2(y-x_0)])\rightarrow L^2([0,2(y-x_0)]),\quad (1-\yop)^{-1}=\sum_{j=0}^\infty\yop^j
\]
is bounded and the series representation converges. 

Verification of the formula (\ref{magic-formula-1}) rests on various facts established earlier concerning the harmonic exponential.  To begin, set 
\[
k=\widehat{K},\quad h=E_{-\alpha_1}^{(x_0,y)}-1,\quad f=\frac{1}{\overline{E_{-\alpha_1}^{(x_0,y)}}}-1,
\]
and note by Prop.~\ref{prop-preliminary-approximation-application}\ref{R1-formula} that 
\[
k=\frac{1}{\left|E_{-\alpha_1}^{(x_0,y)}\right|^2}-1.
\]
These functions are all in $L^2(\real)$: the foregoing analysis of $\yop$ demonstrates that $k\in L^2(\real)$; $h\in L^2(\real)\cap L^\infty(\real)$ by Prop.~\ref{prop-harmonic-2}\ref{h-properties}; that $f\in L^2(\real)$ then follows from the identity
\begin{equation}\label{khf-1}
kh+k+h=f. 
\end{equation}
Taking the inverse Fourier transform of (\ref{khf-1}) yields 
\begin{equation}\label{khf-2}
\rwc{h}+K\ast\rwc{h}=-K+\rwc{f}.
\end{equation}
The function $\rwc{h}$ is supported on $[0,2(y-x_0)]$ by Prop.~\ref{prop-harmonic-2}\ref{h-check-properties}.  Combining the identity 
\[
\rwc{h}+K\ast\rwc{h}=(1-\yop)\rwc{h}
\]
with equation (\ref{khf-2}) therefore implies
\[
\rwc{h}=(1-\yop)^{-1}(-K+\rwc{f})=\chi_{([0,2(y-x_0)])}(1-\yop)^{-1}(-K+\rwc{f}).
\]
Now, $\rwc{f}$ is supported on $\real_-$ by Prop.~\ref{prop-preliminary-approximation-application}\ref{inverse-fourier-support}, and $\rwc{f}\in L^2(\real)$ since $f\in L^2(\real)$.  Thus 
\[
\chi_{([0,2(y-x_0)])}(1-\yop)^{-1}\rwc{f}=0,
\]
and hence 
\begin{equation}\label{h-check}
\rwc{h}=-\chi_{([0,2(y-x_0)])}(1-\yop)^{-1}K.
\end{equation}
The goal is to evaluate $E_{-\alpha_1}^{(x_0,y)}(0)$ using the definition of $h$ and equation (\ref{h-check}), since
\begin{equation}\label{E-square-root-formula}
E_{-\alpha_1}^{(x_0,y)}(0)=\sqrt{\zeta(y)/\zeta(x_0+)}
\end{equation}
(see Prop.~\ref{prop-harmonic-at-zero}).  Since the kernel $K$ is even, it follows from (\ref{yop-1}) that 
\begin{equation}\label{yop-powers}
\begin{split}
&\yop^{j+1}\one(0)\\
&=(-1)^{j+1}\int_{0}^{2(y-x_0)}\cdots\int_{0}^{2(y-x_0)}K(s_{j+1})K(s_{j+1}-s_j)\cdots K(s_2-s_1)\,ds_1\ldots ds_{j+1}\\
&=-\int_{0}^{2(y-x_1)}\yop^jK(s)\,ds.
\end{split} 
\end{equation}
Expressing $E_{-\alpha_1}^{(x_0,y)}(0)$ in terms of $h(0)$,
\[
\begin{split}
E_{-\alpha_1}^{(x_0,y)}(0)=&1+h(0)\\
=&1+\int_{-\infty}^\infty \rwc{h}\\
=&1-\int_{0}^{2(y-x_1)}(1-\yop)^{-1}K(s)\,ds\\
=&1-\int_0^{2(y-x_1)}\sum_{j=0}^\infty\yop^jK(s)\,ds\\
=&1+\sum_{j=0}^\infty\yop^{j+1}\one(0)\quad\mbox{ by (\ref{yop-powers}) }\\
=&\sum_{j=0}^\infty\yop^j\one(0).  
\end{split}
\]
The desired formula (\ref{magic-formula-1}) then follows from (\ref{E-square-root-formula}).  
\end{pf} 
 
Formula (\ref{magic-formula-2}) 
expresses $\zeta(y)$ in terms of the data $R$, inverting the formula 
 $
 R=\Th^{X}_\alpha
$
at short range, i.e., at spatial locations $y$ sufficiently near the reference point $x_0$.  
Iteration of the short-range inversion formula yields the following.
\begin{cor}\label{cor-determinacy}
Fix $X=(x_0,x_1)$, $\zeta\in\pwacplus(X)$ and set $R=g_\zeta(0)$.  If $\zeta$ is continuous on $(x_0,y_1)$ for some $y_1\in X$ and $\int_{x_0}^{y_1}|\alpha|^2<\infty$, then $R$ together with $\zeta(x_0+)$ determines $\zeta|_{(x_0,y_1)}$.  
If in addition $y_1<x_1$ and $\zeta$ has a jump point at $y_1$, then $y_1-x_0$ and $\zeta(y_1+)$ are also determined. 
\end{cor}
\begin{pf}
Consider first the determination of $\zeta|_{(x_0,y_1)}$. If $\min\{x_0+\gamma,y_1\}=y_1$ there is nothing to prove, since the result follows directly from Theorem~\ref{thm-short-range-inversion}. Suppose on the other hand that $y_0:=x_0+\gamma<y_1$. Set 
\[
\gamma_1=\frac{\left(1-\tanh\int_{y_0}^{y_1}|\alpha|\right)^2}{4\int_{y_0}^{y_1}|\alpha|^2},
\]
and observe that $\gamma_1\geq \gamma$.  By Theorem~\ref{thm-short-range-inversion}, $R$ and $\zeta(x_0+)$ determine $g_{\zeta|_{(x_0,y_0)}}^\omega$ and also $\zeta(y_0)=\zeta(y_0-)$; these in turn determine 
\[
R_1(\omega)=g_{\zeta|_{(y_0,x_1)}}^\omega(0)=\bigl(g_{\zeta|_{(x_0,y_0)}}^\omega\bigr)^{-1}\circ g_{\zeta}^\omega(0).
\]
Theorem~\ref{thm-short-range-inversion} applied to $R_1$, with $\gamma_1$ in place of $\gamma$ and $y_0$ in place of $x_0$, determines $\zeta(y)$ for \[y_0<y<\min\{y_0+\gamma_1,y_1\}.\]  Thus $\zeta(y)$ is determined for all
\[
x_0<y<\min\{x_0+\gamma+\gamma_1,y_1\}.\]
Since $\gamma_1\geq\gamma$, at most $\lceil (y_1-x_0)/\gamma\rceil$ iterations of this process ultimately yield $\zeta|_{(x_0,y_1)}$. 

Suppose next that $y_1<x_1$ and $y_1$ is a jump point of $\zeta$, and let 
\[
r_1=\frac{\zeta(y_1-)-\zeta(y_1+)}{\zeta(y_1-)+\zeta(y_1+)}=\tanh\frac{1}{2}\log\bigl(\zeta(y_1-)/\zeta(y_1+)\bigr)
\]
denote the associated reflectivity.  The quantities $y_1-x_0$ and $r_1$ are determined by Corollary~\ref{cor-almost-periodic-computation}.  Since $\zeta(y_1-)$ is known from $\zeta|_{(x_0,y_1)}$, $\zeta(y_1+)$ is given by the formula
\[
\zeta(y_1+)=\zeta(y_1-)e^{-2\tanh^{-1}r_1}. 
\]
\end{pf}

Referring to the decomposition (\ref{reg-sing-decomposition}), in the case where $\alpha_{\regular}\in L^2(X)$, Corollary~\ref{cor-determinacy} sets up a punctuated layer stripping procedure whereby $\zeta$, having jump points $y_1<\cdots<y_n$, may be expressed in terms of its reflection coefficient
\begin{equation}\label{R-composite-form-1}
R=g_\zeta(0)=g_{\zeta_1}\circ\varphi_{1,r_1}\circ\cdots\circ g_{\zeta_n}\circ\varphi_{1,r_n}\circ g_{\zeta_{n+1}}(0).
\end{equation}
Here $r_j$ denotes the reflectivity at $y_j$ $(1\leq j\leq n)$; $y_0=x_0$, $y_{n+1}=x_1$; and $\zeta_j=\zeta|_{(y_{j-1},y_j)}$ for all $1\leq j\leq n+1$. In detail, the layer stripping alternates between the short-range inversion formula and the almost periodic transform used in the proof of Corollary~\ref{cor-determinacy} as follows.  Given $R$ and $\zeta(x_0+)$, the quantities $\zeta_1$, $r_1$ and $\zeta(y_1+)$ are determined as in Corollary~\ref{cor-determinacy}, in turn yielding 
\[
R_1=\varphi_{1,r_1}^{-1}\circ g_{\zeta_1}^{-1}(R)=g_{\zeta_2}\circ\varphi_{1,r_2}\circ\cdots\circ g_{\zeta_n}\circ\varphi_{1,r_n}\circ g_{\zeta_{n+1}}(0).
\]
Given $R_1$ and $\zeta(y_1+)$, the procedure may be repeated to yield $\zeta_j$, $r_j$, $\zeta(y_j+)$, and 
\[
R_j=\varphi_{1,r_j}^{-1}\circ g_{\zeta_j}^{-1}(R_{j-1})=g_{\zeta_{j+1}}\circ\varphi_{1,r_{j+1}}\circ\cdots\circ g_{\zeta_n}\circ\varphi_{1,r_n}\circ g_{\zeta_{n+1}}(0)
\]
for all $1\leq j\leq n$.  Lastly, $R_n$ and $\zeta(y_n+)$ determine $\zeta_{n+1}$ by Corollary~\ref{cor-determinacy}, and thus 
\[
\zeta=\zeta_1^{\;\frown}\cdots^{\;\frown}\zeta_{n+1}
\]
is completely determined by $R$.  This proves

\begin{thm}\label{thm-piecewise-injectivity}
Fix $X=(x_0,x_1)$ and set 
\[
\mathscr{C}(X)=\left\{\left.\zeta\in\pwacplus(X)\,\right|\, \zeta(x_0+)=1\mbox{ and } \alpha_{\regular}\in L^2(X)\mbox{ where }\alpha=-\bigl(\textstyle\frac{1}{2}\log\zeta\bigr)^\prime\right\}.
\]
The mapping 
\[
\mathscr{C}(X)\rightarrow L^\infty(\real)\cap C^\infty(\real),\qquad \zeta\mapsto g_\zeta(0)
\]
is injective. 
\end{thm}

\subsection{A singular trace formula\label{sec-renormalized-trace}}

The classical trace formula (\ref{classical-trace}) is not valid for $R_2=g_\zeta(0)$ when $\zeta\in\step_+(X)$ is a non-constant step function, since in this case $g_\zeta^\omega(0)$ is almost periodic and fails to decay as $|\omega|\rightarrow\infty$.  For step functions the trace formula has a singular counterpart derived in \cite{Gi:Szego2022}, as follows.
\begin{thm}[from \cite{Gi:Szego2022}]\label{thm-Szego}
Let $\zeta_2\in\step_+(X)$ have jump points $y_1<\cdots<y_n$ in $X$, with associated reflectivities
\[
r_j=\frac{\zeta(y_j-)-\zeta(y_j+)}{\zeta(y_j-)+\zeta(y_j+)}\qquad(1\leq j\leq n).
\]
Then
\begin{equation}\label{Szego-0}
\lim_{L\rightarrow\infty}\frac{1}{2L}\int_{-L}^L-\log\bigl(1-|g_{\zeta_2}^\omega(0)|^2\bigr)\,d\omega=\sum_{j=1}^n-\log\bigl(1-r_j^2\bigr).
\end{equation}
\end{thm}
More generally, the classical trace formula fails to hold for any $\zeta\in\pwacplus(X)\setminus C_+(X)$.
In terms of the factorization $\zeta=\zeta_1\zeta_2$, where $\zeta_1$ is continuous and $\zeta_2\in\step_+(X)$, Theorems~\ref{thm-forward-scattering}\ref{forward-scattering-asymptotics} and \ref{thm-Szego} imply 
\[
\begin{split}
\int_{-L}^L-\log\bigl(1-|g_{\zeta}^\omega(0)|^2\bigr)\,d\omega\quad&\sim\quad\int_{-L}^L-\log\bigl(1-|g_{\zeta_2}^\omega(0)|^2\bigr)\,d\omega\\
&\sim\quad 2L\sum_{j=1}^n-\log\bigl(1-r_j^2\bigr)\quad\mbox{ as }\quad L\rightarrow\infty,
\end{split}
\]
whereby
\begin{equation}\label{oh-L}
\int_{-L}^L-\log\bigl(1-|g_{\zeta}^\omega(0)|^2\bigr)\,d\omega\rightarrow\infty\quad\mbox{ as }\quad L\rightarrow\infty. 
\end{equation}
So the classical trace formula blows up if $\zeta\in\pwacplus(X)$ is discontinuous. 
On the other hand, the asymptotic formula (\ref{forward-scattering-asymptotics-2}) implies that the above singular trace formula extends without change to this more general setting, yielding the following alternative to the classical formula. 
\begin{thm}\label{thm-renormalized-trace}
Given $X=(x_0,x_1)$ and $\zeta\in\pwacplus(X)\setminus C_+(X)$, write $\zeta=\zeta_1\zeta_2$, where $\zeta_1$ is continuous and $\zeta_2\in\step_+(X)$, and denote by $r_j$ $(1\leq j\leq n)$ the reflectivities corresponding to the $n\geq 1$ jump points of $\zeta_2$.  Then 
\[
\lim_{L\rightarrow\infty}\frac{1}{2L}\int_{-L}^L-\log\bigl(1-|g_{\zeta}^\omega(0)|^2\bigr)\,d\omega=\sum_{j=1}^n-\log\bigl(1-r_j^2\bigr).
\]
\end{thm}
Thus the positive quantity $\sigma_\zeta=\sum_{j=1}^n-\log(1-r_j^2)$, which measures the strength of the reflectivities due to the singular factor $\zeta_2$, obtains directly from the reflection coefficient.  On the other hand, if $\zeta\in C_+(X)$ then 
\[
\lim_{L\rightarrow\infty}\frac{1}{2L}\int_{-L}^L-\log\bigl(1-|g_{\zeta}^\omega(0)|^2\bigr)\,d\omega=0,
\]
since in this case $g_\zeta^\omega(0)\rightarrow 0$ as $|\omega|\rightarrow\infty$.

\section{Effective computation\label{sec-fast-algorithms}}

For any given $X=(x_0,x_1)$, if $\zeta\in C_+(X)$,  the reflection coefficient is represented by the hyperbolic tangent operator as per Theorem~\ref{thm-hyperbolic}\ref{R-Th}.  Letting $\zeta_n$ denote the standard approximation to $\zeta$, the reflection coefficients $R_n=g_{\zeta_n}(0)$ may be expressed in terms of OPUC as 
\begin{equation}\label{computational-approximation}
g_{\zeta_n}^\omega(0)=\frac{\Psi_n^\ast(e^{2i\Delta\omega})-\Phi_n^\ast(e^{2i\Delta\omega})}{\Psi_n^\ast(e^{2i\Delta\omega})+\Phi_n^\ast(e^{2i\Delta\omega})}\qquad(\omega\in\real),
\end{equation}
by Proposition~\ref{prop-opuc}\ref{equally-spaced-R-representation}.
This OPUC representation underpins simple forward and inverse algorithms that are fast and accurate, exploiting the recurrence relations (\ref{OPUC-recursion}).

\subsection{Forward scattering.\label{sec-forward-computation}}

Theorems~\ref{thm-continuous} and \ref{thm-singular-representation} yield an efficient method to compute the reflection coefficient $R(\omega)=g_\zeta^\omega(0)$ at a given frequency $\omega$ for any $\zeta\in\reg_+(X)$.  Fix an $n$-jump step function $\zeta_n\in\step_+(X)$ sufficiently close to $\zeta$ in $L^\infty(X)$; define $2\times 2$ matrices $M_j$ $(1\leq j\leq n+1)$ in terms of $\zeta_n$ according to Proposition~\ref{prop-matrix-product}; evaluate the matrix product (\ref{matrix-product-E}), and then the expression for $R(\omega)$ given in Theorem~\ref{thm-singular-representation}\ref{singular-R-representation}.  The computation requires $O(n)$ flops (per point $\omega$) and is guaranteed to be accurate if $\|\zeta_n-\zeta\|_{\infty}$ is sufficiently small.

However, in acoustic imaging the measured quantity determined by forward scattering is $\rwc{R}$, not $R$.  Computation of $\rwc{R}$ in the case $\zeta\in C_+^1(X)$ (or more generally, $\zeta\in\pwplus(X)\cap C(X)$) is facilitated by the representation (\ref{OPUC-notation}) of OPUC as a matrix product in combination with uniform convergence on compact sets of (\ref{computational-approximation}) to $R$.  Note that $\zeta$ is determined by the restriction of $\rwc{R}$ to the interval $(0,2(x_1-x_0))$. The following algorithm approximates $\rwc{R}(t)$ for a sequence of $n$ equally-spaced times $t=t_j$,  where $0<t_j<t_{\max}=2(x_1-x_0)$, taking as inputs the impedance function $\zeta$, the parameters $x_0<x_1$ and the index $n$.  (The algorithm extends in an obvious way to arbitrarily large times, in effect by concatenation of $\zeta$ with an impedance function having constant value $\zeta(x_1-)$ at $x\geq x_1$.)  The idea is to compute the initial Fourier coefficients of $R_n=g_{\zeta_n}(0)$, where $\zeta_n$ is the $n$th standard approximant to $\zeta$.  That is, $R_n$ has the form 
\[
R_n(\omega)=\sum_{j=1}^\infty a_jz^j,
\]
where $z=e^{2i\Delta_n\omega}$ and $\Delta_n=(x_1-x_0)/(n+1)$.  Therefore 
\[
\rwc{R}_n(t)=\sum_{j=1}^\infty a_j\delta(t-t_j)\qquad(t_j=2j\Delta_n).
\]
The desired values $\rwc{R}(t_j)$ for $1\leq j\leq n$ are then given by the approximation $\rwc{R}(t_j)\approx a_j/(2\Delta_n)$, discussed further below.  
\begin{alg}[Forward scattering]\label{alg-forward}\ \\ Input: $\zeta$, $x_0$, $x_1$, $n$. \begin{enumerate}
\item Set $\Delta_n=(x_1-x_0)/(n+1)$, $z=e^{2i\Delta_n}$,  and 
\[
y_j=x_0+j\Delta,\quad r_j=\frac{\zeta(y_j-\Delta_n/2)-\zeta(y_j+\Delta_n/2)}{\zeta(y_j-\Delta_n/2)+\zeta(y_j+\Delta_n/2)}\qquad(1\leq j\leq n).
\]
\item Set $\Phi_0=\Phi_0^\ast=\Psi_0=\Psi_0^\ast=1$, and for $j=0:n-1$ compute the coefficients of the polynomials $\Phi_{j+1},\Phi_{j+1}^\ast,\Psi_{j+1},\Psi_{j+1}^\ast$ by the OPUC recursion
\[
\begin{pmatrix}
\Psi_{j+1}(z)&\Psi_{j+1}^\ast(z)\\
-\Phi_{j+1}(z)&\Phi_{j+1}^\ast(z)
\end{pmatrix}
=
\begin{pmatrix}
z\Psi_j(z)+\overline{r_{j+1}}\Psi_j^\ast(z)&\Psi_j^\ast(z)+r_{j+1}z\Psi_j(z)\\
-\left(z\Phi_j(z)-\overline{r_{j+1}}\Phi_j^\ast(z)\right)&\Phi_j^\ast(z)-r_{j+1}z\Phi_j(z)
\end{pmatrix}. 
\]
\item Let $b_1,\ldots,b_n$ and $c_0,\ldots,c_n$ denote the respective coefficients of $(\Psi_n^\ast-\Phi_n^\ast)/2$ and $(\Psi_n^\ast+\Phi_n^\ast)/2$ (so in particular $c_0=1$). To compute the initial coefficients $d_0,\ldots,d_{n-1}$ of $2/(\Psi_n^\ast+\Phi_n^\ast)$, set $d_0=1$ and solve the system
\[
\begin{pmatrix}1&0&0&\cdots&0\\ c_1&1&0&\cdots&0 \\ \vdots&\ddots&\ddots&\ddots&\vdots\\ c_{n-1}&\cdots&c_1&1&0\end{pmatrix}\begin{pmatrix}1\\ d_1\\ \vdots\\ d_{n-1}\end{pmatrix}=\begin{pmatrix}1\\ 0\\ \vdots\\ 0\end{pmatrix}
\]
recursively by back substitution. 
\item To compute the initial coefficients $a_1,\ldots,a_n$ of $(\Psi_n^\ast-\Phi_n^\ast)/(\Psi_n^\ast+\Phi_n^\ast)$, evaluate the discrete convolutional product
\[
\begin{pmatrix}a_1\\ a_2\\ \vdots\\ a_n\end{pmatrix}=
\begin{pmatrix}b_1&0&\cdots&0\\ b_2&b_1&\ddots&\vdots \\ \vdots&\ddots&\ddots&0\\ b_n&\cdots&b_2&b_1\end{pmatrix}\begin{pmatrix}1\\ d_1\\ \vdots\\ d_{n-1}\end{pmatrix}.
\]
\item Set $\widetilde{a}_j=a_j/(2\Delta_n)$ and $t_j=2j\Delta_n$ $(1\leq j\leq n)$. 
\end{enumerate}
Output:  
$\left\{(t_j,\widetilde{a}_j)\,|\,1\leq j\leq n\right\}$, where $\widetilde{a}_j\approx\rwc{R}(t_j)$. 
\end{alg}
Algorithm~\ref{alg-forward} is easily adapted to accommodate $\alpha$ as input in place of $\zeta$ by replacing the definition of $r_j$ in step~(1) with 
\begin{equation}\label{algorithm-one-alternate}
r_j=\tanh\bigl(\Delta_n\alpha(y_j)\bigr).
\end{equation}
The algorithm requires $O(n^2)$ flops, and is inherently stable (see \cite{Gi:NMPDE2018}).  It turns out also to be remarkably accurate in practice.  However, from the theoretical point of view, the fact that $R_n\rightarrow R$ uniformly on compact sets implies only that $\rwc{R}_n\rightarrow\rwc{R}$ weakly.  Thus for any fixed $t\in2(x_1-x_0)$, for example, any small $\varepsilon>0$, and any $C^\infty$ approximation $\varphi$ of $\chi_{(t-\varepsilon,t+\varepsilon)}$,
\[ 
\int_{t-\varepsilon}^{t+\varepsilon}\rwc{R}_n\approx\int\rwc{R}_n\varphi\rightarrow\int\rwc{R}\varphi\approx\int_{t-\varepsilon}^{t+\varepsilon}\rwc{R}.
\]
Weak approximation alone does not guarantee that $\int\rwc{R}_n\varphi$ is close to $\int\rwc{R}\varphi$ when $\Delta_n\approx\varepsilon$.  Rather, $n$ may have to be larger, such that $k\Delta_n\approx \varepsilon$ for some sufficiently large integer $k>1$.  Suppose that $t=t_j$ is a point of continuity of $\rwc{R}$ and $\varepsilon$ is sufficiently small that $\rwc{R}(t)\approx\rwc{R}(t_j)$ for $t\in(t_j-\varepsilon,t_j+\varepsilon)$.  Then
\[
\int_{t_j-\varepsilon}^{t_j+\varepsilon}\rwc{R}_n=\sum_{\nu=-k+1}^{k-1}a_{j+\nu}\quad\mbox{ and }\quad \int_{t_j-\varepsilon}^{t_j+\varepsilon}\rwc{R}\approx 2k\Delta_n\rwc{R}(t_j),
\]
so that for sufficiently large $n$ one has the approximation 
\begin{equation}\label{numerical-approximation-1}
\rwc{R}(t_j)\approx\frac{\sum_{\nu=-k+1}^{k-1}a_{j+\nu}}{2k\Delta_n}.  
\end{equation}
Algorithm~\ref{alg-forward} uses the value $k=1$. That this actually works in numerical experiments, seemingly without exception, calls for an explanation. Indeed in computing the full wave field $k>1$ is necessary. That is, based on Lemma~\ref{lem-forward-stability} one can also use step function approximations $\zeta_n$ of $\zeta$ to compute not just $\rwc{R}$, but the full wave field on $X$ in the time domain. For the full wave field, a value of $k\geq2$ turns out to be necessary (and $k=2$ suffices); otherwise the approximation fails.  The latter fact can be directly inferred from an analysis of the generalized algorithm, which is outside the scope of the present paper; see \cite{Ch:2021}.  

One can compute numerical reflection data also in the setting $\zeta\in\pwplus(X)$, without requiring continuity.  In this case $\rwc{R}$ ceases to be a regular function.  It acquires a purely singular distributional part, the numerical aspects of which are technically more complicated, although not intractable; see \cite{Ch:2021}.  In a physical setting, where a purely impulsive source signal is not possible, the scattering data is again a regular function, and one can modify Algorithm~\ref{alg-forward} to handle essentially arbitrary $\zeta$.

\subsection{Inverse scattering.\label{sec-inverse-computation}}

The same approximation underlying the last step of Algorithm~\ref{alg-forward}, namely 
\[
\rwc{R}(t_j)\approx a_j/(2\Delta_n),\quad\mbox{ where }\quad\rwc{R}_n(t)=\sum_{j=1}^\infty a_j\delta(t-t_j),
\]
serves as a starting point for an effective inverse algorithm, which also involves the moments corresponding to the orthogonal polynomials occurring in (\ref{computational-approximation}).  The defining property of monic polynomials $\Phi_0,\Phi_1,\ldots,\Phi_n$ is that they are orthogonal with respect to a probability measure $d\mu$ on the unit circle (cf. \S\ref{sec-OPUC}), and so can be constructed from the moments of $d\mu$,  
\[
m_j=\int_{S^1}z^j\,d\mu(z)\qquad(j\geq 0),
\]
via Gram-Schmidt orthogonalization of the monomial sequence $z^j$.   
The first $n$ coefficients of $\rwc{R}_n=\sum_{j=1}^\infty a_j\delta(t-t_j)$ determine the $n+1$ moments $m_0,\ldots,m_n$ through the Herglotz formula 
\begin{equation}\label{Herglotz-1}
\frac{1+R_n(\omega)}{1-R_n(\omega)}=\int_{S^1}\frac{\xi+e^{2i\Delta_n\omega}}{\xi-e^{2i\Delta_n\omega}}\,d\mu(\xi)
\end{equation}
as follows.  Setting $z=e^{2i\Delta_n\omega}$, equation (\ref{Herglotz-1}) has the form
\begin{equation}\label{Herglotz-2}
\frac{1+a_1z+a_2z^2+\cdots}{1-a_1z-a_2z^2-\cdots}=1+2\overline{m}_1z+2\overline{m}_2z^2+\cdots
\end{equation}
whereby one can equate coefficients of powers of $z$ to pass from the $a_j$ to the $m_j$.  First, note by (\ref{Herglotz-2}) that the moments $m_j$ are real, since the coefficients $a_j$ are real, and $m_0=1$ since $\mu$ is a probability measure.  Define $A$ to be the $(n+1)\times(n+1)$ lower triangular Toeplitz matrix given by the formula 
\begin{equation}\label{matrix-A}
A(i,j)=a_{\max\{i-j,0\}}\qquad(1\leq i,j\leq n+1).
\end{equation}
Then the vector of Taylor coefficients of $(1+R)/(1-R)$ (i.e., the left-hand side of (\ref{Herglotz-2})) up to degree $n$ is 
\[
\left(I+2\sum_{j=1}^nA^j\right)\begin{pmatrix}1\\ 0\\ \vdots\\ 0\end{pmatrix},
\]
and so by (\ref{Herglotz-2}), 
\[
\left(\sum_{j=0}^nA^j\right)\begin{pmatrix}1\\ 0\\ \vdots\\ 0\end{pmatrix}=\begin{pmatrix}1\\ m_1\\ \vdots\\ m_n\end{pmatrix}. 
\]
But $A$ is nilpotent, with $A^{n+1}=0$, so $\sum_{j=0}^nA^j=\sum_{j=0}^\infty A^j=(1-A)^{-1}$, and it follows that 
\begin{equation}\label{A-m-equation}
(I-A)\begin{pmatrix}1\\ m_1\\ \vdots\\ m_n\end{pmatrix}=\begin{pmatrix}1\\ 0\\ \vdots\\ 0\end{pmatrix}, 
\end{equation}
which may be solved efficiently for the $m_j$ in terms of the $a_j$ by back substitution.  

The moments determine the orthogonal polynomials; the latter determine the recurrence coefficients (i.e., reflectivities) by way of the formula
\begin{equation}\label{rj-via-OPs}
r_{j+1}=\frac{\langle z\Phi_j,1\rangle_{d\mu}}{\langle\Phi_j^\ast,1\rangle_{d\mu}},
\end{equation}
a consequence of the recursion
\[
\Phi_{j+1}=z\Phi_j-r_{j+1}\Phi_j^\ast
\]
(since $\langle\Phi_{j+1},1\rangle_{d\mu}=0$). The reflectivities in turn determine the values of standard approximant $\zeta_n$ to $\zeta$ via (\ref{rnj}), and thus the values of $\zeta$ on the regular spatial grid $y_j=x_0+j\Delta_n$.   The foregoing observations lead to the following. 
\begin{alg}[Inverse scattering]\label{alg-inverse}\ \\
Input: $x_0, \zeta(x_0), \{(t_j,\rwc{R}(t_j))\,|\,1\leq j\leq n\}$, where the $t_j$ are equally spaced of the form $t_j=j\tau$, for some fixed $\tau>0$. 
\begin{enumerate}
\item Set $a_j=\tau\rwc{R}(t_j)$ $(1\leq j\leq n)$, and $A(i,j)=a_{\max\{i-j,0\}}$ $(1\leq i,j.\leq n+1)$. 
\item Set $m_0=1$ and compute moments $m_1,\ldots,m_n$ by solving the lower Toeplitz equation 
\[
(I-A)\begin{pmatrix}1\\ m_1\\ \vdots\\ m_n\end{pmatrix}=\begin{pmatrix}1\\ 0\\ \vdots\\ 0\end{pmatrix}
\]
by back substitution. 
\item Set $\Phi_0=\Phi_0^\ast=1$ and $r_1=m_1$. For $j=1:n-1$, compute the coefficients of polynomials
$
\Phi_j=z\Phi_{j-1}-r_j\Phi_{j-1}^\ast$ and $\Phi_j^\ast=\Phi_{j-1}^\ast-r_jz\Phi_{j-1},
$
and set 
\[
r_{j+1}=\langle z\Phi_j,1\rangle_{d\mu}/\langle\Phi_j^\ast,1\rangle_{d\mu}
\]
(the latter being computed in terms of the coefficients of $z\Phi_j$ and $\Phi_j^\ast$, and the moments $m_0,\ldots,m_{j+1}$). 
\item Set $\widetilde\zeta_j=\zeta(x_0)\exp\bigl(-2\sum_{\nu=1}^j\tanh^{-1}r_\nu\bigr)$ $(1\leq j\leq n)$. 
\end{enumerate}
Output: $\left\{\left.\bigl(y_j,\widetilde\zeta_j\bigr)\,\right|\,1\leq j\leq n\right\}$ where $y_j=x_0+j\tau/2$.  Here $\widetilde\zeta_j\approx\zeta(y_j)$. 
\end{alg}
To compute $\alpha$ rather than $\zeta$, replace step~(4) of Algorithm~\ref{alg-inverse} with 
\begin{itshape}
\begin{enumerate}[label={(\arabic*$^{\,\prime\!}$)},itemindent=0em]
\setcounter{enumi}{3}
\item\  Set $\widetilde\alpha_j=\frac{1}{\Delta_n}\tanh^{-1}r_\nu$ $(1\leq j\leq n)$. 
\end{enumerate}
\end{itshape}
The output then becomes $\left\{\left.\bigl(y_j,\widetilde\alpha_j\bigr)\,\right|\,1\leq j\leq n\right\}$ where $\widetilde\alpha_j\approx\alpha(y_j)$.

A  version of Algorithm~\ref{alg-inverse} adapted to the physically more realistic situation of a non-Dirac acoustic source is presented in \cite{Gi:JCP2018}, along with numerical experiments that show remarkable accuracy and stability. (For a non-Dirac source it is no longer necessary to require $\zeta$ be continuous; one can invert data for essentially arbitrary impedance.)  Algorithm~\ref{alg-inverse} has an operation count of $O(n^2)$ flops and is in practice substantially faster to evaluate than its forward counterpart Algorithm~\ref{alg-forward}---contrary 
to the intuition that inverse scattering is more complicated than forward scattering!

\subsection{An example.\label{sec-computation-example}}
\begin{figure}
\fbox{
\includegraphics[width=.3\textwidth]{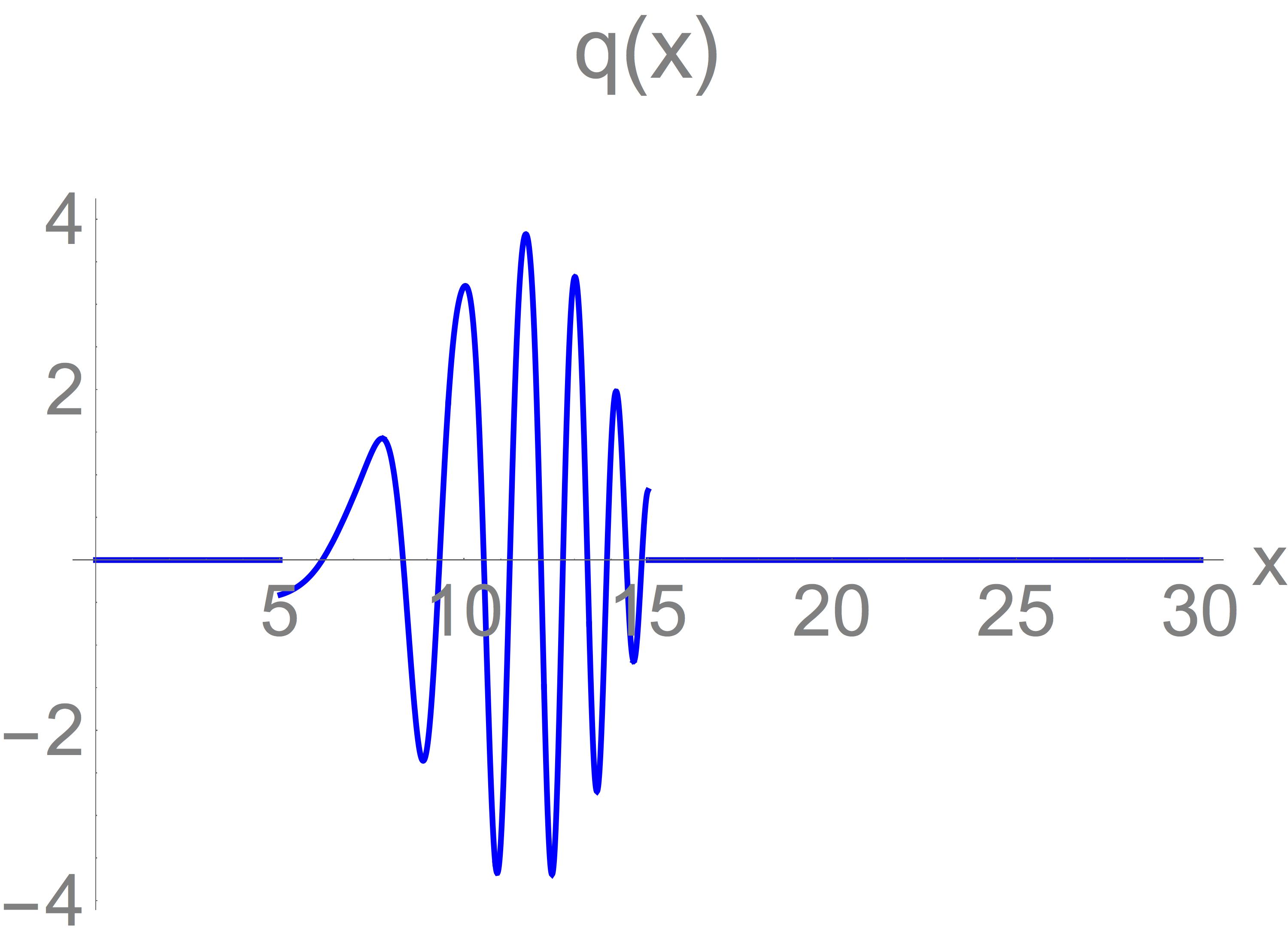}\rule{16pt}{0pt}\includegraphics[width=.3\textwidth]{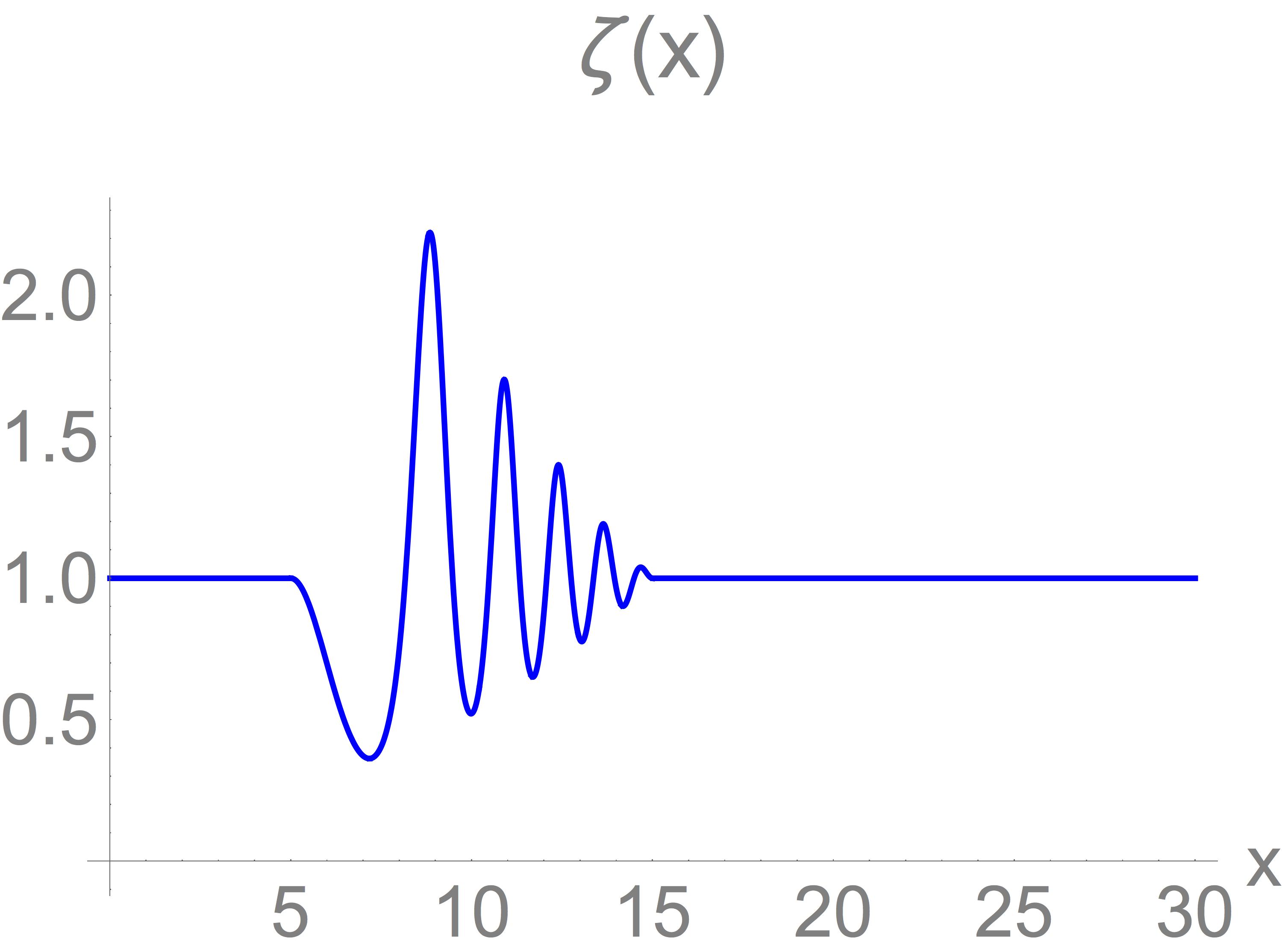}}
\caption{
Plots of the potential $q$ (left) and corresponding impedance $\zeta$ (right).\label{fig-q-zeta}
}
\end{figure}

Fixing parameters $(a,b,c,d)=(5,15,.065,\pi/10)$, consider the potential depicted in Figure~\ref{fig-q-zeta},
\[
q(x)=\left\{
\begin{array}{cc}
0&\mbox{ if }x<a\\
\begin{split}
&\left(2cd(b-x)(x-a)\cos\bigl(d(x-a)^2\bigr)-c\sin\bigl(d(x-a)^2\bigr)\right)^2\\
&+4cd^2(b-x)(x-a)^2\sin\bigl(d(x-a)^2\bigr)\\
&-2cd(b+2a-3x)\cos\bigl(d(x-a)^2\bigr)
\end{split}&\mbox{ if }a\leq x<b\\
0&\mbox{ if }b\leq x
\end{array}\right.,
\]
corresponding to the impedance function
\[
\zeta(x)=\left\{
\begin{array}{cc}
1&\mbox{ if }x<a\\
\exp\left(-2c(b-x)\sin\bigl(d(x-a)^2\bigr)\right)
&\mbox{ if }a\leq x<b\\
1&\mbox{ if }b\leq x
\end{array}\right.
\]
by way of the formula $q=\sqrt{\zeta}^{\,\prime\prime}/\sqrt{\zeta}$.  
The reflection coefficient $R(\omega)$ for the Schr\"{o}dinger equation
\[
-y^{\prime\prime}+qy=\omega^2y
\]
coincides precisely with the reflection coefficient $g_\zeta^\omega(0)$ for $\zeta\in C^1(X)$, where $X=(0,15)$. Thus the acoustic imaging problem and inverse scattering on the line for the Schr\"{o}dinger equation (extending $q$ to all of $\real$ by assigning it the value $0$ outside of $(0,15)$) are equivalent. The acoustic imaging problem for $\zeta$ is to measure the sound field $$d(t)=u(x_0,t)\qquad (0<t<t_{\max})$$ that results from a unit Dirac impulse transmitted rightward from $x_0$ at time $t=0$ toward the support of $\zeta^\prime$, and to compute $\zeta$ from these values of $d(t)$.  The inverse scattering problem for the Schr\"{o}dinger equation is to compute $q$ given $R$.  The methods described in \S\ref{sec-forward-computation} and \S\ref{sec-inverse-computation} allow efficient computation of the scattering data $R$ and $d=\rwc{R}$ for the Schr\"{o}dinger and acoustic equations, as well as reconstruction of $\zeta$ and $q$ from $\rwc{R}$.  More precisely, the restriction of $\zeta$, and hence $q$, to $(0,x_1)$ is determined by $d|_{(0,2x_1)}$. Direct computation of $\zeta$ and $q$ from a sampled version of $d|_{(0,2x_1)}$ is prescribed by Algorithm~\ref{alg-inverse}.  
Results of these computations for $x_1=30$ are presented here as a simple illustration (see \cite{Gi:JCP2018,Ch:2021} for more comprehensive discussion). 

\begin{figure}
\fbox{
\includegraphics[width=.3\textwidth]{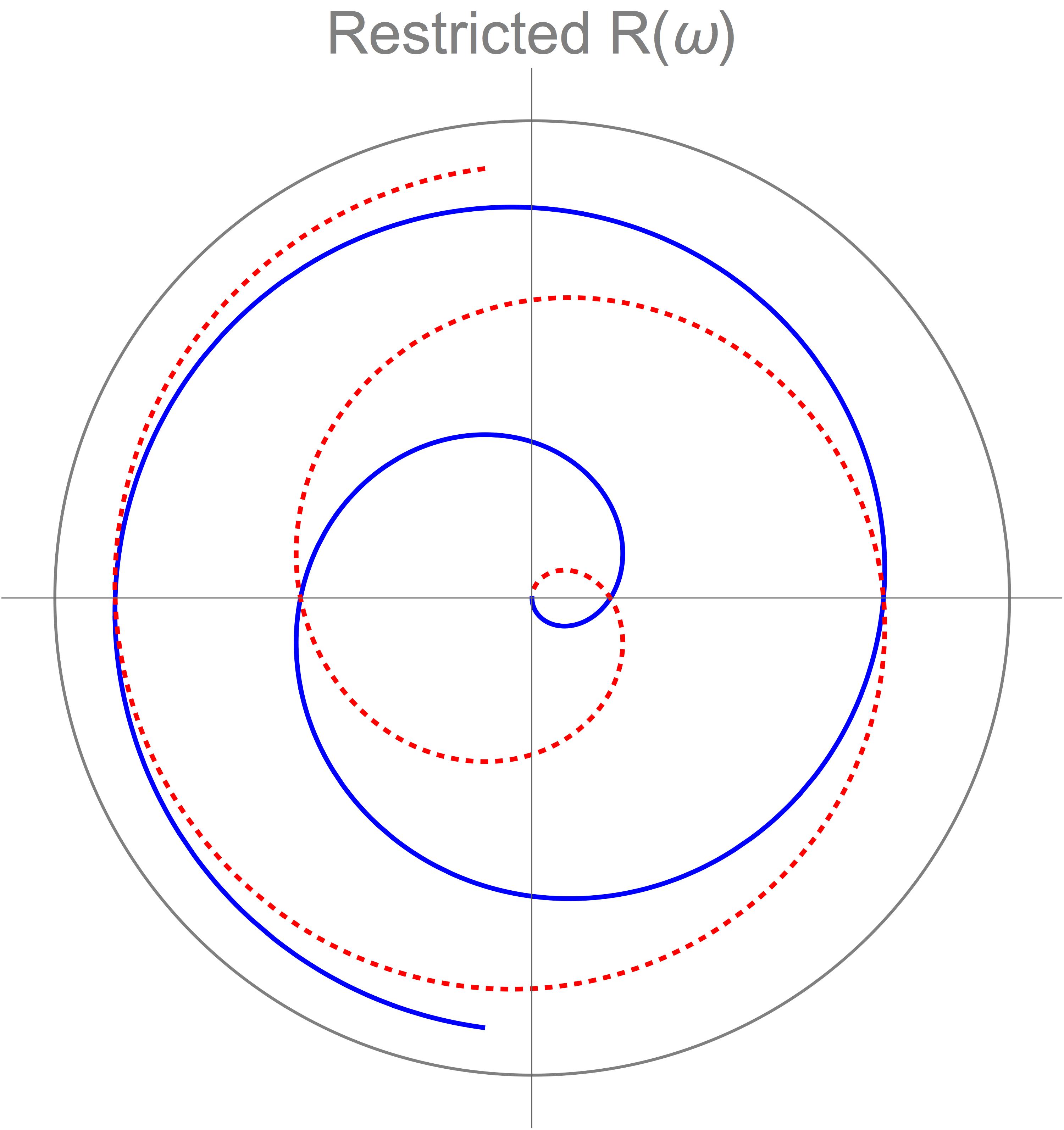}\includegraphics[width=.3\textwidth]{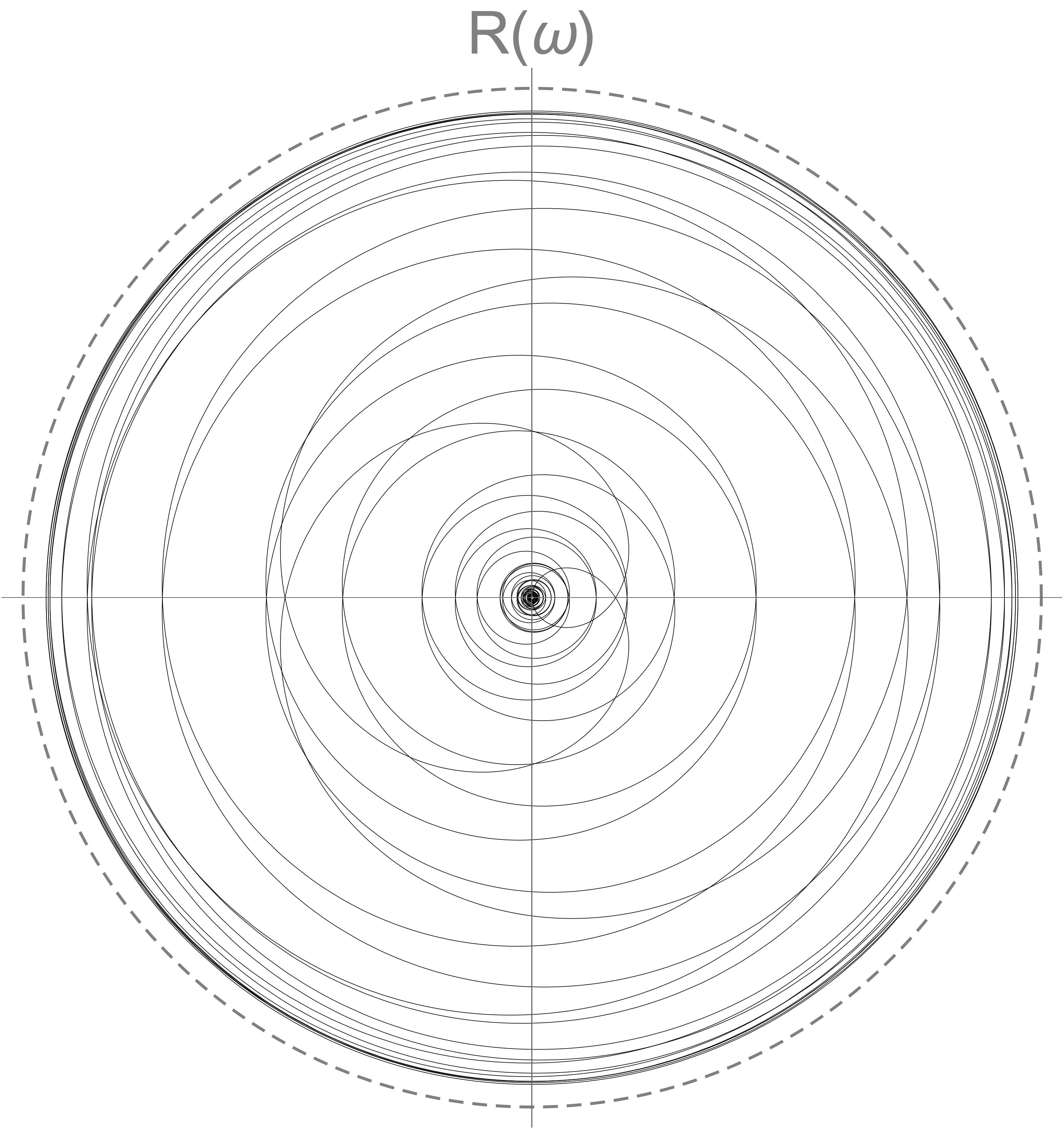}
\parbox[b]{.3\textwidth}{\includegraphics[width=.3\textwidth]{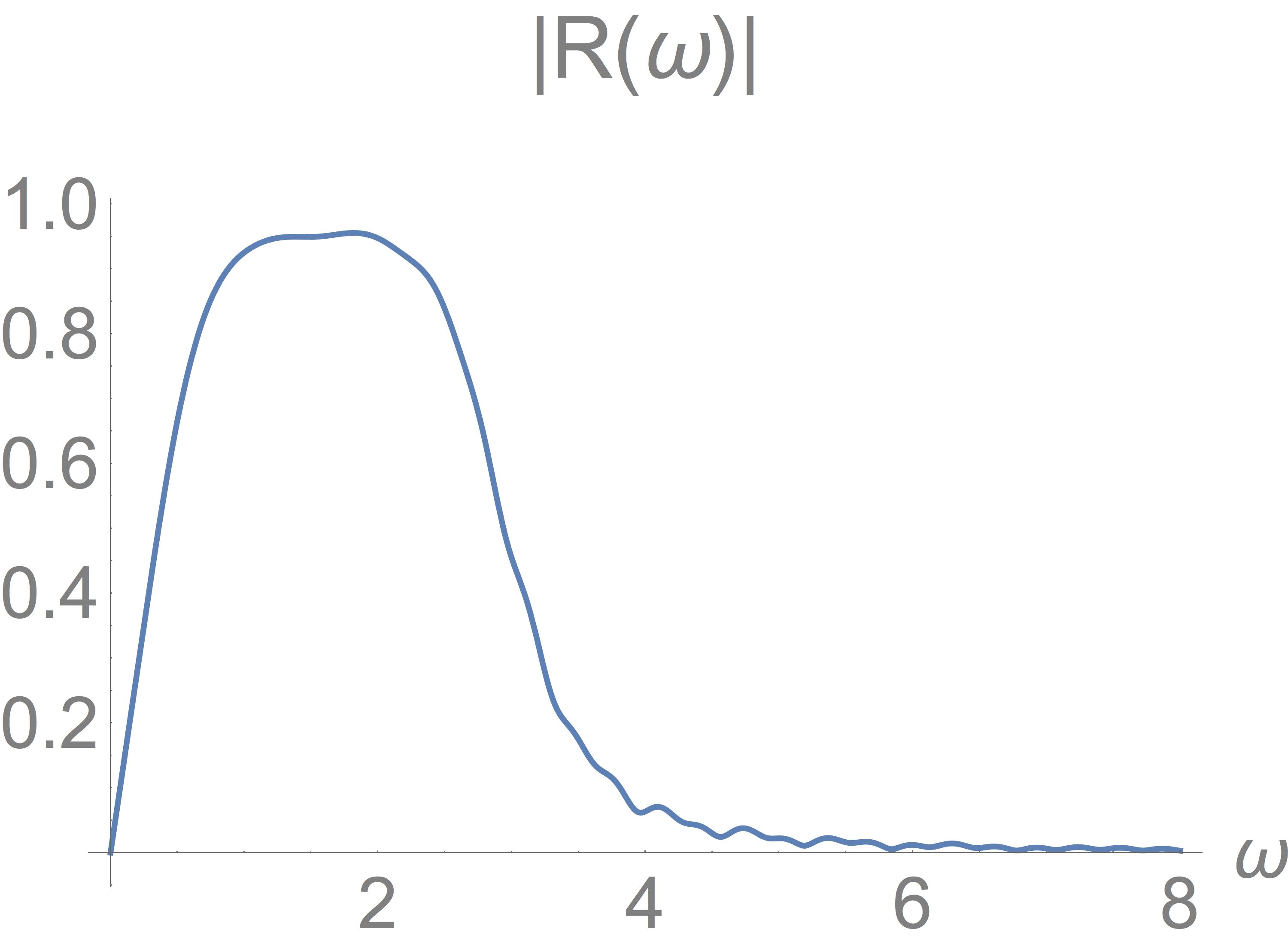}\\[13pt]}}
\caption{Reflection coefficient $R(\omega)\in\complex$, computed using OPUC of degree 4000. Left: $R|_{(-.9,.9)}$; dotted red curve corresponds to $-.9<\omega<0$, solid blue curve to $0<\omega<.9$. Middle: $R|_{(-8,8)}$ plotted at 1000 equally-spaced points; dashed curve is the unit circle; note $R(\omega)\approx0$ if $|\omega|\geq8$. Right: $|R|$ for $\omega>0$; note that $|R|$ is even. Total computing time $8.03$ sec.\label{fig-R}}
\end{figure}

Figure~\ref{fig-R} depicts the reflection coefficient $R$, computed using matrix products as described at the beginning of \S\ref{sec-forward-computation}.  Accuracy of the plots is such that no discernible change results if one doubles the degree of OPUC used in the computation to 8000 from 4000 (the relative $\ell^2$ difference between the two computations is $0.004\%$). 

\begin{figure}
\fbox{
\includegraphics[width=.3\textwidth]{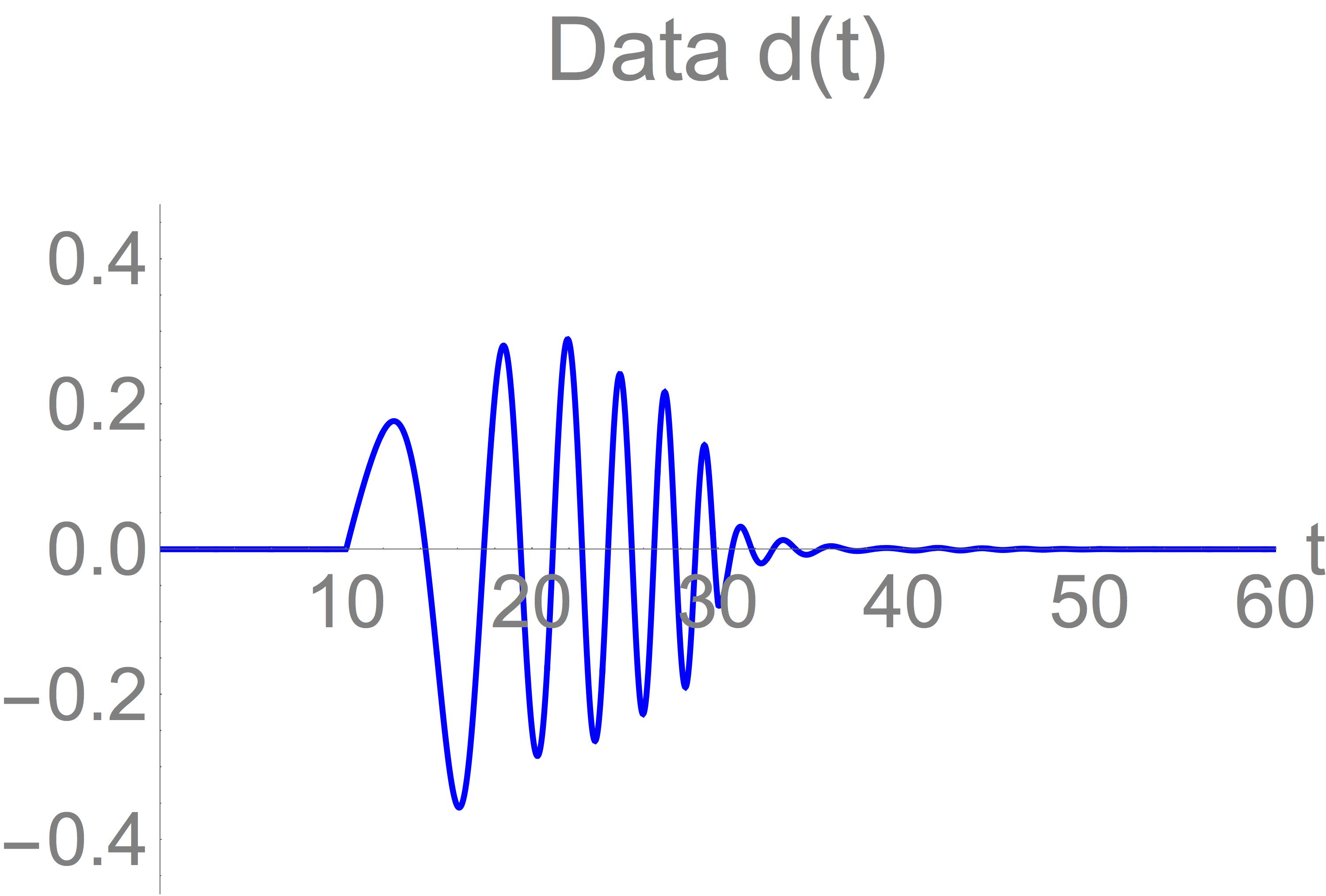}\includegraphics[width=.3\textwidth]{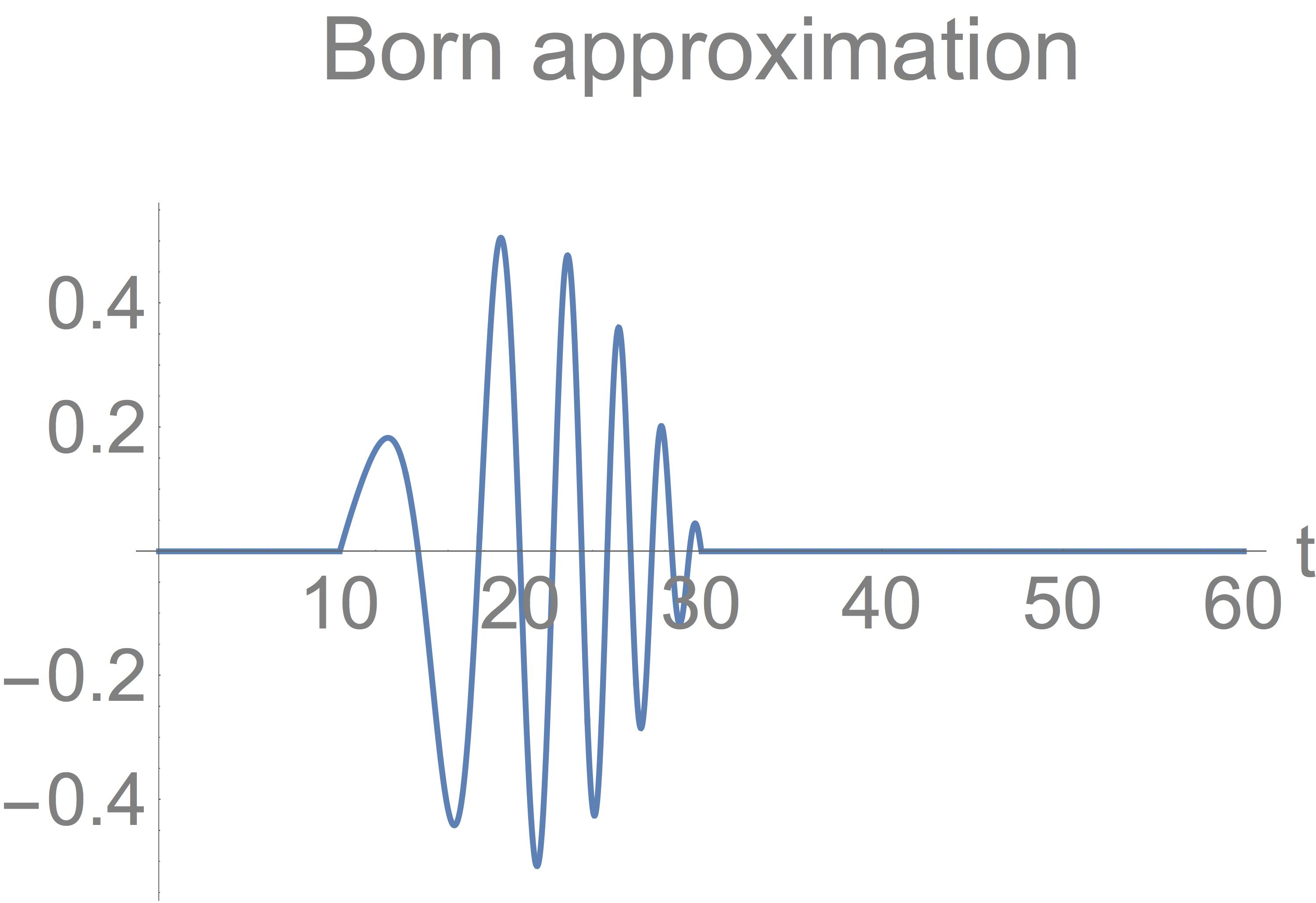}
\includegraphics[width=.3\textwidth]{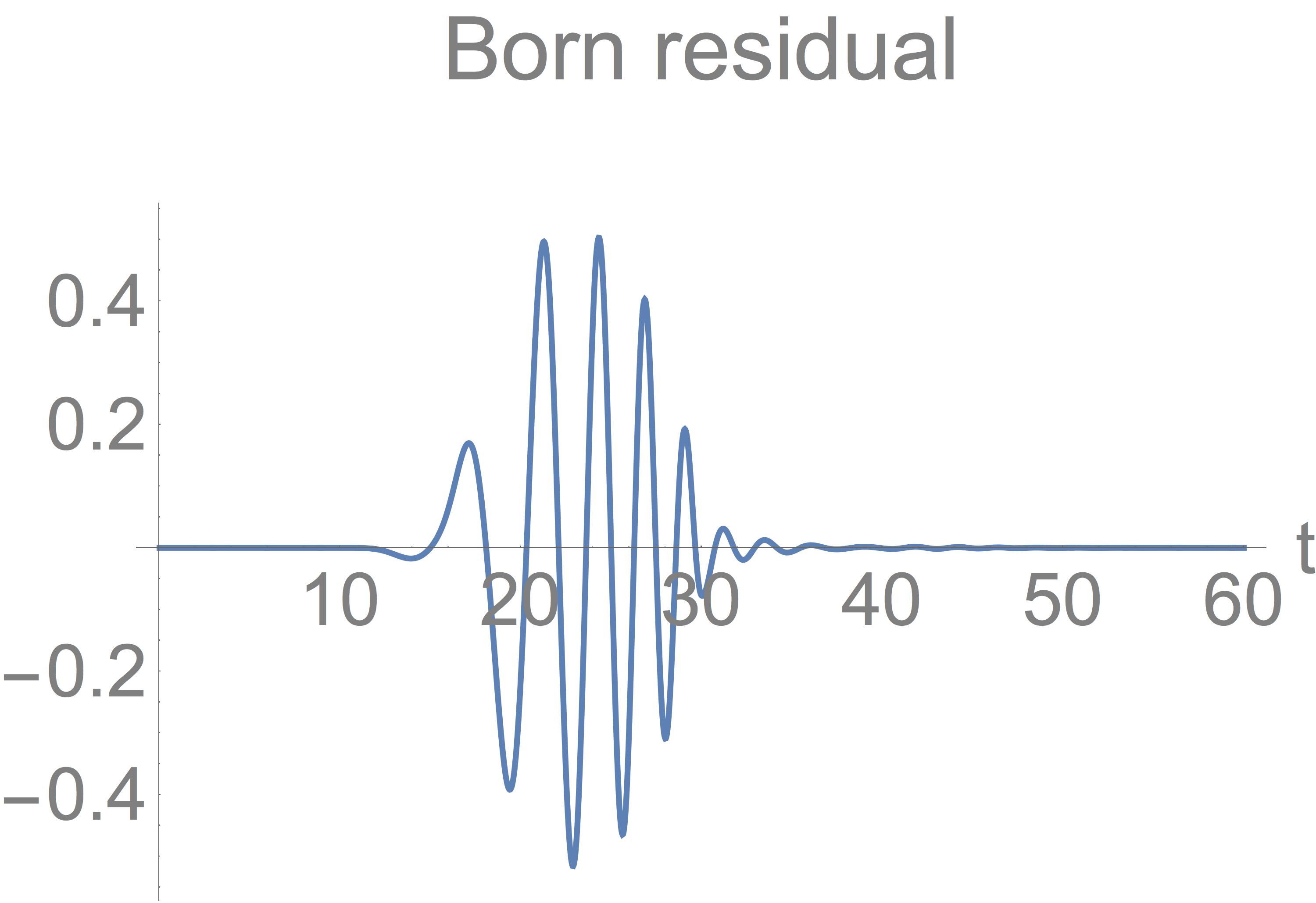}}
\caption{
Measured data $d$, the inverse Fourier transform of $R$. Left: $d$ computed using Algorithm~\ref{alg-forward} with input: $\zeta$, $x_0=0$, $x_1=30$, $n=2000$; computing time $7.97$ sec. Middle: the Born approximation $d(t)\approx\widetilde{\alpha}(t)=\textstyle\frac{1}{2}\alpha\bigl(\textstyle\frac{1}{2}t\bigr)$. Right: the Born residual $d-\widetilde{\alpha}$. The Born residual has larger $\ell^2$ norm than the true data, the relative error in the Born approximation being $130.6\%$.\label{fig-d}
}
\end{figure}

Figure~\ref{fig-d} depicts the measured data $d=\rwc{R}$ as computed using Algorithm~\ref{alg-forward} with input parameter $n=2000$, alongside the classical Born approximation and Born residual.   No perceptible change in the data results if one increases $n$ beyond $2000$: for instance, applying Algorithm~\ref{alg-forward} with $n=32000$ (which takes 38.76 min.~computing time) and downsampling by a factor of 16, the relative $\ell^2$ difference in the two vectors is $0.052\%$ (the absolute $\ell^\infty$ difference is $0.000825$).  

\begin{figure}
\fbox{
\includegraphics[width=.3\textwidth]{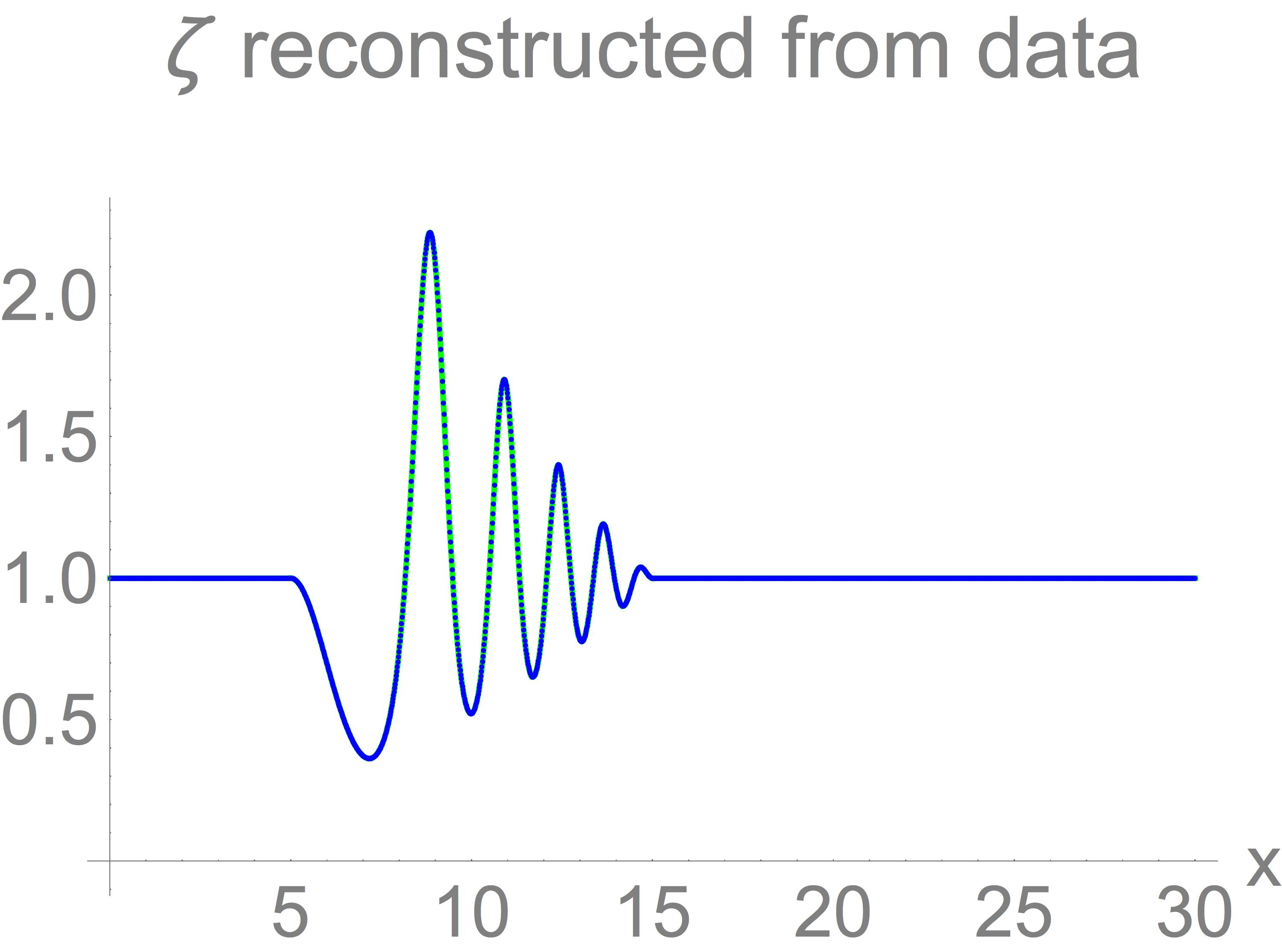}\rule{12pt}{0pt}\includegraphics[width=.3\textwidth]{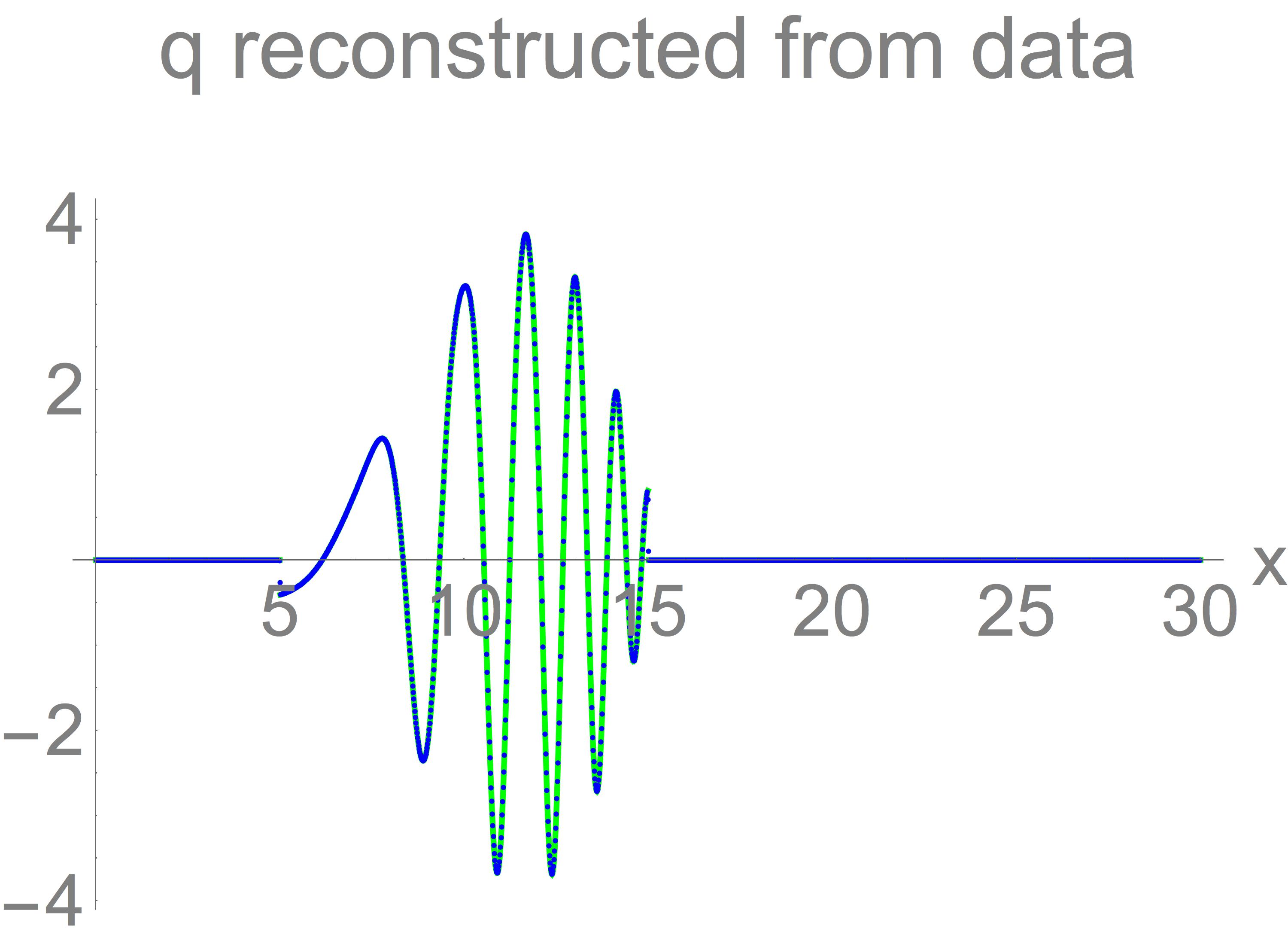}\rule{12pt}{0pt}\includegraphics[width=.3\textwidth]{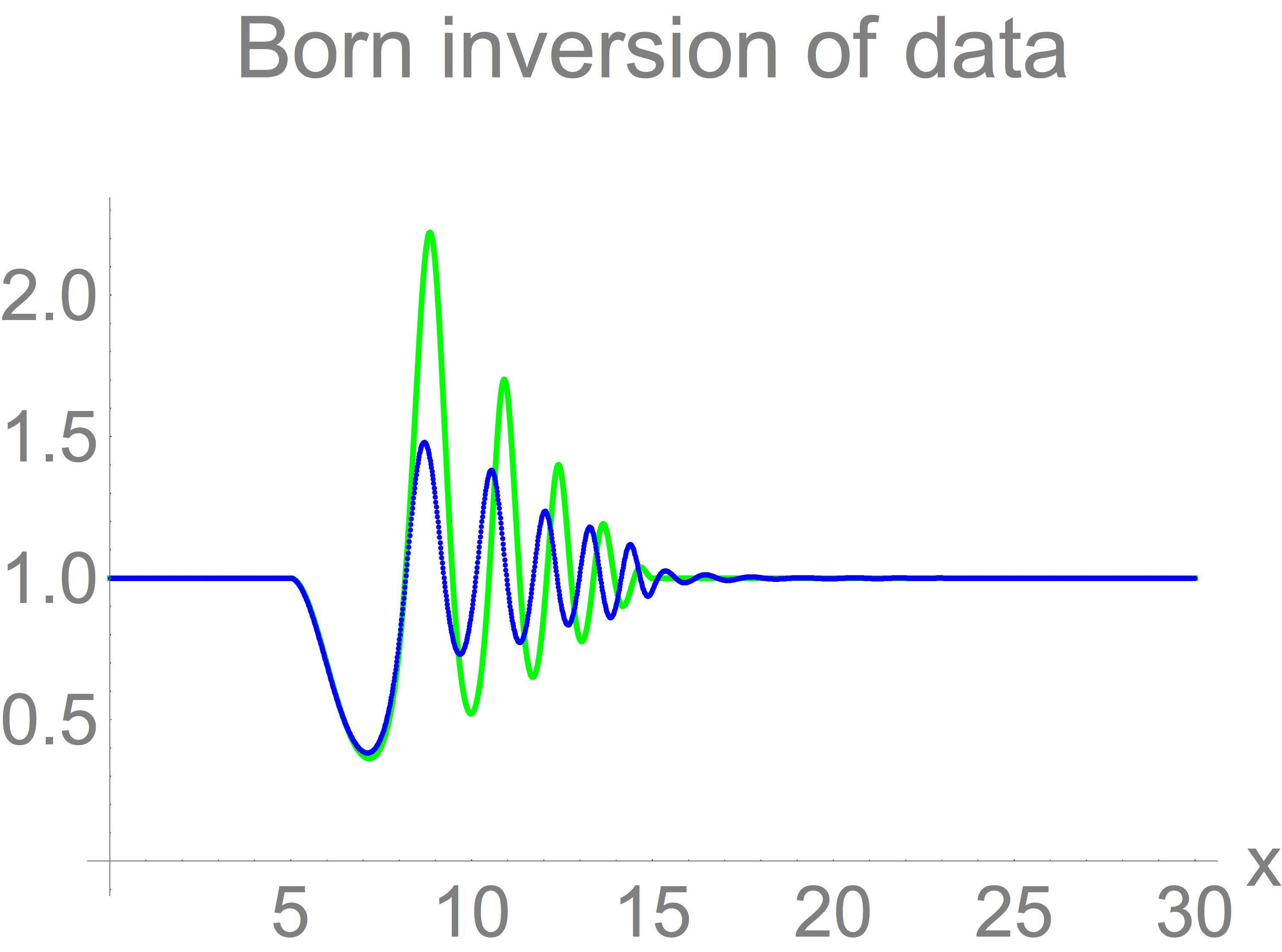}}
\caption{Reconstructed impedance and potential.  Left: in blue, $\zeta$ computed on a grid of 2000 time points by Algorithm~\ref{alg-inverse} from the data $d$ as depicted at left in Figure~\ref{fig-d}; computing time $0.786$ sec; true impedance in green. Middle: in blue, $q$ on a grid of 2000 time points computed from the inversion at left; true potential in green. Right: in blue, approximate $\zeta$ obtained by Born inversion of the measured data; true impedance in green. 
\label{fig-reconstruction}}
\end{figure}

Algorithm~\ref{alg-inverse} inverts the data from Figure~\ref{fig-d} to recover $\zeta$ on a regular grid to essentially perfect accuracy, within a relative $\ell^2$ error of $1.46\times10^{-13}$ percent in a computing time of $0.786$ seconds, $10$ times faster than the forward computation.  A second divided difference applied to $\sqrt{\zeta}$ on the given grid recovers $q$ to less accuracy, with a relative $\ell^2$ error of $0.429\%$. See Figure~\ref{fig-reconstruction}.  A standard inversion method in seismic imaging is to treat the measured data as if the Born approximation were correct, and compute the corresponding $\zeta$; the result of this Born inversion applied to the measured data is displayed at the right of Figure~\ref{fig-reconstruction} for comparison.  The $\ell^2$ relative error of the Born inversion is $17.1\%$---more accurate than the Born approximation of the data, but 14 orders of magnitude larger than error of Algorithm~\ref{alg-inverse}. 

\begin{figure}
\fbox{
\includegraphics[width=.4\textwidth]{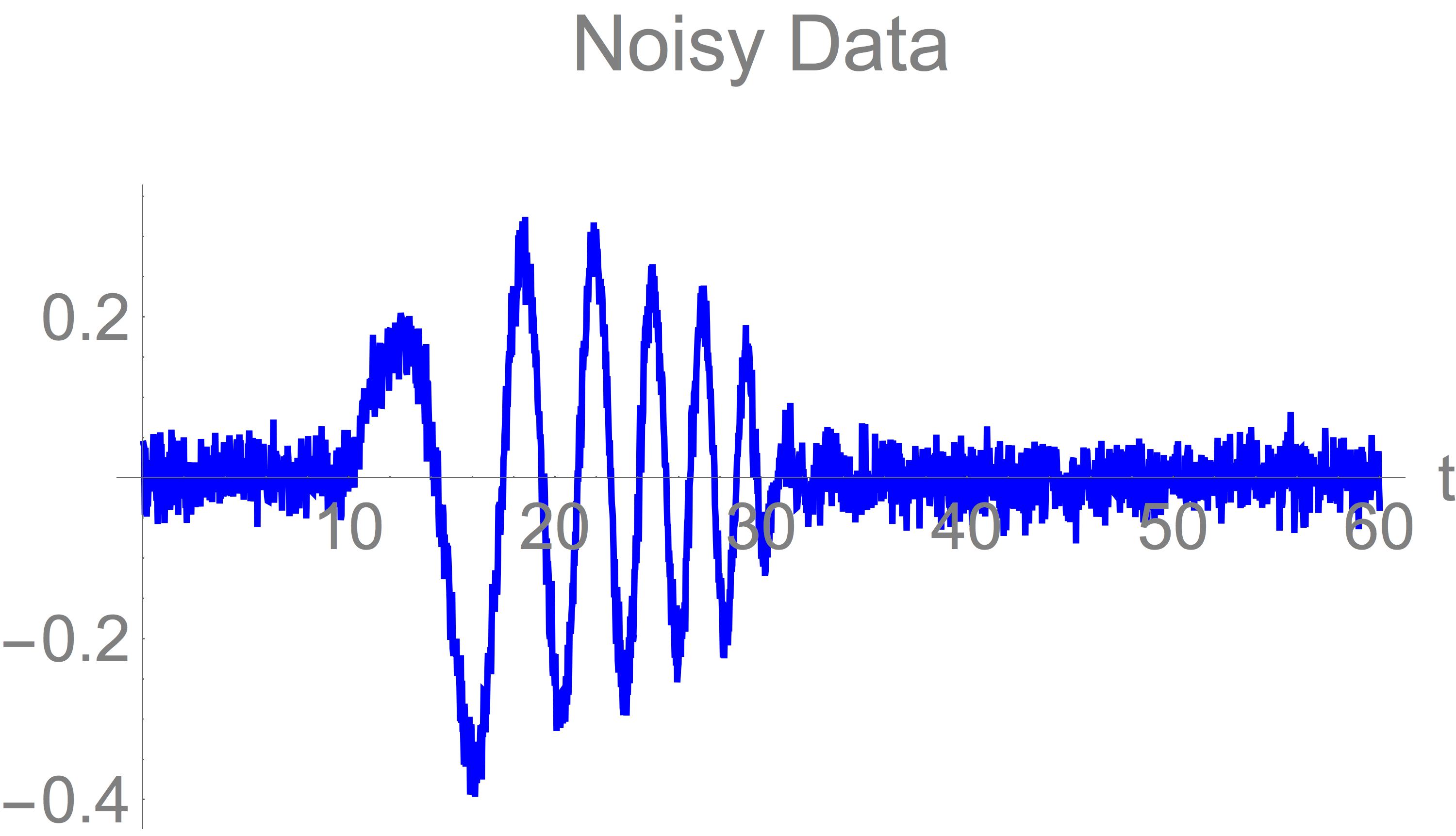}\rule{18pt}{0pt}\includegraphics[width=.4\textwidth]{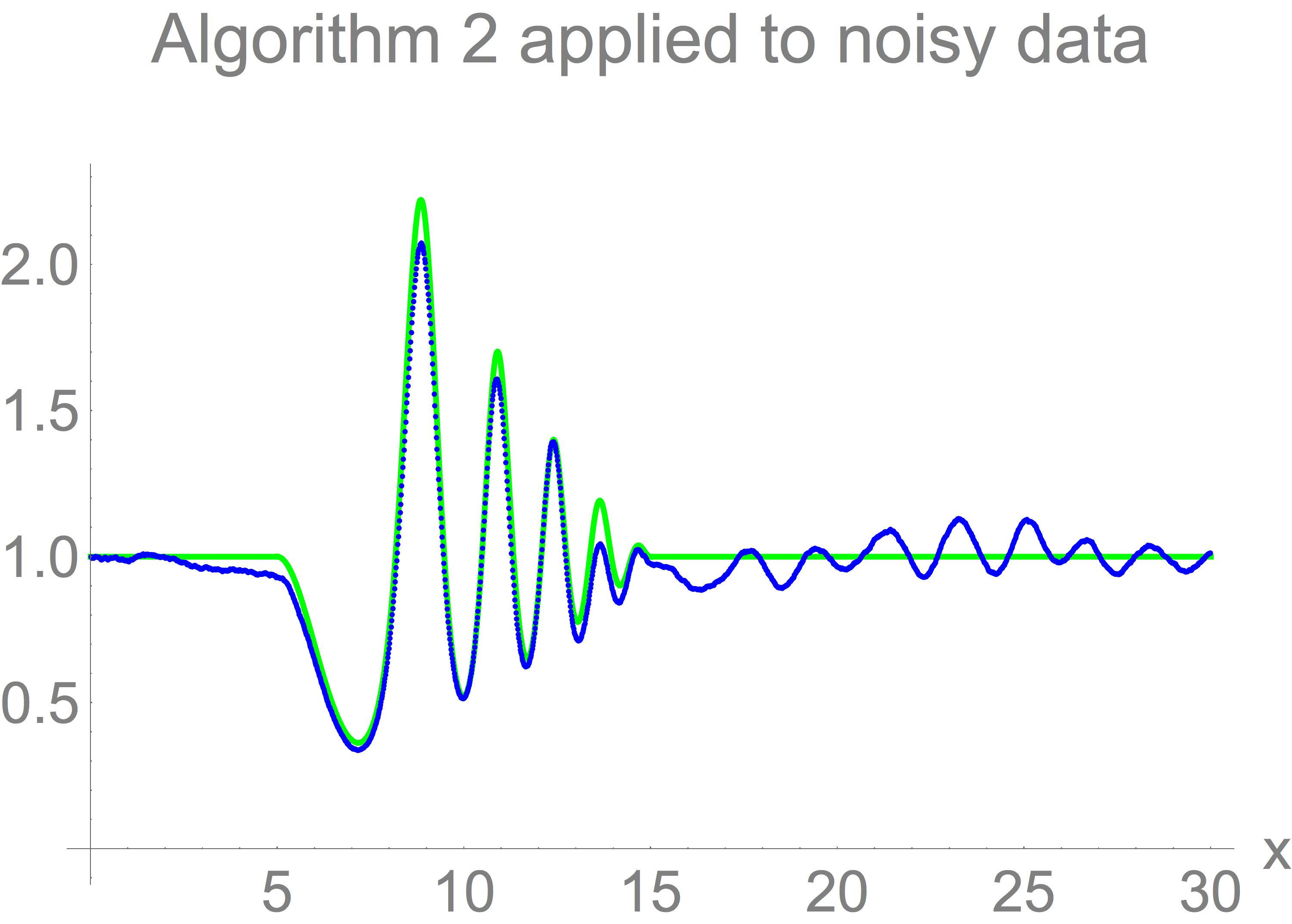}}
\caption{Impedance reconstructed from noisy data by Algorithm~\ref{alg-inverse}. Left: data with $25\%$ iid gaussian noise added. Right: in blue, the output of Algorithm~\ref{alg-inverse} with the noisy data at left as input; true impedance in green for comparison. The $\ell^2$ relative error of the reconstruction is $5.87\%$. 
\label{fig-noisy-reconstruction}}
\end{figure}

A final remark concerns the stability of Algorithm~\ref{alg-inverse}, illustrated by applying it to noisy data. Figure~\ref{fig-noisy-reconstruction} depicts the same data as before with $25\%$ iid gaussian noise added, and the accompanying reconstruction.  More extensive numerical experiments show the roughly $5\%$ error in the given reconstruction, which is far less than that in the noisy data, to be typical behaviour, provided $\zeta_{\max}/\zeta_{\min}$ is not too large.  Thus numerical inversion as encapsulated in Algorithm~\ref{alg-inverse} is not only fast and accurate, it also rather stable.  Large values of $\zeta_{\max}/\zeta_{\min}$ are associated with a maximum value  of $|R(\omega)|$ very close to 1. In this case inversion becomes less stable, much like in the toy example discussed in \S\ref{sec-simple-example}.

\section{Conclusion\label{sec-conclusion}}

This final section remarks on technical and historical aspects of the results obtained above in
sections~\ref{sec-scattering-for-piecewise-continuous} and \ref{sec-fast-algorithms}, and states some open problems. 

\subsection{Remarks\label{sec-conclusion-remarks}}
\subsubsection{The Riemannian structure on $K$\label{sec-conclusion-riemannian}}
The precise way in which the Riemannian structure on the manifold $K$ (representing the closure of PSL$(2,\mathbb{R})$) is manifest in the scattering map is worth reiterating, as it constitutes a fundamentally new insight. Recall the factorization of $\zeta=\zeta_1\zeta_2$ into components $\zeta_1\in C_+(X)$ and $\zeta_2\in \step(X)$. Part~\ref{forward-scattering-asymptotics} of Theorem~\ref{thm-forward-scattering} implies that, as almost periodic functions, the reflection coefficients $g_\zeta(0)$ and $g_{\zeta_2}(0)$ are identical, in the sense that 
\[
\|g_\zeta(0)-g_{\zeta_2}(0)\|_{\rm ap}=0. 
\]
Indeed  $g_{\zeta_2}(0)$ \emph{is}
the almost periodic part of $g_\zeta(0)$ as defined in \S\ref{sec-almost-periodic}.  Moreover, according to Theorem~\ref{thm-fourier-series}, it may be expressed in terms of the eigenfunctions $\pi^{(p,q)}$ of the Laplace-Beltrami operator (\ref{laplace-beltrami}) for $K$ as 
\begin{equation}\label{purely-singular-part}
g_{\zeta_2}^\omega(0)=\sum_{k\in\{1\}\times\mathbb{Z}^{n-1}_+}\mathfrak{a}_k(r,0)e^{2i\langle k,\mathbf{y}\rangle\omega}\quad\mbox{ where }\quad\mathfrak{a}_k(r,0)=\prod_{j=1}^{n}\pi^{(k_j,k_{j+1})}(r_j).
\end{equation}
Here $r_j$ are the reflectivities corresponding to the $n$ jump points $y_j$ of $\zeta_2$, and
\[
\mathbf{y}=(y_1-x_0,y_2-y_1,\ldots,y_n-y_{n-1})\in\mathbb{R}^{n}_+. 
\]
Thus for each $\lambda\in\mathbb{R}$, 
\begin{equation}\label{key-connection}
\langle g_\zeta(0),e^{i\lambda\omega}\rangle_{\rm ap}=\sum_{\stackrel{k\in\{1\}\times\mathbb{Z}^{n-1}_+}{2\langle k,\mathbf{y}\rangle=\lambda}}\prod_{j=1}^{n}\pi^{(k_j,k_{j+1})}(r_j),
\end{equation}
where the right-hand side is understood to  be  zero if $2\langle k,\mathbf{y}\rangle=\lambda$ has no solution $k\in\{1\}\times\mathbb{Z}^{n-1}_+$.  Equation (\ref{key-connection}) shows explicitly how the Riemannian structure of $K$ relates to the scattering map $\zeta\mapsto g_{\zeta}(0)$.

\subsubsection{Computation in the presence of singularities\label{sec-conclusion-computation}}
Section~\ref{sec-fast-algorithms} illustrates in practical terms the implications for digital signal processing of the theoretical results in \S\ref{sec-scattering-for-piecewise-continuous}. In the case of continuous impedance, where $\zeta_2$ is constant, both forward and inverse scattering can be computed efficiently and accurately using simple direct algorithms---with the additional requirement for inverse scattering that $\alpha$ be square integrable.  (Indeed Algorithm~\ref{alg-inverse} is faster than the Durbin algorithm \cite[\S4.7]{GoVa:1996}, details of which will be presented elsewhere.)  Discontinuities in $\zeta$ present additional complications, however, since the inverse Fourier transform of (\ref{purely-singular-part}), 
\begin{equation}\label{purely-singular-time-domain-data}
d_{\rm sing}(t)=\sum_{k\in\{1\}\times\mathbb{Z}^{n-1}_+}\mathfrak{a}_k(r,0)\delta\bigl(t-{2\langle k,\mathbf{y}\rangle}\bigr),
\end{equation}
emerges as a singular component of the time-domain data.
Numerical implementation of the forward map in this case is detailed in \cite{Ch:2021}.  Computation of the inverse map when a non-trivial singular component (\ref{purely-singular-time-domain-data}) is present can be carried out by following a punctuated layer-stripping procedure, roughly as described in \S\ref{sec-inverse-scattering} prior to Theorem~\ref{thm-piecewise-injectivity}, intertwining Algorithm~\ref{alg-inverse} with numerical  inversion of (\ref{purely-singular-time-domain-data}).  This has not been pursued in detail, since real-world measurements sample a mollified version of (\ref{purely-singular-time-domain-data}); and in the case of mollified data 
one can take a more direct approach, as described in \cite{Gi:JCP2018}. 

\subsubsection{Singular approximation\label{sec-conclusion-singular-approximation}}  
In the absence of a continuum limit, and without any attempted analysis of the limiting process, waves in piecewise constant media have been studied by many authors, e.g., \cite{Br:1951,BeGoWa:1958,WaAk:1969,BuBu:1983,FoFr:1990,CaMo:2006}; thus the idea of working with piecewise constant media is far from new.  What is new in the present work is the central technique of singular approximation, whereby a rigorous notion of continuum limit is established and the limiting process analyzed. 

At the core of singular approximation is the harmonic exponential operator as a renormalized continuum limit of OPUC. Proposition~\ref{prop-opuc}\ref{singular-harmonic-OPUC} combines with Lemma~\ref{lem-singular-approximation} to yield the pivotal uniform convergence on compact sets, 
\begin{equation}\label{precise-expression}
\Psi_n^\ast(e^{2i\Delta_n\omega})\rightarrow E_\alpha^{(x_0,y)}(\omega)\quad\mbox{ as }\quad n\rightarrow\infty,
\end{equation}
where $\zeta_n$ is the standard approximation to $\zeta\in C^1_+(X)$, and $\alpha=-\zeta^\prime/(2\zeta)$. 
More generally, continuity of the harmonic exponential with respect to $\alpha\in L^1_\real(X)$ implies $E_\alpha^{(x_0,y)}$ is a continuum limit of OPUC for any absolutely continuous $\zeta\in C_+(X)$, although the approximating polynomials no longer come directly from the standard approximation. 

Holomorphicity of the $m$-function in the upper half plane, as discussed in the introduction, implies the reflection coefficient itself is the boundary function of a function holomorphic in the upper half plane, a fact central to earlier treatments of scattering on the line \cite{Fa:1964,DeTr:1979}.  The limit (\ref{precise-expression}) provides a useful alternative, allowing one to work instead with holomorphicity of polynomials and rational functions on the unit disk, as in the proof of Proposition~\ref{prop-opuc}. The recurrence relation for orthogonal polynomials makes it easy to infer behaviour of $\Psi_n^\ast(z)$ for $z\in\ddd$, and in particular to prove that $\Psi^\ast_n$, $\Phi^\ast_n$ and $\Psi^\ast_n+\Phi^\ast_n$ are holomorphic and zero-free on $\overline{\ddd}$, properties which then carry over to the harmonic exponential and hyperbolic tangent operators.  Singular approximation thus provides an alternative to direct complexification of the spectral parameter $\omega$ in the system (\ref{equation},\ref{bc}). 
 
\subsubsection{Singular approximation as a compositional integral\label{sec-conclusion-compositional-integral}}
On the level of the reflection coefficient, singular approximation maybe viewed as a projective version of multiplicative integration (also known as product integration, or an ordered matrix exponential). Multiplicative integration is typically formulated as a continuum limit of matrix products, such as arises in soliton theory \cite[\S\,I.2,\S\,III.1]{FaTa:2007} (which is loosely related to present considerations) and also survival analysis \cite{GiJo:1990}. Singular approximation has a slight twist in that the continuum limit underlying Theorem~\ref{thm-forward-scattering} (see Proposition~\ref{prop-opuc}\ref{equally-spaced-R-representation}),
\begin{equation}\label{projective-multiplicative-integral}
R_n(\omega)=\frac{\Psi_n^\ast(e^{2i\Delta_n\omega})-\Phi_n^\ast(e^{2i\Delta_n\omega})}{\Psi_n^\ast(e^{2i\Delta_n\omega})+\Phi_n^\ast(e^{2i\Delta_n\omega})}\quad\rightarrow\quad
\Th^{(x_0,y)}_\alpha(\omega)=\frac{E^{(x_0,y)}_\alpha(\omega)-E^{(x_0,y)}_{-\alpha}(\omega)}{E^{(x_0,y)}(\omega)_\alpha+E^{(x_0,y)}_{-\alpha}(\omega)},
\end{equation}
incorporates the elementary projective relation (\ref{multiplicative}) between matrix products and linear fractional transformations, and so is not directly a limit of matrix products.  Properly speaking, singular approximation corresponds to a compositional, rather than a multiplicative, integral.

\subsubsection{Probability measures\label{sec-probability}}
The renormalized limit (\ref{precise-expression}) is also interesting from the perspective of probability measures. 
With $X=(x_0,x_1)$, for each $n\geq 1$, set $$\Delta_n=(x_1-x_0)/(n+1),$$ and let $y_{n,j}$ $(1\leq j\leq n)$ denote the jump points of $\zeta_n$ and $r_{n,j}$ $(1\leq j\leq n)$ the associated reflectivities. Then, according to formula (\ref{unit-circle-measures}) of \S\ref{sec-OPUC}, for each $n$, the sequence of Verblunsky parameters $-r_{n,j}$ corresponds to a probability measure on the circle,
\begin{equation}\label{circle-measure-2}
d\nu_n(e^{i\theta})=\frac{\prod_{j=1}^n(1-|r_{n,j}|^2)}{\bigl|\Psi_n^\ast(e^{i\theta})\bigr|^2}\,\frac{d\theta}{2\pi}=\frac{\prod_{j=1}^n(1-|r_{n,j}|^2)}{\bigl|\Psi_n^\ast(e^{2i\Delta_n\omega})\bigr|^2}\,\frac{\Delta_n\,d\omega}{\pi},
\end{equation}
the latter equation resulting from the change of variable $\theta=2\Delta_n\omega$.  Setting
\[
W_n(\omega)=\frac{\prod_{j=1}^n(1-|r_{n,j}|^2)}{\bigl|\Psi_n^\ast(e^{2i\Delta_n\omega})\bigr|^2}\qquad(n\geq 1),
\]
it follows from (\ref{precise-expression}) and Corollary~\ref{cor-log-rnj-limit}\ref{rnj-square-limit} that 
\begin{equation}\label{Wn-approximation}
W_n\rightarrow\bigl|E_{\alpha}^X\bigr|^{-2}\quad\mbox{ as }\quad n\rightarrow\infty
\end{equation}
uniformly on compact sets.  Turning this around, approximation of $W_n$ by $\bigl|E_{\alpha}^X\bigr|^{-2}$ yields 
\begin{equation}\label{measure-approximation}
d\nu_n(e^{i\theta})\approx\bigl|E_{\alpha}^X\bigl(\theta/(2\Delta_n)\bigr)\bigr|^{-2}\frac{d\theta}{2\pi}.
\end{equation}
Made rigorous, this constitutes a sort of concentration of measure result:  fluctuations in the measure $d\nu_n$ are concentrated nearer and nearer to the point $\theta=0$ as $n\rightarrow\infty$, since for any fixed $\varepsilon>0$, 
\[
E_{\alpha}^X\bigl(\theta/(2\Delta_n)\bigr)\rightarrow 1\;\mbox{ uniformly for $|\theta|>\varepsilon$ as }\; n\rightarrow\infty,
\]
by Corollary~\ref{cor-harmonic-exponential-bounded}.  In the limit, the fluctuations recede into a purely singular object supported at $\theta=0$.  To make the argument rigorous, however, requires a finer approximation result than in \S\ref{sec-approximation}, namely that (\ref{Wn-approximation}) holds uniformly on intervals $\bigl(-\pi/(2\Delta_n),\pi/(2\Delta_n)\bigr)$.  The latter can be proved by  separately analyzing convergence on the the variable windows $\frac{1}{\sqrt{\Delta_n}}\leq|\omega|\leq\frac{\pi}{2\Delta_n}$ and $|\omega|\leq1/\sqrt{\Delta_n}$; however, the details of this analysis are more complicated than the proof of Lemma~\ref{lem-singular-approximation} and have therefore been omitted. 

\subsubsection{Integral equations and inverse scattering\label{sec-integral-equations}}
Beyond the original goal of analyzing scattering on the line for coefficients of low regularity, the results in \S\ref{sec-scattering-for-piecewise-continuous} give a description of forward and inverse scattering in the regular case that is more explicit than classical results---especially notable is the formulation in Theorem~\ref{thm-short-range-inversion} of impedance purely in terms of the reflection coefficient.  Whereas the integral equations of Borg-Marchenko-Gelfand-Levitan, surveyed by Faddeev in \cite{Fa:1963} and by Newton in \cite{Ne:1980}, show that the potential of the Schr\"odinger equation is determined by the reflection coefficient under sufficient regularity, the actual inversion procedure involves more than one equation and does not readily lead to an explicit formula.  The integral equation that arises in the proof of Theorem~\ref{thm-short-range-inversion}, namely (\ref{khf-2}), stems from singular approximation and the probability measures (\ref{measure-approximation}). To be precise, (\ref{khf-2}) is the continuum counterpart of the orthogonality relations 
\[
\langle\Psi_n^\ast(z),z^j\rangle_{d\nu_n(z)}=0\quad(1\leq j\leq n).
\]
Thus the singular case of piecewise constant impedance underlies the short-range inversion formula, even though short-range inversion applies to the absolutely continuous part, i.e.~the regular part, of the impedance function.

\subsection{Open problems\label{sec-conclusion-problems}}
The present paper raises a number of questions that warrant further consideration, some of a technical nature and others more fundamental.  As mentioned earlier, the question of inverse stability will be addressed in separate work. 

A key technical fact concerning piecewise absolutely continuous impedance is that the associated reflection coefficient has modulus bounded away from 1. Proposition~\ref{prop-bounded-variation} requires bounded variation of the impedance function, and it is an open question whether $|g_\zeta^\omega(0)|=1$ is possible for $\zeta\in\reg_+(X)$ not of bounded variation. In other words, the question is whether the constant functions in $K$ are attained by the generalized reflection coefficient $g_\zeta$. 
Constant elements  are a likely obstruction to injectivity of the forward mapping $\zeta\mapsto R=g_\zeta(0)$ on $\reg_+(X)$ with $\zeta(x_0+)$ fixed.  A more fundamental problem is to determine a maximal domain of injectivity, and to characterize the corresponding range. 

A technical problem relating to absolutely continuous impedance functions was mentioned already in \S\ref{sec-classical-trace}:  to extend the trace formula (\ref{classical-trace}) to $\alpha\in L^2_\real(X)$.  More generally, what is the maximal domain of validity of the singular trace formula of Theorem~\ref{thm-renormalized-trace}?  Does it extend to arbitrary discontinuous $\zeta\in\reg_+(X)$? 

As discussed in the introduction, inverse scattering for the Helmholtz equation carries over to the Schr\"odinger equation provided the potential belongs to $H^{-1}$, allowing singularities of Dirac type, for example. On a formal level, the results go further, since discontinuities in $\zeta$ correspond to singularities in the potential that are derivatives of Dirac functions. But in this case the Schr\"odinger equation no longer has an interpretation in the sense of distributions. It is an open question whether scattering for the Schr\"odinger equation can be extended to a sufficiently general setting, such as that of Colombeau algebras \cite{Co:1990,Ob:1992}, so that results in the present paper for $\zeta\in\pwacplus(X)$ may be carried over directly.


\end{document}